\documentclass[11pt]{article}
\usepackage[american]{babel}
\usepackage{microtype}
\usepackage{graphicx}
\usepackage{float}
\usepackage{epstopdf}
\usepackage{pdfsync}
\usepackage{amsmath,amssymb,amsthm,amsfonts,mathrsfs}
\usepackage{amssymb}
\usepackage{verbatim}
\usepackage{setspace}
\usepackage{hyperref}
\usepackage[a4paper]{geometry}
\newtheorem{thm}{Theorem}[section]
\newtheorem{cor}[thm]{Corollary}
\newtheorem{lem}[thm]{Lemma}
\newtheorem{prop}[thm]{Proposition}

\theoremstyle{definition}
\newtheorem{defn}{Definition}[section]
\theoremstyle{remark}
\newtheorem{rem}{Remark}[section]

\numberwithin{equation}{section}

\DeclareMathSymbol{\C}{\mathalpha}{AMSb}{"43}

\newcommand{\eps}{\varepsilon}

\newcommand{\alp}{\alpha}
\newcommand{\x}{{\mathcal{X}}}

\newcommand{\dx}{\,\mathrm{d}x}

\newcommand{\R}{{\mathbb{R}}}
\newcommand{\h}{{\mathcal{H}}}
\newcommand{\m}{{\mathcal{M}}}
\newcommand{\inte}{\int_{\mathbb{R}^2}}
\newcommand{\intB}{\int _{B_\delta (x_{2,k})}}
\newcommand{\intPB}{\int _{\partial B_\delta (x_{2,k})}}

\newcommand{\bsub}{\begin{subequations}}
\newcommand{\esub}{\end{subequations}$\!$}

\title{Ground States of Two-Component Attractive Bose-Einstein Condensates I: Existence and Uniqueness}
\author{
Yujin Guo\thanks{Wuhan Institute of Physics and Mathematics,
    Chinese Academy of Sciences, P.O. Box 71010, Wuhan 430071,
    P. R. China.  Email: \texttt{yjguo@wipm.ac.cn}. Y. J.  Guo is partially supported by NSFC under Grant No. 11671394 and MOST under Grant No. 2017YFA0304500.
    },
    \, Shuai Li\thanks{University of Chinese Academy of Sciences, Beijing 100190, P. R. China;  Wuhan Institute of Physics and Mathematics,
    Chinese Academy of Sciences, P.O. Box 71010, Wuhan 430071,
    P. R. China.  Email: \texttt{lishuai\_wipm@outlook.com}.  },
\,   Juncheng Wei\thanks{Department of Mathematics, University of British Columbia, Vancouver, BC V6T 1Z2,
    Canada.  Email: \texttt{jcwei@math.ubc.ca}. J. C. Wei is partially supported by NSERC of Canada.}
    \, and\, Xiaoyu Zeng\thanks{Department of Mathematics, Wuhan University of Technology, Wuhan 430070, P. R. China  Email: \texttt{xyzeng@whut.edu.cn}. X. Y. Zeng is partially supported by NSFC grant 11501555.}
    }
\date{\today}

\begin{document}

\maketitle
\begin {abstract}
We study ground states of two-component Bose-Einstein condensates (BEC) with trapping potentials in $\R^2$, where the intraspecies interaction $(-a_1,-a_2)$ and the interspecies interaction $-\beta$ are both attractive, $i.e,$ $a_1$, $a_2$ and $\beta $ are all positive. The existence and non-existence of ground states are classified completely by investigating equivalently the associated $L^2$-critical constraint variational problem. The uniqueness and symmetry-breaking of ground states are also analyzed under different types of trapping potentials as $\beta \nearrow \beta ^*=a^*+\sqrt{(a^*-a_1)(a^*-a_2)}$, where $0<a_i<a^*:=\|w\|^2_2$ ($i=1,2$) is fixed and $w$ is the unique positive solution of $\Delta w-w+w^3=0$ in $\R^2$. The semi-trivial limit behavior of ground states is tackled in the companion paper \cite{GLWZ}.
\end {abstract}

\vskip 0.2truein

\tableofcontents

\section{Introduction}
In this paper, we consider the following coupled nonlinear Gross-Pitaevskii equations
\begin{equation}\label{equ:CGPS}
\begin{cases}
-\Delta u_{1} +V_1(x)u_{1} =\mu u_{1} +a_1u_{1}^3 +\beta u_{2}^2 u_{1}   \,\ \mbox{in}\,\  \R^2,\\
-\Delta u_{2} +V_2(x)u_{2} =\mu u_{2} +a_2u_{2}^3 +\beta u_{1}^2 u_{2}   \,\ \mbox{in}\,\  \R^2,\,\
\end{cases}
\end{equation}
where $(u_{1} ,u_{2})\in\mathcal{X}=\h_1(\R ^2)\times \h_2(\R ^2)$  and
the space 
\begin{equation*}
 \h_i(\R ^2) = \Big \{u\in  H^1(\R ^2):\ \int _{\R ^2}  V_i(x)|u(x)|^2\dx<\infty \Big\},\,\
\end{equation*}
is equipped with the norm $\|u\|_{_{\h_i}}=\big(\int _{\R ^2} \big[|\nabla u|^2+ V_i(x)|u(x)|^2\big] \dx\big)^{\frac{1}{2}}$ for $i=1, 2$. The system \eqref{equ:CGPS} arises in describing two-component Bose-Einstein condensates (BEC) with trapping potentials $V_1(x) $ and $V_2(x) $ (cf. \cite{EGBB,HMEWC,LWCMP,LW,LW2,PW,Royo}), and $\mu\in\R$ is a chemical potential.  From the physical point of view, we assume
that the trapping potential $0\le V_i(x)\in C^\alpha_{\rm loc}(\R^2)$ (where $0<\alpha<1$) satisfies for $i=1,2,$
\begin{equation}\label{cond:V.1}
\lim_{|x|\to\infty} V_i(x) = \infty,\,\   \text{both}\,\ \inf_{x\in\R^2}V_i(x)=0
\,\ \text{and}\ \inf\limits_{x\in \R^2} \big(V_1(x)+V_2(x)\big)\,\ \text{are attained.}
 \end{equation}
Here $a_i>0$ ($resp.$ $<0$) represents that the intraspecies interaction  of the atoms inside each component is attractive  ($resp.$ repulsive), and $\beta>0$ ($resp.$ $<0$) denotes that the interspecies interaction  between two components is attractive ($resp.$ repulsive).


The main aim of the present paper is to analyze  ground states of the system (\ref{equ:CGPS}) for the case where the intraspecies interaction and the interspecies interaction are both attractive, $i.e,$ $a_1>0$, $a_2>0$ and $\beta >0$.  As illustrated by Proposition A.1 in the Appendix A.1, ground states of the system (\ref{equ:CGPS}) in this case can be described {\em equivalently} by nonnegative minimizers of
the following $L^2-$critical constraint variational problem
\begin{equation}\label{def:e}
e(a_1,a_2,\beta):=\inf_{(u_1,u_2)\in \mathcal{M}} E_{a_1,a_2,\beta}(u_1,u_2),\,\ a_1>0,\,\ a_2>0,\,\ \beta>0,
\end{equation}
where
\begin{equation*}
\mathcal{M}:=\Big\{(u_1,u_2)\in \mathcal{X}:\, \int_{\R ^2}(|u_1|^2+|u_2|^2 )\dx=1\Big\},
\end{equation*}
and the Gross-Pitaevskii (GP)  energy functional $ E_{a_1,a_2,\beta}(u_1,u_2)$ is given by
\begin{equation}\label{def:E}
\begin{split}
E_{a_1,a_2,\beta}(u_1,u_2)
=&\int_{\R ^2} \big(|\nabla u_1|^2+|\nabla u_2|^2\big)\dx +\int_{\R ^2}\big( V_1(x)|u_1|^2+V_2(x)|u_2|^2\big) \dx\\
&-\int_{\R ^2} \Big(\frac{a_1}{2}|u_1|^4+\frac{a_2}{2}|u_2|^4 +\beta |u_1|^2|u_2|^2\Big) \dx \,,\quad (u_1,u_2)\in\mathcal{X}.
\end{split}
\end{equation}
To discuss equivalently ground states of (\ref{equ:CGPS}), throughout the whole paper we shall therefore focus on investigating  \eqref{def:e}, instead of (\ref{equ:CGPS}).

As another type of motivations, the studies of \eqref{def:e} are also stimulated by the recent works \cite{BC,GS,GWZZ,GZZ} and the references therein, where the authors investigated ground states of a single component attractive BEC, $i.e.$ the following single minimization problem
\begin{equation}\label{1:two}
e_i(a):=\inf_{\{u\in\h_i,\int_{\mathbb{R}^2} |u|^2dx=1\}} E_{a}^i(u),\,\ a>0,
\end{equation}
and the GP energy functional $E_{a}^i(u)$ satisfies
\begin{equation*}
E_{a}^i(u):=\int_{\R ^2} \big(|\nabla
  u(x)|^2+V_i(x)|u(x)|^2\big)dx-\frac{a}{2}\int_{\R ^2}|u(x)|^4dx, \,\  i=1 \,\  \mbox{or}\,\  2.
\end{equation*}
Actually, it was proved in \cite{BC,GS} that (\ref{1:two}) admits minimizers if and only if $0<a<a^*:=\|w\|^2_2$, where $w=w(|x|)>0$ denotes  (cf.  \cite{GNN,Kwong,LN}) the unique positive solution of the following nonlinear scalar field
equation
\begin{equation}\label{equ:w}
\Delta w-w+w^3=0,\,\,\  w \in H^1(\R^2).
\end{equation}
The blow-up behavior of minimizers of $e_i(a)$ as $a\nearrow a^*$ was also analyzed in \cite{GS,GWZZ,GZZ} under different types of trapping potentials. The above mentioned results  show that the analysis of $e_i(a)$ makes full use of the following classical Gagliardo-Nirenberg inequality
\begin{equation}\label{ineq:GNQ}
\frac {\|w\|_2^{2}}{2}
  =\inf\limits_{\{u(x)\in H^1(\R^2)\setminus \{0\}\}}
\frac{\inte |\nabla u(x) |^2 \dx \inte |u(x)|^2 \dx}{\inte |u(x)|^4  \dx},
\end{equation}
where the equality is attained at $w$ (cf. \cite{W}). Applying \eqref{equ:w} and \eqref{ineq:GNQ}, one can get the following identifies
\begin{equation}\label{ide:w}
  \|w\|_2^2=\|\nabla w\|_2^2=\frac{1}{2} \|w\|_4^4.
\end{equation}
Note also from \cite[Proposition 4.1]{GNN} that $w(x)$ decays exponentially in the sense that
 \begin{equation} \label{decay:w}
w(x) \, , \ |\nabla w(x)| = O(|x|^{-\frac{1}{2}}e^{-|x|}) \,\
\text{as} \,\ |x|\to \infty.
\end{equation}
Since the interspecies interaction
between the components leads to more elaborate physical phenomena, one may expect that multiple-component BECs present more complicated characters than one-component BEC, and the corresponding analytic investigations are more challenging. Therefore, comparing with the analysis of (\ref{1:two}), investigating $e(a_1,a_2,\beta)$ needs more involved analytic methods and variational arguments to overcome above difficulties.

When the constraint condition of \eqref{def:e} is replaced by
$\big\{(u_1,u_2)\in\mathcal{X}:\, \int_{\R ^2}|u_1|^2 \dx=\int_{\R ^2}|u_2|^2 \dx=1\big\}$,
the existence, nonexistence, mass concentration, and other analytic properties of minimizers for \eqref{def:e} are analyzed recently in \cite{GZZ2}. Different from \cite{GZZ2}, in the present paper and  \cite{GLWZ} we shall analyze the  minimizers of \eqref{def:e} under a different constraint assumption that $(u_1,u_2)\in\mathcal{M}$.
In this paper, we shall not only address a {\em complete} classification of the existence and nonexistence of minimizers for \eqref{def:e}, but also discuss  more importantly the uniqueness,  mass concentration, and symmetry-breaking of minimizers for \eqref{def:e}. In the companion paper \cite{GLWZ}, we shall analyze the semi-trivial limit behavior of minimizers for \eqref{def:e}.

\subsection{Main results}

The first main result of the present paper is concerned with the following existence and nonexistence of minimizers:

\begin{thm}\label{Th:existence1}
Suppose $V_i(x)$ satisfies \eqref{cond:V.1}  for $i=1,\,2$, and set
\begin{equation}\label{thm1.1:1}
\beta^*=\beta^*(a_1,a_2):= a^* + \sqrt {(a^*-a_1)(a^*-a_2)}, \,\ \text{where}\,\ 0< a_1, a_2<a^*:=\|w\|_2^2.
\end{equation}
Then we have
\begin{enumerate}
\item   If $0<a_1,a_2<a^*$ and  $0<\beta<\beta^*$, then problem (\ref{def:e}) has at least one minimizer.
\item   If either $a_1> a^*$ or $a_2>a^*$ or $\beta> \beta^*$, then problem (\ref{def:e}) has no minimizer.
\end{enumerate}
\end{thm}

We remark that Theorem \ref{Th:existence1} was already proved in \cite[Theorems 2.1 and 2.2]{BC}. In this paper we shall address the proof of Theorem \ref{Th:existence1} in a simpler way, which depends strongly on the following Gagliardo-Nirenberg type inequality
 \begin{equation}\label{Ineq:GN}
 \int_{\R ^2} \big(|u_1|^2+|u_2|^2  \big)^2 \dx
 \le  \frac{2}{ \|w\|_2^2}
   \int_{\R ^2} \big(|\nabla u_1|^2+|\nabla u_2|^2\big) \dx  \int_{\R ^2}\big(|u_1|^2+|u_2|^2\big) \dx,
 \end{equation}
where $(u_1,u_2)\in H^1(\R^2)\times H^1(\R^2)$. The proof of Lemma \ref{Lem:w1w2w} in the appendix shows that $\frac{2}{ \|w\|_2^2} $ is the best constant of (\ref{Ineq:GN}), where the equality is attained at $(w\sin \theta  , w\cos \theta )$ for any $\theta \in [0,2\pi)$. The following theorem is concerned with the existence of  minimizers for \eqref{def:e} at the threshold where either $a_1=a^*$ or $a_2=a^*$ or $\beta=\beta^*$:

\begin{thm}\label{Th:existence2}
Suppose $V_i(x)$ satisfies \eqref{cond:V.1}  for $i=1,\,2$, and assume $\beta^*=\beta^*(a_1,a_2)\ge a^*$ is defined by (\ref{thm1.1:1}).
Then we have the following results:
\begin{enumerate}
\item   Suppose either $a_1=a^*$ and $\beta\leq\beta^*(=a^*)$, or $a_2=a^*$ and $ \beta\leq\beta^*(=a^*)$. Then there is no minimizer for problem (\ref{def:e}).
Furthermore, we have $\lim\limits_{a_1\nearrow a^*} e(a_1,a_2,\beta)=e(a^*,a_2,\beta)=0$
and $\lim\limits_{a_2\nearrow a^*} e(a_1,a_2,\beta)=e(a_1,a^*,\beta)=0$, where $\beta\leq a^*$.

\item  Suppose  $\beta= \beta^*$, $0< a_1,a_2<a^*$, and   $V_1(x)$ and $V_2(x)$ have at least one common minimum  point. Then there is no minimizer for problem (\ref{def:e}).
Furthermore,  \begin{equation}\label{eq1.15}
\lim\limits_{\beta\nearrow\beta^*} e(a_1,a_2,\beta)=e(a_1,a_2,\beta^*)=0.\end{equation}

\item
Suppose  $\beta= \beta^*$, $a_1<a^*$ and $a_2<a^*$, and assume  $V_1(x)$ and $V_2(x)$ have no common minimum  points.
If
\begin{equation}\label{thm1.2:1A}
e(a_1,a_2,\beta)<\frac{1}{\sqrt{a^*-a_1}+\sqrt{a^*-a_2}}
\inf\limits_{x\in\R^2}\big(\sqrt{a^*-a_2}V_1(x)+\sqrt{a^*-a_1}V_2(x)\big),
\end{equation}
then there exists at least one minimizer for problem (\ref{def:e}).
\end{enumerate}
\end{thm}

\begin{figure}
\centering
{\includegraphics[width = 10cm,height=6cm,clip]{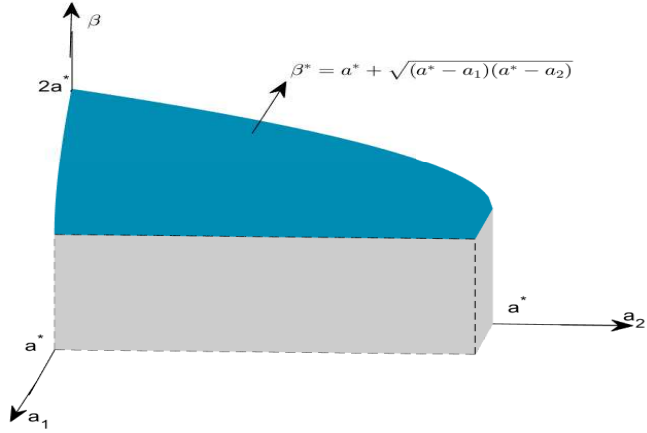}}
    \caption{\em  If $V_1(x)$ and $V_2(x)$ have at least one common minimum  point, then $e(a_1,a_2,\beta)$ has minimizers  if and only if the point $(a_1,a_2,\beta)$ lies within the cuboid.}
 \end{figure}

If $V_1(x)$ and $V_2(x)$ have at least one common minimum  point,
Theorems \ref{Th:existence1} and \ref{Th:existence2} give {\em a complete classification} of the existence and nonexistence of minimizers for $e(a_1,a_2,\beta)$. More precisely, (\ref{def:e}) has minimizers in this case, if and only if the point $(a_1,a_2,\beta)$ lies within the cuboid as illustrated by  Figure 1. On the other hand,
when $V_1(x)$ and $V_2(x)$ have no common minimum  point, we expect that both the existence and nonexistence of minimizers may occur at the threshold where $\beta= \beta^*$, $a_1<a^*$ and $a_2<a^*$, which depends on the shapes of both $V_1(x)$ and $V_2(x)$. Under the additional assumption (\ref{thm1.2:1A}), Theorem \ref{Th:existence2}(3) is proved by applying Ekeland's variational principle (cf. \cite[Theorem 5.1]{S}).

Inspired by \cite{GS,GWZZ,GZZ,GZZ2}, in the following we focus on analyzing the limit behavior of minimizers as $\beta\nearrow \beta^*$.
Since $|\nabla |u||\leq |\nabla u|$ holds $a.e.$ in $\R^2$, any minimizer $(u_1, u_2)$ of \eqref{def:e} satisfies either  $u_i\ge 0$  or $u_i\le 0$ in $\R^2$ for $i=1 $ and $ 2$. Without loss of generality, we therefore remark that  minimizers  of \eqref{def:e} can be restricted to nonnegative  vector functions.
For convenience, we next rewrite the functional  $E_{a_1,a_2,\beta}(\cdot)$ as
\begin{equation}\label{exp2:ea}
\begin{split}
  &E_{a_1,a_2,\beta}(u_1,u_2) \\
  = &\int_{\R^2} |\nabla u_1(x)|^2 + |\nabla u_2(x)|^2\dx - \frac{a^*}{2}\int_{\R^2}\big(|u_1(x)|^2+|u_2(x)|^2\big)^2\dx \\
  & + \int_{\R^2} V_1(x)|u_1(x)|^2\dx + \int_{\R^2} V_2(x)|u_2(x)|^2\dx \\
  & + \frac{1}{2}\int_{\R^2}\big(\sqrt{a^*-a_1}|u_1(x)|^2-\sqrt{a^*-a_2}|u_2(x)|^2\big)^2\dx \\
  & + (\beta^* - \beta) \int_{\R^2} |u_1(x)|^2|u_2(x)|^2\dx,
\end{split}
\end{equation}
and define the positive constant $\gamma$ as
\begin{equation}\label{def:beta.k}
\gamma:=\gamma(a_1, a_2)=\frac{\sqrt{a^*-a_2}}{\sqrt{a^*-a_1}+\sqrt{a^*-a_2}}\in (0,1), \text{ where }\,0< a_1,a_2<a^*.
\end{equation}
 Through analyzing (\ref{exp2:ea}), we shall establish the following limit behavior of nonnegative  minimizers of \eqref{def:e} as $\beta\nearrow \beta^*$.

\begin{thm}\label{Thm1.3}
Suppose that $V_1(x)$ and $V_2(x)$ satisfy \eqref{cond:V.1}  and have at least one common minimum  point.
Let $(u_{1k}(x),u_{2k}(x))$ be a nonnegative minimizer of $e(a_1,a_2,\beta_k)$, where $0< a_1, a_2<a^*$ and $\beta_k\nearrow\beta^*$ as $k\to\infty$. Then there exists a subsequence of $\{\beta _k\}$, still denoted by
$\{\beta _k\}$, such that  $(u_{1k},u_{2k})$ satisfies
\begin{equation}\label{lim:beta.u.exp}
\begin{cases}
 \lim\limits_{k\to\infty}\varepsilon_{k} u_{1{k}}(\varepsilon_{k} x+\bar{z}_{1k})
 =\frac{\sqrt{\gamma}}{\|w\|_2}w(x),\\
 \lim\limits_{k\to\infty}\varepsilon_{k} u_{2{k}}(\varepsilon_{k} x+\bar{z}_{2k})
 =\frac{\sqrt{1-\gamma}}{\|w\|_2}w(x),
\end{cases} \text{ strongly in $H^1(\R^2)$},
\end{equation}
where $0<\gamma =\gamma(a_1, a_2)<1$ is given by (\ref{def:beta.k}), $\varepsilon_{k}>0$ satisfies
\begin{equation}\label{thm1.3:1}
\varepsilon_{k}^{-2}=\int_{\R^2}\big(|\nabla u_{1{k}}(x)|^2+|\nabla u_{2{k}}(x)|^2\big)\dx\to+\infty\,\ \text{as}\,\ k\to\infty,
\end{equation}
and $\bar{z}_{ik}$ is the unique maximum point of
 $u_{ik}$ satisfying
\begin{equation}\label{lim:beta.z}
  \lim\limits_{k\to\infty}\bar{z}_{1k}=\lim\limits_{k\to\infty}\bar{z}_{2k}=\bar{x}_0,
 \,\, \text{where}\,\, \bar{x}_0\in\R^2\,\, \text{satisfies}\,\, V_1(\bar{x}_0)=V_2(\bar{x}_0)=0,
\end{equation}
and
\begin{equation}\label{eq1.22}
\lim\limits_{k\to\infty}\frac{|\bar{z}_{1k}-\bar{z}_{2k}|}{\varepsilon_{k}}=0.
\end{equation}
\end{thm}

Theorem \ref{Thm1.3} shows that the minimizers of $e(a_1,a_2,\beta_k)$ blow up and concentrate at a common minimum point of $V_1(x)$ and $V_2(x)$  as $\beta_k \nearrow\beta^*$.  Unfortunately,
one can note from (\ref{thm1.3:1}) that Theorem \ref{Thm1.3} cannot give the explicit blow-up information  of $(u_{1k}, u_{2k})$ as $k\to\infty$, due to the fact that the specific local profiles near the common minimum  points of potentials $V_1(x)$ and $V_2(x)$ are unavailable in Theorem \ref{Thm1.3}. For this reason, in the following we consider some typical classes of special potentials $V_1(x)$ and $V_2(x)$ to get the precise blow-up behavior of minimizers. We thus define
\begin{defn} \label{defn:1.1}
A function $ f(x)$   is homogeneous of degree $p\in\R^+$ (about the origin), if there holds that
 \begin{equation*}
f(tx)=t^pf(x)\ \, \mbox{in}\ \, \R^2 \ \mbox{for any}\ t>0.
\end{equation*}
\end{defn}
\noindent If $f(x)\in C(\R^2)$ is homogeneous of degree $p>0$,  the above definition then implies that
 \begin{equation*}
0\le f(x)\le C|x|^{p}\,\ \mbox{in}\,\ \R^2,
\end{equation*}
where $C>0$ denotes the maximum of $f(x)$ on $\partial B_1(0)$. Moreover, if $\lim_{|x|\to\infty} f(x) = \infty$, then $0$ is the unique minimum point of $f(x)$.

For generality, in what follows we assume that  $V_1(x)$ and $V_2(x)$ have exactly $l$ common minimum points, namely,
\begin{equation}\label{def:beta.z}
  Z:=\big\{x\in\R^2:V_1(x)=V_2(x)=0\big\}=\big\{x_{1}, x_{2}, \cdots , x_{l}\big\}, \,\ \text{where}\,\ l\geq1.
\end{equation}
We also assume that for $i=1, 2$ and $j=1,2,\cdots,l$, $V_i(x)$ is almost  homogeneous of degree $p_{ij}>0$ around each $x_j$  in the sense that 
\begin{equation}\label{eq1.26}
\lim_{x\to0}\frac{V_i(x+x_j)}{ V_{ij}(x)}=1,
\end{equation}
where $V_{ij}(x)$ satisfies
\begin{equation}\label{eq1.27}
\text{$V_{ij}(x)\in C^\alpha_{\rm loc}(\R^2)$ is homogeneous of degree $p_{ij}>0$}\, \text{ and }\, \lim_{|x|\to\infty} V_{ij}(x) = +\infty.
\end{equation}
Additionally, we also define $H_j(y) $ $(1\leq j\leq l)$  by
\begin{equation}\label{def:unique.Hy}
H_j(y)=\begin{cases}\gamma\displaystyle
\inte V_{1j}(x+y)w^2(x)dx\,\  &\text{if }\,\  p_{1j}<p_{2j},\\[3mm]
\displaystyle\inte \big[\gamma V_{1j}(x+y)+(1-\gamma)V_{2j}(x+y)\big]w^2(x)dx\,\  &\text{if }\,\  p_{1j}=p_{2j},\\[3mm]
(1-\gamma)\displaystyle\inte V_{2j}(x+y)w^2dx\,\  &\text{if }\,\  p_{1j}>p_{2j},
\end{cases}
\end{equation}
where $0<\gamma <1$ is given by \eqref{def:beta.k}.
Define
\begin{equation}\label{def:beta.p0}
  {p}_0:=\max_{1\leq j\leq l}{p}_j \text{ with }{p}_j:=\min\big\{{p}_{1j},{p}_{2j}\big\}, \text{ and }\bar{Z}:=\big\{x_{j}\in Z: p_j= p_0\big\} \subset Z.
\end{equation}
Set
\begin{equation}\label{def:beta.gamma}
\bar{\lambda}_0:=\min\limits_{ j\in \Gamma }\bar\lambda_j, \,\text{ where }\bar\lambda_j:=\min_{y\in \R^2}H_j(y)\ \text{ and }\,\Gamma:=\big\{j: x_{j}\in \bar Z\big\}.
\end{equation}
Denote
\begin{equation}\label{def:beta.z0}
Z_0:=\big\{x_{j}\in \bar{Z}: \bar\lambda_j= \bar{\lambda}_0\big\}
\end{equation}
 the set of the flattest common minimum points of $V_1(x)$ and $V_2(x)$.
Under above assumptions, we  have the following blow-up behavior of nonnegative minimizers as $\beta\nearrow\beta^*$.

\begin{thm}\label{Thm1.4}
Suppose that $V_1(x)$ and $V_2(x)$ satisfy  (\ref{def:beta.z})--(\ref{eq1.27}).
Let $(u_{1k}(x),u_{2k}(x))$ be the  convergent minimizer subsequence of $e(a_1,a_2,\beta_k)$ obtained in Theorem \ref{Thm1.3}, where $0< a_1, a_2<a^*$ and $\beta_k\nearrow\beta^*$ as $k\to\infty$.  Then we have
\begin{equation}\label{lim:beta.V.u.exp}
\begin{cases}
 \lim\limits_{k\to\infty} \bar{\varepsilon}_{k} u_{1{k}}( \bar{\varepsilon}_{k} x+\bar{z}_{1k})
=\frac{\sqrt{\gamma}}{\|w\|_2}w(x),\\
 \lim\limits_{k\to\infty} \bar{\varepsilon}_{k} u_{2{k}}( \bar{\varepsilon}_{k} x+\bar{z}_{2k})
=\frac{\sqrt{1-\gamma}}{\|w\|_2}w(x)
\end{cases}\text{ strongly in } H^1(\R^2),
\end{equation}
where $\bar \varepsilon_{k}>0$ is given by
\begin{equation}\label{def:beta.V.eps}
  \bar{\varepsilon}_{k}:=\Big[\frac{4\gamma(1-\gamma)}{\bar{p}_0\bar{\lambda}_0}(\beta^*-\beta_k)\Big]^\frac{1}{\bar{p}_0+2},
\end{equation}
and $0<\gamma =\gamma(a_1, a_2)<1$ is given by (\ref{def:beta.k}).
Moreover, up to a subsequence if necessary, we have
\begin{equation}\label{lim:beta.V.y0}
  \lim\limits_{k\to\infty}\frac{\bar{z}_{1k}- x_{j_0}}{\bar{\varepsilon}_{k}}
  =\lim\limits_{k\to\infty}\frac{\bar{z}_{2k}-x_{j_0}}{\bar{\varepsilon}_{k}}
  =y_0,
\end{equation}
where $x_{j_0}\in Z_0$ and $y_0\in \R^2$ satisfies  $H_{j_0}(y_0)=\min_{y\in\R^2}H_{j_0}(y)=\bar\lambda_0$.
\end{thm}

\begin{rem}
Let $V_1(x)$ and $V_2(x)$ be given by
\begin{equation*}
  V_i(x)=g_i(x)\prod^{n_i}_{j=1}|x-x_{ij}|^{{p}_{ij}},\,\, {p}_{ij}\in\R^+\,\,\text{ and                                                   }\,\,n_i\in \mathbb{N}^+, \,\, i=1,\,2,
\end{equation*}
where  $g_i(x)\in C_{loc}^\alpha(\R^2)$ satisfies $\frac{1}{C}<g_i(x)<C$ in $\R^2$ for some positive constant $C$.
We also assume that there exists a constant $l\in \mathbb{N}^+$ satisfying $1\leq l\leq\min\{n_1,n_2\}$ such that
\begin{equation*}
\begin{cases}
x_{1j}=x_{2j},\,\ \text{where} \,\  j=1, \cdots , l, \\
x_{1j}\neq x_{2k}, \,\ \text{where} \,\  j\in\{l+1,\cdots ,n_1\} \, \text{ and }\, k\in\{l+1,\cdots ,n_2\}.
\end{cases}
\end{equation*}
Then, $V_1(x)$ and $V_2(x)$ satisfy all assumptions of  Theorem \ref{Thm1.4}. In this setting, we have
$$V_{ij}(x)=g_i(x_{ij})|x|^{p_{ij}}, \text{ where }i=1,2 \text{ and }j=1,2,\cdots, l.$$
Moreover, one can check that each $H_j(y)$, where $j=1, \cdots , l$, attains its minimum at the unique point $y=0$. Thus, we have $y_0=0$ in \eqref{lim:beta.V.y0}.
\end{rem}
Theorem \ref{Thm1.4} shows that  minimizers of $e(a_1,a_2,\beta)$ must concentrate at a flattest common minimum point  of potentials $V_1(x)$ and $V_2(x)$. Moreover, if the trapping potentials $ V_1(x)=V_2(x)$  have a symmetry, e.g., $V_1(x)=V_2(x)=\prod
_{i=1}^n|x-x_i|^p$ with $p>0$ and $x_i$ arranged on the vertices
of a regular polygon, Theorem \ref{Thm1.4} then implies that the following  {\em symmetry
breaking and multiplicity} occur for the minimizers of $e(a_1,a_2,\beta)$ as $\beta \nearrow \beta^*$: for any fixed $0< a_1, a_2<a^*$,  there exists $\beta_*<\beta ^*$
such that for $\beta_* < \beta < \beta^*$, $e(a_1,a_2,\beta)$ has
at least $n$ different nonnegative minimizers, each of which
concentrates at a specific common minimum point $x_i$.

The uniqueness of ground states for single component BEC was addressed recently in \cite[Theorem 1.1]{GLW} by deriving Pohozaev identities, see also \cite{Cao,Deng,Grossi} for related works. Motivated by this fact, we finally investigate the uniqueness of nonnegative minimizers for $e(a_1,a_2,\beta )$ under some further assumptions on $V_1(x)$ and $V_2(x)$, where  $0< a_1, a_2<a^*$ are fixed and $\beta \nearrow \beta^*$.
We shall assume that $V_1(x)$ and $V_2(x)$ have a unique flattest common minimum point, i.e., $Z_0$ defined in (\ref{def:beta.z0}) contains only one element. Our uniqueness can be then stated as follows.


\begin{thm}\label{Thm1.5} Suppose that $V_i(x)\in C^2(\R^2)$ satisfies  (\ref{def:beta.z})--(\ref{eq1.27}) for $i=1,\,2$, $Z_0=\{x_1\}$ holds in (\ref{def:beta.z0}) and
\begin{equation}\label{1:H}
y_0 \,\ \text{is the unique and non-degenerate critical point of}\,\ H_1(y),
\end{equation}
where $H_1(y)$ is given by (\ref{def:unique.Hy}).
Assume also that
there exist $\kappa>0$ and $r_0>0$ such that
\begin{equation}\label{eq1.39}
 V_i(x)\leq C e^{\kappa |x|}\,  \text{  if }\,  |x| \text{ is large},
\end{equation}
and
\begin{equation}\label{eq1.40}
\frac{\partial V_i(x+x_1)}{\partial x_j}=\frac{\partial  V_{i1}(x)}{\partial x_j}+R_{ij}(x)\,  \text{ and }\,\  |R_{ij}(x)|\leq C|x|^{q_i} \,  \text{ in }\,    B_{r_0}(0),\,\
\end{equation}
where $q_i>p_{i1}-1$ and $i, j=1,\,2$.
Then for any given $a_1\in (0, a^*)$ and $a_2\in (0, a^*)$, there exists a unique nonnegative minimizer of $e(a_1,a_2,\beta)$ as $\beta \nearrow \beta^*$.
\end{thm}

\begin{figure}
\centering
{\includegraphics[width = 10cm,height=6cm,clip]{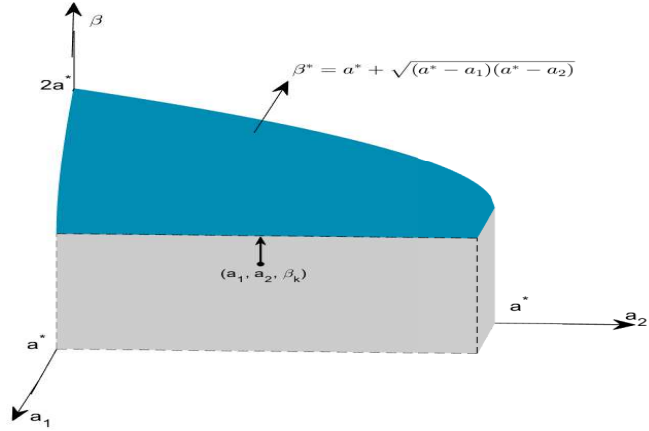}}
    \caption{\em Theorems \ref{Thm1.3}-\ref{Thm1.5} describe the spike profiles and uniqueness of nonnegative minimizers for  $e(a_1,a_2,\beta)$  as the point $(a_1,a_2,\beta)$ approaches uprightly to the top of the cuboid.}
 \end{figure}


Theorem \ref{Thm1.5} and Proposition \ref{prop:A1} imply that for any given $0<a_i<a^*$, where $i=1,2$, there exists a unique $\mu\in\R$ such that there exists a unique ground state of (\ref{equ:CGPS}) as $\beta \nearrow \beta^*$.   Theorem \ref{Thm1.5} therefore seems to be the first result on the uniqueness of  ground states for  (\ref{equ:CGPS}) with trapping potentials. On the other hand, we also note that if $V_i(x)\in C^2(\R^2)$ is homogeneous of degree $p_i\ge 2$ for $i=1,2$, then Theorem \ref{Thm1.5} can be simplified as Corollary  \ref{prop4.4}  in Section 4, which shows that  $y_{0}\not =0$ may occur in (\ref{1:H}) in view of Example 1.1 in \cite{GLW}.

To prove Theorem \ref{Thm1.5}, by contradiction suppose $(u_{1,k}, v_{1,k})$ and $(u_{2,k}, v_{2,k})$ to be two different nonnegative minimizers of $e(a_1,a_2,\beta_k)$ as $\beta _k\nearrow \beta^*$. Comparing with discussing single equations, see \cite{Cao,Deng,Grossi,GLW} and the references therein, we however need to overcome some extra difficulties. Actually, one needs to seek for a suitable difference function $( \hat\xi _{1,k},  \hat\xi _{2,k})$ between $(u_{1,k}, v_{1,k})$ and $(u_{2,k}, v_{2,k})$. By establishing Lemma \ref{lem4.2}, our analysis turns out that
$( \hat\xi _{1,k},  \hat\xi _{2,k})$ can be well defined as
\begin{equation*}
\arraycolsep=1.5pt
 \begin{array}{lll}
  \hat\xi_{1,k}(x)&=\displaystyle\frac{ u_{2,k}(x)-  u_{1,k}(x)}{\| u_{2,k}-  u_{1,k}\|^\frac{1}{2}_{L^\infty(\R^2)}\| v_{2,k}- v_{1,k}\|^\frac{1}{2}_{L^\infty(\R^2)}},\\[4.5mm]
 \hat\xi_{2,k}(x)&=\displaystyle\frac{ v_{2,k}(x)-  v_{1,k}(x)}{\| u_{2,k}-  u_{1,k}\|^\frac{1}{2}_{L^\infty(\R^2)}\| v_{2,k}-  v_{1,k}\|^\frac{1}{2}_{L^\infty(\R^2)}}.
\end{array}
\end{equation*}
It is then necessary to study carefully the limit structure of $( \hat\xi _{1k}, \hat\xi _{2k})$ as $k\to\infty$, for which we need to make full use of the non-degeneracy results in \cite{DW}.
In order to employ the non-degeneracy assumption of (\ref{1:H}) to derive Pohozaev identities, some delicate estimates and new ideas are also needed to handle with the crossing terms in BEC systems, see Step 2 in the proof for details.
We also remark that Theorems \ref{Thm1.3}-\ref{Thm1.5} describe the spike profiles and uniqueness of nonnegative minimizers for  $e(a_1,a_2,\beta)$ as the point $(a_1,a_2,\beta)$ approaches uprightly to the top of the cuboid as illustrated by Figure 2.

This paper is organized as follows.
Section 2 is focussed on the proof of Theorem \ref{Th:existence1} and Theorem \ref{Th:existence2} on classifying completely the existence and nonexistence of the minimizers for $e(a_1,a_2,\beta)$.
In Section 3 we first address the proof of Theorem \ref{Thm1.3} on the mass concentration of nonnegative minimizers for $e(a_1,a_2,\beta)$, based on which Theorem \ref{Thm1.4} is then proved in Subsection 3.1. In Section 4 we shall prove Theorem \ref{Thm1.5} on the uniqueness of nonnegative minimizers for $e(a_1,a_2,\beta)$ as $(a_1,a_2,\beta)\nearrow (a_1,a_2,\beta ^*)$. As stated in Proposition A.1, we shall derive in Appendix A.1 the equivalence between ground states  of \eqref{equ:CGPS} and constraint minimizers of $e(a_1,a_2,\beta)$. In Appendix A.2 we prove the Gagliardo-Nirenberg type inequality \eqref{Ineq:GN}, and  the proof of Lemma \ref{lem4.3} is given finally in Appendix A.3.

\section{Existence and Nonexistence of Minimizers}\label{sect:existence}
By making full use of the inequality (\ref{Ineq:GN}), this section is devoted to the proof of Theorems \ref{Th:existence1} and  \ref{Th:existence2} on the existence and non-existence of minimizers for problem \eqref{def:e}.
We start with introducing the following compactness lemma, which can be derived in a similar way of \cite[Theorem  XIII.67]{RS} or \cite[Theorem 2.1]{BW}.

\begin{lem}\label{Lem:compact}
Suppose $V_i\in L_{loc}^\infty(\R^2)$ satisfies
$\lim\limits_{|x|\to\infty}V_i(x)=\infty$ for $i=1,\,2$. Then the embedding $\mathcal{X}:=\h_1\times \h_2\hookrightarrow
L^{q}(\R^2)\times L^{q}(\R^2)$ is compact for all $2\leq q<\infty$.
\end{lem}


\noindent {\bf Proof of Theorem \ref{Th:existence1}.} Since Theorem \ref{Th:existence1}(2) is proved in  \cite[Theorem 2.2]{BC}, here for the reader's convenience we only provide a simpler proof of Theorem \ref{Th:existence1}(1) by applying the inequality (\ref{Ineq:GN}).

For any given   $0\le a_1 < a^*$, $0\le  a_2< a^*$ and $\beta<\beta^*$, there exists an $a\in\R^+$ such that $$\max\big\{a_1,a_2\big\}\leq a<a^*,\,\ \beta<a+\sqrt{(a-a_1)(a-a_2)}.$$
We thus rewrite \eqref{def:E} as
\begin{equation}\label{trans:E}
\begin{split}
E_{a_1,a_2,\beta}(u_1,u_2)
=& \int_{\R ^2} \big(|\nabla u_1|^2+|\nabla u_2|^2\big) \dx-\frac{a}{2}\int_{\R^2}(|u_1|^2+|u_2|^2)^2 \dx\\
 & +\int_{\R ^2}\big(V_1(x)|u_1|^2+V_2(x)|u_2|^2 \big)\dx\\
 &+\frac{1}{2}\int_{\R ^2} \big(\sqrt{a-a_1}|u_1|^2-\sqrt{a-a_2}|u_2|^2\big)^2 \dx\\
 &+\big(a+\sqrt {(a-a_1)(a-a_2)}-\beta\big) \int_{\R ^2}|u_1|^2|u_2|^2 \dx,\,\ (u_1,u_2)\in\mathcal{X}.
\end{split}
\end{equation}
Let $\{(u_{1n},u_{2n})\}\in\mathcal{X}$ be a minimizing sequence of $e(a_1,a_2,\beta)$, 
we then get from  (\ref{trans:E}) and the Gagliardo-Nirenberg type inequality \eqref{Ineq:GN}  that $(u_{1n},u_{2n})$ is bounded uniformly in $\mathcal{X}$.
Applying Lemma \ref{Lem:compact}, this yields that there exists $(u_{10},u_{20})\in\x$ such that, up to a subsequence if necessary,
\begin{equation}\label{lim:uin.Lq}
(u_{1n},u_{2n}) \xrightarrow{n} (u_{10},u_{20})\,\  \mbox {strongly in}\,\  L^q(\R^2),\,\ \mbox {where}\,\  2\leq q<\infty.
\end{equation}
Then, similar to \cite[Theorem 1]{GS}, one can easily deduce that  $(u_{10},u_{20})\in\mathcal{M}$ is a minimizer of $e(a_1,a_2,\beta)$.\qed




\vskip 0.1truein

\noindent {\bf Proof of (1) and (2) of Theorem \ref{Th:existence2}.}
Choose a cutoff function $0\leq \varphi \in C_0^\infty(\R^2)$ such that $\varphi(x)=1$ for $|x|\leq 1$, and $\varphi(x)=0$ for $|x|\geq 2$.
For any $\tau>0$, set
\begin{equation}\label{def:trial}
\begin{split}
  u_{1\tau}(x)&=\sqrt{\theta}A_\tau\frac{\tau}{\| w\|_2}\varphi(x-x_0) w(\tau|x-x_0|), \\
  u_{2\tau}(x)&=\sqrt{1-\theta}A_\tau\frac{\tau}{\| w\|_2}\varphi(x-x_0) w(\tau|x-x_0|),
\end{split}
\end{equation}
where $A_{\tau}>0$ is chosen so that $(u_{1\tau},u_{2\tau})\in\mathcal{M}$, $\theta \in [0,1]$ and $x_0\in\R^2$ are to be determined later. Using (\ref{decay:w}), one can verify that
\begin{equation*}
 1 \leq A_\tau^2 \leq 1+Ce^{-2\tau} \, \ \mbox{as}\,\ \tau\to\infty.
\end{equation*}
Similar to  \cite[Theorem 1]{GS}, one can employ \eqref{ide:w} to derive that
\begin{equation}\label{sup:E.w}
\begin{split}
e(a_1,a_2,\beta)\le & E_{a_1,a_2,\beta}(u_{1\tau},u_{2\tau}) \\
\leq  &  \frac{A_\tau^2 \tau^2}{\|w\|_2^2}\Big[\int_{\R ^2} |\nabla w|^2 \dx - \frac{a_1\theta^2 +a_2 (1-\theta)^2
+2\beta\theta(1-\theta)}{2a^*}\int_{\R ^2} w^4 \dx\Big]  \\
& +\int_{\R ^2}\Big[(1-\theta)\varphi^2\Big(\frac{x}{\tau}\Big) V_2\Big(\frac{x}{\tau}+x_0\Big)\\
&\qquad \quad +\theta\varphi^2\Big(\frac{x}{\tau}\Big)V_1\Big(\frac{x}{\tau}+x_0\Big)\Big]w^2 \dx+Ce^{-2\tau}\\
\leq  & \tau^2 \Big[1 - \frac{a_1\theta^2 +a_2 (1-\theta)^2+2\beta\theta(1-\theta)}{a^*}\Big]\\
&+\theta V_1(x_0)+(1-\theta)V_2(x_0)+o(1)\ \ \text{as}\,\ \tau\to\infty.
\end{split}
\end{equation}
We next follow (\ref{sup:E.w}) to continue the rest proof.

\vskip 0.1truein

\noindent{\bf  1.}
We first consider the case where $a_1=a^*$, $a_2\leq a^*$ and $\beta\leq a^*$.
Taking $\theta=1$ and $x_0\in\R^2$ such that $V_1(x_0)=0$,
we then derive from \eqref{exp2:ea} and \eqref{sup:E.w} that
\begin{equation}\label{val:e0}
 0\leq e(a^*,a_2,\beta)\leq V_1(x_0)=0.
\end{equation}
Suppose now that $e(a^*,a_2,\beta)$ has a minimizer $(u_1^*,u_2^*)$. It then follows  from  \eqref{exp2:ea} and \eqref{val:e0} that
\begin{equation}\label{equality:e0}
\begin{split}
\int_{\R ^2} (|\nabla u_1^*|^2+|\nabla u_2^*|^2) \dx&=\frac{a^*}{2}\int_{\R ^2} (|u_1^*|^2+|u_2^*|^2)^2 \dx,\\
\int_{\R ^2}V_1(x)|u_1^*|^2 \dx& = \int_{\R ^2}V_2(x)|u_2^*|^2 \dx = 0.
\end{split}
\end{equation}
Applying  \eqref{Ineq:GN}, the first equality of (\ref{equality:e0}) gives that $(u_1^*,u_2^*)$ must be a suitable scaling of $(w\sin \theta^*  , w\cos \theta^*  )$ for some $\theta^* \in [0,2\pi)$.
However, the second equality of (\ref{equality:e0}) shows that it has a  compact  support, which is impossible. Therefore, there is no minimizer for \eqref{def:e} in this case.
Furthermore, when $a_2\leq a^*$, $\beta\leq a^*$ and $a_1\nearrow a^*$, take $\theta=1$, $\tau=(a^*-a_1)^{-\frac{1}{4}}$ and $x_0\in\R^2$ satisfying $V(x_0)=0$. Using \eqref{sup:E.w}, we then obtain from (\ref{val:e0}) that
\begin{equation}\label{lim:e.0}
0\leq\lim\limits_{a_1\nearrow a^*} e(a_1,a_2,\beta)\leq \lim\limits_{a_1\nearrow a^*}E_{a_1,a_2,\beta}(u_{1\tau},u_{2\tau})=0=e(a^*,a_2,\beta),
\end{equation}
and we are done.

Similarly, one can also establish  Theorem \ref{Th:existence2} (1) for the other  case.

\vskip 0.1truein

\noindent{\bf 2.} To establish Theorem \ref{Th:existence2} (2), we denote $x_1$ a  common minimum  point of $V_1(x)$ and $V_2(x)$.
Taking $\theta=\frac{\sqrt{a^*-a_2}}{\sqrt{a^*-a_1}+\sqrt{a^*-a_2}}$  and $x_0=x_1$ for \eqref{sup:E.w}, together with \eqref{exp2:ea} we then derive  that
\begin{equation*}
 0\leq e(a_1,a_2,\beta^*)
\leq \frac{\sqrt{a^*-a_1}V_2(x_0)+\sqrt{a^*-a_2}V_1(x_0)}{(\sqrt{a^*-a_1}+\sqrt{a^*-a_2})^2}=0.
\end{equation*}
Arguing as above \eqref{equality:e0}, we further conclude that $e(a_1,a_2,\beta^*)$ has no minimizer in this case.

Furthermore, similar to \eqref{lim:e.0}, we then derive from \eqref{sup:E.w}, where we choose $\tau=(\beta^*-\beta)^{-\frac{1}{4}}$, that $\lim\limits_{\beta\nearrow \beta^*} e(a_1,a_2,\beta)=e(a_1,a_2,\beta^*)=0$.
This  completes the proof of Theorem \ref{Th:existence2} (2).
\qed

\subsection{Proof of  Theorem \ref{Th:existence2} (3) } \label{ap2}
In this subsection, we shall prove Theorem \ref{Th:existence2} (3) on the existence of  minimizers under the additional assumption (\ref{thm1.2:1A}).

\vskip 0.1truein
\noindent {\bf Proof of Theorem \ref{Th:existence2} (3).}
For any $(u_1,u_2),(v_1,v_2)\in \m$, we define
\begin{equation*}
\begin{split}
  d\big((u_1,u_2),(v_1,v_2)\big):&=\|(u_1,u_2)-(v_1,v_2)\|_\x.
  \end{split}
\end{equation*}
One can verify that the space $(\mathcal{M},d)$ is a complete metric space.
Hence, by Ekeland's variational principle \cite[Theorem  5.1]{S},
there exists a minimizing sequence  $\{(u_{1n},u_{2n})\}\subset \m$ of $e(a_1,a_2,\beta^*)$ such that
\begin{eqnarray}\label{ineq:Ekeland1}
e(a_1,a_2,\beta^*)\leq  E_{a_1,a_2,\beta^*} (u_{1n},u_{2n})\leq  e(a_1,a_2,\beta^*) +\frac{1}{n},
\end{eqnarray}
and for any $(v_1,v_2) \in \mathcal{M}$,
\begin{eqnarray}\label{ineq:Ekeland2}
  E_{a_1,a_2,\beta^*} (u_{1n},u_{2n})- E_{a_1,a_2,\beta^*} (v_1,v_2) \leq \frac{1}{n}\big\|(u_{1n},u_{2n})-(v_1,v_2)\big\|_\mathcal{X}.
\end{eqnarray}
By  Lemma \ref{Lem:compact}, in order to show the existence of the minimizers for $ e(a_1,a_2,\beta^*) $,
it suffices to prove that  $\{(u_{1n},u_{2n})\}$ is bounded uniformly in $\mathcal{X}$.
Motivated by  Theorem 1.3(ii) in \cite{GZZ2}, by contradiction we suppose that
\begin{equation}\label{lim:un.infty}
 \|(u_{1n},u_{2n})\|_\mathcal{X}\to\infty \, \ \mbox{as} \,\ n\to\infty.
\end{equation}
\noindent \textbf{Step 1.}
We have
\begin{equation}\label{lim:nablaun.infty}
\liminf\limits_{n\to\infty}\int_{\R ^2}(|\nabla u_{1n}|^2+|\nabla u_{2n}|^2) \dx=+\infty,
\end{equation}
\begin{equation}\label{ineq:u1nu2n}
  \frac{1}{2}\int_{\R ^2} \Big(\sqrt{a-a_1}|u_{1n}|^2-\sqrt{a-a_2}|u_{2n}|^2\Big)^2 \dx\leq  e(a_1,a_2,\beta^*) +\frac{1}{n},
\end{equation}
and
\begin{equation}\label{lim:nablau4}
  \liminf\limits_{n\to\infty}\frac{\int_{\R ^2} (|\nabla u_{1n}|^2+|\nabla u_{2n}|^2) \dx} {\int_{\R ^2} (|u_{1n}|^2+|u_{2n}|^2)^2 \dx}= \frac{a^*}{2}.
\end{equation}

To derive above results, we write $E_{a_1,a_2,\beta^*}(u_{1n},u_{2n})$ in the  form of (\ref{exp2:ea}).
By \eqref{Ineq:GN}, we then deduce from \eqref{ineq:Ekeland1} that \eqref{ineq:u1nu2n} holds and also
\begin{equation}\label{sup:Vun}
\int_{\R^2} \big(V_1(x)|u_1|^2+V_2(x)|u_2|^2\big) \dx\leq  E_{a_1,a_2,\beta^*} (u_{1n},u_{2n})\leq  e(a_1,a_2,\beta^*) +\frac{1}{n},
\end{equation}
which thus implies \eqref{lim:nablaun.infty} in view of \eqref{lim:un.infty}.
Together with \eqref{lim:nablaun.infty} and \eqref{ineq:u1nu2n}, we obtain from  \eqref{sup:Vun} that \eqref{lim:nablau4} holds. 
This completes the proof of Step 1.

\vskip 0.1truein

\noindent \textbf{Step 2.} Limit behavior of  minimizers.
Set
\begin{equation*}
  \varepsilon_n  := \Big(\int_{\R^2}|\nabla u_{1n}|^2+|\nabla u_{2n}|^2\dx\Big)^{-\frac{1}{2}}>0,
\end{equation*}
so that $\varepsilon_n  \to 0$ as $n\to\infty$. Similar to the proof  of \cite[Lemma 5.3]{GZZ2},
one can derive from \eqref{lim:nablaun.infty}--\eqref{lim:nablau4} that there exist a sequence $\{y_{\varepsilon_n}\}\subset\R^2$ and positive constants $R_0$ and $\eta$ such that
\begin{equation}\label{sub:w}
   \liminf\limits_{n\to\infty}\int_{B_{R_0}(0)} |w_{in}|^2 \dx \geq \eta>0,
\end{equation}
where  $w_{in}$ is defined as
\begin{equation*}
  w_{in}(x) := \varepsilon_n u_{in}(\varepsilon_nx+\varepsilon_ny_{\varepsilon_n}),\,\,\ i=1,2.
\end{equation*}
Also, we have \begin{equation}\label{val:nablawn}
  \int_{\R^2}(|\nabla w_{1n}|^2+|\nabla w_{2n}|^2)dx=1,\ \lim\limits_{n\to\infty}\int_{\R ^2} (|w_{1n}|^2+|w_{2n}|^2)^2 \dx = \frac{2}{a^*},
\end{equation}
and
\begin{equation}\label{lim:u1u24}
 \lim\limits_{n\to\infty} \int_{\R ^2} \Big(\sqrt{a^*-a_1}|w_{1n}|^2-\sqrt{a^*-a_2}|w_{2n}|^2\Big)^2 \dx=0.
\end{equation}
Moreover, we deduce from \eqref{sup:Vun} that
\begin{equation*}
\int_{\R^2} V_i(\varepsilon_n x+\varepsilon_n y_{\varepsilon_n})w_{in}(x)^2\dx\leq  e(a_1,a_2,\beta^*) +\frac{1}{n},\,\,i=1,2.
\end{equation*}
Since $\lim\limits_{|x|\to\infty}V_i(x)=+\infty$, we conclude from (\ref{sub:w}) that
 the sequence $\{\varepsilon_ny_{\varepsilon_n}\}$ is bounded uniformly in $\R^2$, and hence, up to a subsequence if necessary, we have
\begin{equation}\label{lim:y}
  \varepsilon_ny_{\varepsilon_n} \to z_0\,\ \text{for some}\,\  z_0\in\R^2 \,\   \mbox{as} \,\ n\to\infty.
\end{equation}
\noindent \textbf{Step 3.} In this step, we shall complete the proof by deriving a contradiction.
For any $\varphi_i(x)\in C_0^{\infty}(\R^2)$, set $\bar{\varphi}_i (x):=\varphi_i\big(\frac{x-\varepsilon_ny_n}{\varepsilon_n}\big)$ for $i=1, 2$
and define
\begin{equation*}
  f(\tau, \sigma):=\frac{1}{2}\int_{\R^2}\big[(u_{1n}+\tau u_{1n}+\sigma\bar{\varphi_1})^2+(u_{2n}+\tau u_{2n}+\sigma\bar{\varphi}_2 )^2\big] \dx.
\end{equation*}
Direct calculations give that
\begin{equation*}
  f(0, 0)=\frac{1}{2},\,\,\ \frac{\partial f(0,0)}{\partial \tau}=1,
\,\ \mbox{and }\,\ \frac{\partial f(0,0)}{\partial \sigma}=\int_{\R^2}\big(u_{1n}\bar{\varphi}_1+u_{2n}\bar{\varphi}_2 \big) \dx.
\end{equation*}
Applying the implicit function theorem then gives that there exist a constant $\delta_n>0$
and a function $\tau(\sigma)\in C^1\big((-\delta_n,\delta_n), \R\big)$ such that
$$\tau(0)=0,\,\ \tau'(0)=-\int_{\R^2}\big(u_{1n}\bar{\varphi}_1 +u_{2n}\bar{\varphi}_2 \big) \dx,
\,\ \text{and}\ \ f(\tau(\sigma),\sigma)=f(0,0)=\frac{1}{2}.$$
This implies that
$$\Big(u_{1n}+\tau(\sigma) u_{1n}+\sigma\bar{\varphi}_1, u_{2n}+\tau(\sigma) u_{2n}+\sigma\bar{\varphi}_2 \Big)\in \mathcal{M},
\,\ \text{where}\,\ \sigma\in (-\delta_n,\delta_n).$$
Note from \eqref{ineq:Ekeland2} that
\begin{equation*}
\begin{split}
  & E_{a_1,a_2,\beta^*} (u_{1n},u_{2n})
  - E_{a_1,a_2,\beta^*} \big(u_{1n}+\tau(\sigma) u_{1n}+\sigma\bar{\varphi}_1 , u_{2n}+\tau(\sigma) u_{2n}+\sigma\bar{\varphi}_2\big ) \\
\leq & \frac{1}{n}\|(\tau(\sigma) u_{1n}+\sigma\bar{\varphi}_1, \tau(\sigma) u_{2n}+\sigma\bar{\varphi}_2 )\|_\mathcal{X}.
\end{split}
\end{equation*}
Taking $\sigma\to0^+$ and $\sigma\to0^-$ for the above estimate, respectively, it  yields that
 \begin{equation*}
 \begin{split}
 &\Big|\langle E'_{a_1,a_2,\beta^*}(u_{1n},u_{2n}),(\tau'(0) u_{1n}+\bar{\varphi}_1 , \tau'(0) u_{2n}+\bar{\varphi}_2 )\rangle\Big|\\
\leq & \frac{1}{n}\|(\tau'(0) u_{1n}+\bar{\varphi}_1, \tau'(0) u_{2n}+\bar{\varphi}_2 )\|_\mathcal{X}.
\end{split}\end{equation*}

Since
\begin{equation*}
  \tau'(0)=-\int_{\R^2}\big(u_{1n}\bar{\varphi}_1 +u_{2n}\bar{\varphi}_2 \big) \dx
  =-\varepsilon_n\int_{\R^2}\big(w_{1n}{\varphi_1}+w_{2n}{\varphi_2}\big) \dx,
\end{equation*}
one can verify that
\begin{equation*}
  \|(\tau'(0) u_{1n}+\bar{\varphi}_1,\tau'(0) u_{2n}+\bar{\varphi}_2) \|_\mathcal{X}\leq C.
\end{equation*}
Combining with the definitions of $w_{in}$ and $\bar{\varphi_{i}}$, where $i=1,2$, we then have
\begin{equation}\label{val:tau.eqnw}
\begin{split}
      & \Big| \int_{\R^2}\big[\nabla w_{1n} \nabla\varphi_1+\nabla w_{2n} \nabla\varphi_2+\varepsilon_n^2V_1(\varepsilon_n x+\varepsilon_n y_{\varepsilon_n})w_{1n}\varphi_1\\
      &+\varepsilon^2_nV_2(\varepsilon_n x+\varepsilon_n y_{\varepsilon_n})w_{2n}\varphi_2-a^*(w_{1n}^2+w_{2n}^2)( w_{1n}\varphi_1+w_{2n}\varphi_2) \\
      &+w_{1n}\varphi_1+w_{2n}{\varphi_2}
      +\big(\sqrt{a^*-a_1}|w_{1n}|^2-\sqrt{a^*-a_2}|w_{2n}|^2\big)\\
      &\quad\big(\sqrt{a^*-a_1}w_{1n}\varphi_1-\sqrt{a^*-a_2}w_{2n}\varphi_2\big)\big] \dx\Big|
      \leq \frac{C\varepsilon_n}{n}.
\end{split}
\end{equation}
From  \eqref{sub:w} and \eqref{val:nablawn}, we know that there exists $(0,0)\not\equiv(w_1,w_2)\in H^1(\R^2)\times H^1(\R^2)$ such that   $(w_{1n},w_{2n})\overset{n}\rightharpoonup(w_1,w_2)$ in $\mathcal{X}$. Using (\ref{lim:u1u24}),
we then deduce from (\ref{val:tau.eqnw}) by letting $n\to\infty$ that
$(w_1,w_2)$ is a weak solution of the following system
\begin{equation}\label{equ:CGPS.w1w2}
\begin{cases}
-\Delta w_{1}+ w_{1}= a^*w_{1}^3+ a^*w_{2}^2w_{1}\,\,\ \mbox{in}\,\ \R^2,\\
-\Delta w_{2}+ w_{2}= a^*w_{2}^3+ a^*w_{1}^2w_{2}\,\,\ \mbox{in}\,\ \R^2.
\end{cases}
\end{equation}
Note from \eqref{lim:u1u24} that $\sqrt{a^*-a_1}|w_{1}|^2=\sqrt{a^*-a_2}|w_{2}|^2$  $a.e.$ in $\R^2$.
Set
\begin{equation}\label{def:w0}
  w_0(x)^2:=(a^*-a_2)^{-\frac{1}{2}}|w_{1}|^2=(a^*-a_1)^{-\frac{1}{2}}|w_{2}|^2.
\end{equation}
It then follows from \eqref{equ:CGPS.w1w2} that $w_0$ satisfies the following equation
\begin{equation*}
  -\Delta w_0(x)+w_0(x)=a^*\big(\sqrt{a^*-a_1}+\sqrt{a^*-a_2}\big)w_0^3(x)\,\ \text{in} \,\  \R^2.
\end{equation*}
Similar to (79) in \cite{GZZ2}, utilizing (\ref{ineq:GNQ}), one can prove that
$
 \|w_0 \|^2_2 =  \frac{1}{\sqrt{a^*-a_1}+\sqrt{a^*-a_2}}
$.
Applying  Fatou's Lemma, we then derive from  \eqref{sup:Vun}, (\ref{lim:y}) and \eqref{def:w0} that
\begin{equation*}
\begin{split}
   e(a_1,a_2,\beta^*)
 &\geq\lim\limits_{n\to\infty}\int_{\R ^2}\Big[V_1(\varepsilon_n x+\varepsilon_n y_{\varepsilon_n})|w_{1n}|^2
+V_2(\varepsilon_n x+\varepsilon_n y_{\varepsilon_n})|w_{2n}|^2\Big] \dx\\
 & \geq\frac{1}{\sqrt{a^*-a_1}+\sqrt{a^*-a_2}}\Big[\sqrt{a^*-a_2}V_1(z_0)+\sqrt{a^*-a_1}V_2(z_0)\Big],
\end{split}
\end{equation*}
which however contradicts to the assumption (\ref{thm1.2:1A}). The proof of Theorem \ref{Th:existence2} (3) is thus complete.
\qed

\section{Mass Concentration and Symmetry Breaking}\label{sect:limit1}
In this section, we focus on the proof of Theorems \ref{Thm1.3} and \ref{Thm1.4} on the mass concentration behavior and symmetry-breaking of nonnegative minimizers for $e(a_1, a_2, \beta)$ as $\beta\nearrow\beta^*$, where $0< a_1<a^*$ and $0< a_2<a^*$ are fixed. We begin with the following lemma.


\begin{lem}\label{lem:beta.u}
Suppose that $V_1(x)$ and $V_2(x)$ satisfy \eqref{cond:V.1} and have at least one common minimum point. Let $(u_{1\beta}(x),u_{2\beta}(x))$ be a nonnegative minimizer of $e(a_1,a_2,\beta)$  as $\beta\nearrow\beta^*$. Then we have
\begin{enumerate}
  \item  $(u_{1\beta}(x),u_{2\beta}(x))$ blows up in the sense that for $i=1,2$,
\begin{equation}\label{lim:beta.kinetic}
   \lim\limits_{\beta\nearrow\beta^*}\int_{\R^2}|\nabla u_{i\beta}(x)|^2\dx=+\infty \ \text{ and }\ \lim\limits_{\beta\nearrow\beta^*}\int_{\R^2}|u_{i\beta}(x)|^4\dx=+\infty.
\end{equation}
  \item $(u_{1\beta}(x),u_{2\beta}(x))$ also satisfies
\begin{equation}\label{lim.beta.V}
  \lim\limits_{\beta\nearrow\beta^*}\int_{\R^2}V_1(x)u_{1\beta}(x)^2\dx=  \lim\limits_{\beta\nearrow\beta^*}\int_{\R^2}V_2(x)u_{2\beta}(x)^2\dx=0,
\end{equation}
\begin{equation}\label{lim:beta.minus}
\lim\limits_{\beta\nearrow\beta^*}\int_{\R^2}\big(\sqrt{a^*-a_1}|u_{1\beta}(x)|^2-\sqrt{a^*-a_2}|u_{2\beta}(x)|^2\big)^2\dx=0,
\end{equation}
and further,
\begin{equation}\label{lim:beta.GN}
\lim\limits_{\beta\nearrow\beta^*}\frac{\int_{\R^2}(|\nabla u_{1\beta}(x)|^2+|\nabla u_{2\beta}(x)|^2)\dx}{\int_{\R^2}(| u_{1\beta}(x)|^2+|u_{2\beta}(x)|^2)^2\dx}=\frac{a^*}{2},
\end{equation}
\begin{equation}\label{lim:beta.u14:u24}
     \lim\limits_{\beta\nearrow\beta^*}
     \frac{\int_{\R^2}|u_{1\beta}(x)|^4\dx}{\int_{\R^2}|u_{2\beta}(x)|^4\dx}
     = \frac{{a^*-a_2}}{{a^*-a_1}}.
\end{equation}
\end{enumerate}
\end{lem}

\noindent {\bf Proof.} By rewriting  $E_{a_1,a_2,\beta}(\cdot)$ as the form of (\ref{exp2:ea}),
 (\ref{lim.beta.V}) and (\ref{lim:beta.minus}) follow  directly from \eqref{eq1.15}. Moreover, similar to  \cite[Lemma 3.1 (i)]{GWZZ}, applying \eqref{eq1.15} and Lemma \ref{Lem:compact}, one can prove by contradiction that
\begin{equation}\label{lim:beta.sumnabla}
   \lim\limits_{\beta\nearrow\beta^*}\int_{\R^2}\big(|\nabla u_{1\beta}(x)|^2+|\nabla u_{2\beta}(x)|^2\big)\dx=+\infty.
\end{equation}
It then implies from (\ref{lim:beta.sumnabla}) that (\ref{lim:beta.GN}) holds by using (\ref{exp2:ea}), (\ref{lim.beta.V}) and (\ref{lim:beta.minus}). Also,  the conclusion that $\inte|u_{i\beta}|^4dx\to+\infty$ follows from (\ref{lim:beta.minus}), \eqref{lim:beta.GN} and (\ref{lim:beta.sumnabla}). Moreover,  (\ref{lim:beta.kinetic}) can be derived by applying (\ref{ineq:GNQ}). Finally, one can show that (\ref{lim:beta.u14:u24}) holds in view of (\ref{lim:beta.kinetic}) and (\ref{lim:beta.minus}).
\qed

Define
\begin{equation*}
\varepsilon_\beta:=\Big(\int_{\R^2}\big(|\nabla u_{1\beta}(x)|^2+|\nabla u_{2\beta}(x)|^2\big)\dx\Big)^{-\frac{1}{2}}>0,
\end{equation*}
so that $\varepsilon_\beta\to0$ as $\beta\nearrow\beta^*$ by Lemma \ref{lem:beta.u}.

\begin{lem}\label{lem:beta.w}
Under the assumptions of  Lemma \ref{lem:beta.u},  we have
\begin{enumerate}
\item
  There exist a sequence $\{y_{\varepsilon_\beta}\}\subseteq\R^2$ and positive constants $R_0$ and $\eta$ such that
\begin{equation}\label{sub:wibeta}
  \liminf\limits_{\beta\nearrow\beta^*}\int_{B_{R_0}(0)}w_{i\beta}(x)^2\dx\geq \eta>0, \,\ i=1,2,
\end{equation}
where $w_{i\beta}(x)$ is defined as
\begin{equation}\label{def:wbeta}
  w_{i\beta}(x):=\varepsilon_\beta u_{i\beta}(\varepsilon_\beta x+\varepsilon_\beta y_{\varepsilon_\beta}), \,\ i=1,2.
\end{equation}
Moreover,
\begin{equation}\label{lim:wbeta}
\int_{\R^2}\big(|\nabla{w}_{1\beta}|^2+|\nabla{w}_{2\beta}|^2\big)\dx=1,\,\
\lim\limits_{\beta\nearrow\beta^*}{\int_{\R^2}\big(|{w}_{1\beta}|^2+|{w}_{2\beta}|^2\big)^2\dx}=\frac{2}{a^*},
\end{equation}
and
\begin{equation}\label{lim:wbeta.minus}
\lim\limits_{\beta\nearrow\beta^*}\int_{\R^2}\big(\sqrt{a^*-a_1}|w_{1\beta}(x)|^2-\sqrt{a^*-a_2}|w_{2\beta}(x)|^2\big)^2\dx=0.
\end{equation}

\item
For any sequence $\{\beta_k\}$ with $\beta_k\nearrow\beta^*$ as $k\to\infty$, there exists a subsequence, still denoted by $\{\beta_k\}$, such that $\varepsilon_{\beta_k} y_{\varepsilon_{\beta_k}} \to \bar{x}_0$ as $k\to\infty$, where $\bar{x}_0\in\R^2$ satisfies $V_1(\bar{x}_0)=V_2( \bar{x}_0)=0$.

  \item
  For any sequence $\{\beta_k\}$ satisfying $\beta_k\nearrow\beta^*$ as $k\to\infty$, there exist  a subsequence, still denoted by $\{\beta_k\}$, and $x_0\in\R^2$ such that
\begin{equation}\label{lim:beta.w1w2.exp}
\begin{cases}
 \lim\limits_{k\to\infty} w_{1k}(x) = \frac{\sqrt{\gamma}}{\|w\|_2}w(x-x_0),\\
 \lim\limits_{k\to\infty} w_{2k}(x) =  \frac{\sqrt{1-\gamma}}{\|w\|_2}w(x-x_0)
\end{cases}
\text{strongly in }H^1(\R^2),
\end{equation}
where $(w_{1k},w_{2k})=(w_{1\beta_k},w_{2\beta_k})$.
\end{enumerate}
\end{lem}

\noindent {\bf Proof.}
Since the proof of (1) and (2) in Lemma \ref{lem:beta.w} is similar to that of \cite[Lemma 5.3]{GZZ2}, here we only address the proof of (3) in Lemma \ref{lem:beta.w}. Note that the nonnegative minimizer $(u_{1\beta},u_{2\beta})$ solves the system
\begin{equation}\label{3:sys}
 \begin{cases}
-\Delta u_{1\beta} +V_1(x)u_{1\beta} =\mu_\beta u_{1\beta} +a_1u_{1\beta}^3 +\beta u_{2\beta}^2 u_{1\beta}   \,\ \mbox{in}\,\  \R^2,\\
-\Delta u_{2\beta} +V_2(x)u_{2\beta} =\mu_\beta u_{2\beta} +a_2u_{2\beta}^3 +\beta u_{1\beta}^2 u_{2\beta}   \,\ \mbox{in}\,\  \R^2,\,\
\end{cases}
\end{equation}
where $\mu _{\beta}\in \R$ is a suitable Lagrange multiplier. By (\ref{def:wbeta}), we derive from (\ref{3:sys}) that $(w_{1\beta} ,w_{2\beta})$ satisfies the following system
\begin{equation}\label{equ:CGPS.wbeta}
\begin{cases}
 -\Delta w_{1\beta} +\varepsilon_\beta^2V_1(\varepsilon_\beta x+\varepsilon_\beta y_{\varepsilon_\beta})w_{1\beta}
 =\varepsilon_\beta^2\mu_{\beta} w_{1\beta} +a_1w_{1\beta} ^3+\beta w_{2\beta} ^2w_{1\beta}    \,\ \mbox{in}\,\  \R^2,\\
 -\Delta w_{2\beta} +\varepsilon_\beta^2V_2(\varepsilon_\beta x+\varepsilon_\beta y_{\varepsilon_\beta})w_{2\beta}
 =\varepsilon_\beta^2\mu_{\beta} w_{2\beta} +a_2w_{2\beta} ^3+\beta w_{1\beta} ^2w_{2\beta}    \,\ \mbox{in}\,\  \R^2.
\end{cases}
\end{equation}

We next follow the system (\ref{equ:CGPS.wbeta}) to proceed  the further analysis of $(w_{1\beta} ,w_{2\beta})$  as $\beta\nearrow\beta^*$. From (\ref{lim:wbeta}) we see that
there exist a sequence $\{\beta_k\}$, where $\beta_k\nearrow\beta^*$ as $k\to\infty$,
and $(w_{10},w_{20})\in H^1(\R^2)$ such that
$(w_{1\beta_k},w_{2\beta_k})\overset{k}\rightharpoonup(w_{10},w_{20})$ in $H^1(\R^2)$, where $w_{10}\geq0$ and $w_{20}\geq0$ are due to the nonnegativity  of $ w_{1\beta_k}$ and $ w_{2\beta_k}$.
Using \eqref{lim.beta.V}, \eqref{lim:wbeta} and \eqref{3:sys}, we have
\begin{equation}\label{lim:beta.mu}
\begin{split}
  \varepsilon_\beta^2\mu_{\beta}
   =& 2\varepsilon_\beta^2e(a_1,a_2,\beta) -\int_{\R^2}\big(|\nabla w_{1\beta}(x)|^2+ |\nabla w_{2\beta}(x)|^2\big)\dx\\
  &- \sum_{i=1}^2\int_{\R^2} V_i(\varepsilon_\beta x+\varepsilon_\beta y_{\varepsilon_\beta})|w_{i\beta}(x)|^2\dx\to-1\ \text{ as }\, \beta\nearrow\beta^*.
\end{split}
\end{equation}
Set now $w_{1k}:=w_{1\beta_k}$, $w_{2k}:=w_{2\beta_k}$, $\varepsilon_{k}:=\varepsilon_{\beta_k}>0$ and $\mu _k:=\mu_{\beta_k}$.
Moreover, since it yields from (\ref{lim:wbeta.minus})  that
 $ \sqrt{a^*-a_1}w_{10}(x)^2=\sqrt{a^*-a_2}w_{20}(x)^2$ $a.e.$ in $ \R^2$,  we set
\begin{equation}\label{def:beta.w0}
  w_0(x):=(a^*-a_2)^{-\frac{1}{4}}w_{10}(x)=(a^*-a_1)^{-\frac{1}{4}}w_{20}(x)>0 \,\ \text{a.e. in} \,\ \R^2.
\end{equation}
It then follows  from \eqref{sub:wibeta}, \eqref{equ:CGPS.wbeta} and Lemma \ref{lem:beta.w} (2) that $w_0>0$ satisfies
\begin{equation*}
  -\Delta w_0 +w_0 =a^*\big(\sqrt{a^*-a_1}+\sqrt{a^*-a_2}\big)w_0 ^3\,\ \text{a.e. in} \,\ \R^2.
\end{equation*}
Furthermore, we obtain from  the uniqueness (up to translations) of  positive solutions of \eqref{equ:w} that
\begin{equation*}
  w_0(x)=\frac{1}{(a^*)^\frac{1}{2}(\sqrt{a^*-a_1}+\sqrt{a^*-a_2})^\frac{1}{2}}w(x-x_0)\,\ \text{for some}\,\  x_0\in\R^2.
\end{equation*}
We thus deduce from \eqref{def:beta.w0} that
\begin{equation}\label{exp:beta.w}
 w_{10}(x)=\frac{\sqrt{\gamma}}{\|w\|_2}w(x-x_0) \,\,\text{ and } \,\, w_{20}(x)=\frac{\sqrt{1-\gamma}}{\|w\|_2}w(x-x_0).
\end{equation}
This indicates that
$\|w_{10}\|^2_2 + \|w_{20}\|^2_2=1$.
Therefore,
\begin{equation*}
  \lim\limits_{k\to\infty}(w_{1k}, w_{2k})=(w_{10}, w_{20}) \,\ \text{strongly in}\,\ L^2(\R^2)\times L^2(\R^2).
\end{equation*}
Using the interpolation inequalities  and combining  \eqref{equ:CGPS.wbeta} with \eqref{def:beta.w0}, one can further derive that \eqref{lim:beta.w1w2.exp} holds true, and  the proof of the lemma is therefore complete.
\qed

\vskip 0.1truein

\noindent {\bf Proof of Theorem \ref{Thm1.3}.}
From \eqref{equ:CGPS.wbeta}, one can derive that $w_{ik}$ satisfies
\[-\Delta w_{ik}\leq c_{ik}(x) w_{ik}\,\ \text{in} \,\ \R^2,\,\   i=1,2,\]
where
$$c_{1k}(x)=a_1w_{1k}^2+\beta_kw_{2k}^2\, \text{ and }\, c_{2k}(x)=a_2w_{2k}^2+\beta_kw_{1k}^2.$$
Applying the De Giorgi--Nash--Moser theory (cf. \cite[Theorem 4.1]{HL} or \cite[Theorem 8.15]{GT}), it then yields that
\begin{equation}\label{decay:beta.w}
  w_{ik}(x)\to 0\,\ \text{as}\,\ |x|\to\infty\,\ \text{uniformly on}\,\ k, \,\  \ i=1, 2.
\end{equation}
Thus,  $w_{ik}(x)$ admits at least one global maximum point for $i=1, 2$.

Let $\bar{z}_{ik}$ be a global maximum point of $u_{ik}(x)$ and set ${z}_{k}:=\varepsilon_k y_{\varepsilon_{\beta_k}}$,  where $i=1, 2$.
Since $w_{ik}(x)$ attains its global maximum at the point
$x = \frac{\bar{z}_{ik}-{z}_{k}}{\varepsilon_{k}}$, one can deduce from \eqref{sub:wibeta} and \eqref{decay:beta.w} that
\begin{equation}\label{lim:beta.zz.bdd}
\frac{\bar{z}_{ik}-{z}_{k}}{\varepsilon_{k}} \,\ \text{is bounded uniformly as}\,\  k\to\infty, \,\ \text{where}\,\  i=1, 2.
\end{equation}
Define
\begin{equation}\label{def:beta.w.bar}
  \bar{w}_{ik}(x)
 :=\varepsilon_{k} u_{ik}(\varepsilon_{k} x+\bar{z}_{ik})
  =w_{ik}\Big(x+\frac{\bar{z}_{ik}-{z}_{k}}{\varepsilon_{k}}\Big), \,\ i=1, 2.
\end{equation}
Using Lemma \ref{lem:beta.w} (3) yields that
\begin{equation}\label{lim:beta.barw}
  \bar{w}_{ik}(x) \to \bar{w}_{i0}(x)={w}_{i0}(x+y_i)\,\
\text{strongly in $H^1(\R^2)$ as $k\to\infty$},
\end{equation}
where ${w}_{i0}$ is defined by \eqref{exp:beta.w}
and
\begin{equation}\label{eq3.24}
y_i=\lim\limits_{k\to\infty}\frac{\bar{z}_{ik}-{z}_{k}}{\varepsilon_{k}},\,\,\,\, i=1,2.
\end{equation}
Thus, we get from \eqref{exp:beta.w} and (\ref{lim:beta.barw}) that
\begin{equation}\label{lim:beta.barw1w2.exp1}
\bar w_{10}(x)=\frac{\sqrt{\gamma}}{\|w\|_2}w(x+y_1-x_0) \text{ and }  \bar w_{20}(x)=\frac{\sqrt{1-\gamma}}{\|w\|_2}w(x+y_2-x_0).
\end{equation}
Furthermore,
since $V_i(x)$ is locally Lipschitz continuous in $\R^2$, using the standard elliptic regularity theory, one can deduce from \eqref{equ:CGPS.wbeta}, \eqref{lim:beta.zz.bdd} and \eqref{def:beta.w.bar} that
\begin{equation}\label{lim:beta.wbar.C}
\bar{w}_{i{k}}(x) \to \bar{w}_{i0}(x)\,\ \text{in}\,\  C^2_{loc}(\R^2)\,\ \text{as}\,\  k\to\infty,\,\ \text{where}\,\  i=1, 2.
\end{equation}
See \cite[Lemma 3.1]{GZZ} for the details of  a similar proof.

Finally, we prove the local uniqueness of maximum points for $u_{ik}$.
Suppose $x_{ik}$ is any local maximum point of $\bar{w}_{ik}(x)$. Replacing $(w_{1k},w_{2k})$ by $(\bar{w}_{1k},\bar{w}_{2k})$ for  \eqref{equ:CGPS.wbeta}, one can derive from \eqref{lim:beta.mu} that for sufficiently large $k>0$,
$$a_1\bar w_{1k}^2(x_{1k})+\beta_k\bar w_{2k}^2(x_{1k})\geq\frac12 \text{ and }a_2\bar w_{2k}^2(x_{2k})+\beta_k\bar w_{1k}^2(x_{2k})\geq\frac12.$$
It then follows from (\ref{decay:beta.w}) and \eqref{lim:beta.zz.bdd} that  $\{x_{ik}\}$ is bounded uniformly  as $k\to\infty$.
Similar to \cite[Theorem 1.2]{GZZ}, one can employ \eqref{lim:beta.wbar.C} to further obtain the uniqueness of maximum points of  $\bar{w}_{ik}$, $i.e.$, the uniqueness of maximum points of $u_{ik}$, where  $i=1, 2$.
Moreover, since the origin is the unique maximum point of
both $\bar{w}_{1k}$ and $\bar{w}_{2k}$, it  is also the unique maximum point of both $\bar{w}_{10}$ and $\bar{w}_{20}$, which thus implies that $y_1=y_2=x_0$,  and thus \eqref{eq1.22} follows from \eqref{eq3.24}.  Moreover, \eqref{lim:beta.u.exp} follows directly from \eqref{def:beta.w.bar}, \eqref{lim:beta.barw} and \eqref{lim:beta.barw1w2.exp1}. Also, we derive from \eqref{lim:beta.zz.bdd} and Lemma \ref{lem:beta.w} (2) that \eqref{lim:beta.z} holds.
The proof of Theorem \ref{Thm1.3} is therefore finished.
\qed

\subsection{Proof of Theorem \ref{Thm1.4}}\label{sect:limit1.beta.V}
In this subsection, we shall complete the proof of Theorem \ref{Thm1.4} under the assumptions of (\ref{def:beta.z})--(\ref{eq1.27}). Following Theorem \ref{Thm1.3}, the key part of proving Theorem \ref{Thm1.4} is to derive the explicit blow-up rate (\ref{def:beta.V.eps}). For this purpose, we first  analyze the refined  estimates of the  energy $ e(a_1,a_2, \beta)$ as $\beta \nearrow \beta^*$.

\begin{lem}\label{4A:lem1}
Suppose $V_1(x)$ and $V_2(x)$ satisfy  (\ref{def:beta.z})--(\ref{eq1.27}). Then for any given $0<a_1, a_2<a^*$, we have
\begin{equation}\label{psup:beta.e}
e(a_1,a_2,\beta)\leq\frac{{p}_0+2}{{p}_0a^*}\Big(\frac{{p}_0}{2}\bar{\lambda}_0\Big)^\frac{2}{{p}_0+2}\Big[2\gamma(1-\gamma)+o(1)\Big]^\frac{{p}_0}{{p}_0+2}
(\beta^*-\beta)^\frac{{p}_0}{{p}_0+2}\, \ \mbox{as}\, \ \beta \nearrow \beta^*,
\end{equation}
where $\gamma\in(0,1)$, ${p}_0>0$ and $\bar{\lambda}_0>0$  are defined in \eqref{def:beta.k}, (\ref{def:beta.p0}) and (\ref{def:beta.gamma}), respectively.
\end{lem}

\noindent {\bf Proof.}
Let $(u_{1\tau},u_{2\tau})$ be the trial function defined by \eqref{def:trial}. Set $\theta=\gamma$ and take $x_0=x_{j_0}+y_{j_0}/\tau$, where
\begin{equation}\label{eq3.26}
x_{j_0}\in Z_0 \text{ with }Z_0 \text{ given by \eqref{def:beta.z0} and } y_{j_0}\in \R^2 \text{ such that } H_{j_0}(y_{j_0})=\bar\lambda_0.\end{equation}
Similar to \eqref{sup:E.w}, we have
\begin{equation*}
\begin{split}
 &\int_{\R^2} \big(|\nabla u_{1\tau}|^2+ |\nabla  u_{2\tau}|^2\big)\dx-\int_{\R^2}\big(\frac{a_1}{2}| u_{1\tau}|^4
    +\frac{a_2}{2}| u_{2\tau}|^4
   +\beta | u_{1\tau}|^2| u_{1\tau}|^2\big)\dx\\
&\leq    \frac{2\gamma(1-\gamma)}{a^*}(\beta^*-\beta)\tau^2+C\tau^2 e^{-2 \tau}\,\  \mbox{as}\,\  \tau\to\infty.
\end{split}
\end{equation*}
On the other hand, in view of  (\ref{eq1.26}) and (\ref{eq3.26}), we derive that
\begin{equation*}
\begin{split}
  \sum_{i=1}^2\int_{\R^2} V_i(x)|u_{i\tau}|^2\dx
  \leq\frac{(1+o(1))H_{j_0}(y_{j_0})}{a^*\tau^{p_{0}}}=\frac{(1+o(1))\bar \lambda_0}{a^*\tau^{p_{0}}}\,\  \mbox{as}\,\  \tau\to\infty.
\end{split}
\end{equation*}
Therefore,
\begin{equation}\label{sup:beta.V5}
e(a_1,a_2,\beta) \le E_{a_1,a_2,\beta}(\phi_1,\phi_2)
  \leq\frac{2\gamma(1-\gamma)}{a^*}(\beta^*-\beta)\tau^2 +
  \frac{(1+o(1))\bar{\lambda}_{0}}{a^*\tau^{p_{0}}}\  \text{ as }\tau\to\infty.
\end{equation}
Setting $\tau=\big[\frac{{p}_0\bar{\lambda}_0}{4\gamma(1-\gamma)(\beta^*-\beta)}\big]^\frac{1}{{p}_0+2}$ into \eqref{sup:beta.V5} then yields   the estimate  (\ref{psup:beta.e}), and the proof is complete.\qed

We next follow Theorem \ref{Thm1.3} and Lemma \ref{4A:lem1} to prove Theorem \ref{Thm1.4}.

\vskip 0.1 truein

\noindent {\bf Proof of Theorem \ref{Thm1.4}.} Let $(u_{1\beta_k},u_{2\beta_k})$ be the convergent subsequence obtained in Theorem \ref{Thm1.3}. From \eqref{exp2:ea} and \eqref{def:beta.w.bar}, we have
\begin{equation}\label{sub:beta.e.exp}
\begin{split}
  e(a_1,a_2,\beta_k)
  \geq& (\beta^* - \beta_k) \int_{\R^2} |u_{1\beta_k}|^2|u_{2\beta_k}|^2\dx+\sum_{i=1}^2\int_{\R^2} V_i(x)|u_{i\beta_k}(x)|^2\dx\\
  =&\frac{(\beta^* - \beta_k)}{\varepsilon_k^2}\int_{\R^2} |\bar{w}_{1k}|^2|\bar{w}_{2k}|^2\dx+\sum_{i=1}^2\int_{\R^2} V_i(\varepsilon_k x+\bar{z}_{ik})|\bar{w}_{ik}(x)|^2\dx.
\end{split}
\end{equation}
Applying \eqref{ide:w} and \eqref{lim:beta.u.exp}, there holds that
\begin{equation}\label{lim:beta.mult.exp}
  \lim\limits_{k\to\infty}\int_{\R^2} |\bar{w}_{1k}(x)|^2|\bar{w}_{2k}(x)|^2\dx =\frac{2\gamma(1-\gamma)}{a^*}.
\end{equation}
We get from \eqref{lim:beta.z} that  $\bar{x}_0\in Z$, where $Z$ is defined in \eqref{def:beta.z}. Thus,
  $\bar{x}_0=x_{j_0}$ for some $1\leq j_0\leq l$. Without loss of generality,  suppose that ${p}_{j_0}={p}_{1j_0}\leq {p}_{2j_0}$.
We now claim that
\begin{equation}\label{eq3.35}
\text{$\frac{\bar{z}_{1k}- x_{1j_0} }{\varepsilon_{k}}$ \,\,\,is bounded uniformly as $k\to\infty$,}\end{equation}
and \begin{equation} \label{eq3.36}  {p}_{j_0}= {p}_0.
\end{equation}
By contradiction, suppose that either  $p_{j_0}<{p}_0$ or \eqref{eq3.35} does not hold.  Then, for any given $M>0$, it follows from (\ref{eq1.26}) and Fatou's Lemma that
\begin{equation}\label{sub:c.beta.V1}
\begin{split}
 & \liminf\limits_{k\to\infty}\varepsilon_{k}^{-{p}_{0}}\int_{\R^2}  V_1(\varepsilon_{k} x+\bar{z}_{1k})|\bar{w}_{1k}(x)|^2\dx\\
\geq &
  \int_{\R^2}\liminf\limits_{k\to\infty}\Big(\frac{ V_1(\varepsilon_{k} x+\bar{z}_{1{k}})}{V_{1j_0}(\varepsilon_{k} x+\bar{z}_{1{k}}-x_{1j_0})}
  V_{1j_0}\big(x+\frac{\bar{z}_{1{k}}-x_{1j_0}}{\varepsilon_{k}}\big)|\bar{w}_{1{k}}(x)|^2\Big)\dx
 \geq C M.
\end{split}
\end{equation}
We then derive from \eqref{sub:beta.e.exp}, (\ref{lim:beta.mult.exp}) and \eqref{sub:c.beta.V1} that
\begin{equation*}\label{sub:c.beta.e}
  e(a_1,a_2,\beta_k)
   \geq C'(\beta^* - \beta_k)\varepsilon_{k}^{-2} + CM \varepsilon_{k}^{{p}_{0}}
   \geq C''M^{\frac{2}{{p}_{0}+2}}(\beta^*-\beta_k)^{\frac{{p}_{0}}{{p}_{0}+2}},
\end{equation*}
which however contradicts to \eqref{psup:beta.e}. Thus, the claims \eqref{eq3.35} and \eqref{eq3.36} are proved.

Following \eqref{eq1.22} and \eqref{eq3.35}, up to a subsequence of $\{\beta_k\}$, there exists  $z_0\in\R^2$ such that
\begin{equation}\label{sub:c.beta.ze}
\lim\limits_{k\to\infty}\frac{\bar{z}_{ik}-x_{1j_0}}{\varepsilon_{k}}=z_0, \text{ where } i=1,2.
\end{equation}
Similar to \eqref{sub:c.beta.V1},  we also get from \eqref{eq3.36} that
\begin{equation}\label{sub:beta.V1}
\begin{split}
 & \liminf\limits_{k\to\infty}\varepsilon_{k}^{-{p}_{0}}
 \Big(\int_{\R^2}  V_1(\varepsilon_{k} x+\bar{z}_{1k})|\bar{w}_{1k}(x)|^2\dx
 +\int_{\R^2}  V_2(\varepsilon_{k} x+\bar{z}_{2k})|\bar{w}_{2k}(x)|^2\dx\Big)\\
 &= \frac{1}{a^*}H_{j_0}(z_0)\geq \frac{\bar{\lambda}_{j_0}}{a^*}\geq\frac{\bar{\lambda}_{0}}{a^*},
\end{split}
\end{equation}
where  $\bar{\lambda}_{j_0}$ and $ \bar{\lambda}_0$ are defined in (\ref{def:beta.gamma}).
Here it needs to note that all equalities in \eqref{sub:beta.V1} hold  if and only if $H_{j_0}(z_0)=\bar \lambda_{j_0}$ and $x_{j_0}\in Z_0$.
Together with \eqref{sub:beta.e.exp}, we deduce from \eqref{lim:beta.mult.exp} and \eqref{sub:beta.V1} that
\begin{equation}\label{psub:beta.e.eta}
  e(a_1,a_2,\beta_k)
\geq\Big(\frac{2\gamma(1-\gamma)}{a^*}+o(1)\Big)(\beta^* - \beta_k)\varepsilon_{k}^{-2}+\Big(\frac{\bar{\lambda}_0}{a^*}+o(1)\Big)\varepsilon_{k}^{{p}_0}\,\  \mbox{as}\,\  \tau\to\infty.
\end{equation}
Taking the infimum over $\varepsilon_{k}>0$ for the right hand side of (\ref{psub:beta.e.eta}), 
it thus yields that
\begin{equation}\label{psub:beta.e}
\liminf\limits_{k\to\infty}\frac{e(a_1,a_2,\beta_k)}{(\beta^*-\beta_k)^\frac{{p}_0}{{p}_0+2}}
\geq\frac{{p}_0+2}{{p}_0a^*}\Big(\frac{{p}_0}{2}\bar{\lambda}_0\Big)^\frac{2}{{p}_0+2}\Big[2\gamma(1-\gamma)\Big]^\frac{{p}_0}{{p}_0+2},
\end{equation}
where the equality holds if and only if $H_{j_0}(z_0)=\bar \lambda_{j_0}=\bar \lambda_0$ and
\begin{equation}\label{eq3.42}
\lim_{k\to\infty}\varepsilon_{k}/\bar \varepsilon_{k}=1 \text{ with }\bar{\varepsilon}_k>0 \text{ given by } \eqref{def:beta.V.eps}.\end{equation}
Combining  with \eqref{psup:beta.e}, we deduce that all equalities  in \eqref{sub:beta.V1}--\eqref{eq3.42} hold and $H_{j_0}(z_0)=\bar \lambda_{j_0}=\bar \lambda_0$.
Thus,  \eqref{lim:beta.V.u.exp} and \eqref{def:beta.V.eps} follow from \eqref{lim:beta.u.exp}, \eqref{thm1.3:1} and \eqref{eq3.42}. Also, \eqref{lim:beta.V.y0} follows from \eqref{sub:c.beta.ze}.\qed

\section{Uniqueness of Minimizers as $\beta\nearrow \beta^*$}

Under the non-degeneracy assumption  (\ref{1:H}), in this section we shall prove Theorem \ref{Thm1.5} on the uniqueness of nonnegative minimizers of $e(a_1,a_2,\beta)$ as $\beta \nearrow\beta^*$, where $0< a_1<a^*$ and $0< a_2<a^*$ are fixed. Applying Theorem \ref{Thm1.4}, we first derive the following blow-up estimates of nonnegative minimizers.

\begin{lem}\label{lem4.1}
Suppose that $V_i(x)\in C^2(\R^2)$ satisfies  (\ref{def:beta.z})--(\ref{eq1.27}) and (\ref{eq1.39}) for $i=1,\,2$, $Z_0=\{x_1\}$ holds in (\ref{def:beta.z0}) and
\begin{equation}\label{y:HH}
y_0 \,\ \text{is the unique  critical point of}\,\ H_1(y),
\end{equation}
where $H_1(y)$ is given by (\ref{def:unique.Hy}).
Let $(u_k,v_k)$ be a nonnegative minimizer of $e(a_1,a_2,\beta_k)$, where   $0< a_1, a_2<a^*$ are fixed and  $\beta_k\nearrow\beta^*$ as $k\to\infty$.
Then there exists a subsequence of $\{\beta _k\}$, still denoted by
$\{\beta _k\}$, such that
\begin{enumerate}
\item The minimizer
$(u_k,v_k)$ satisfies
\begin{equation}\label{lem4.1:1}
 \bar{u}_k:=\varepsilon_{k} u_k(\varepsilon_{k} x+x_k)
\to \frac{\sqrt{\gamma}}{\|w\|_2}w(x)\,\ \text{and}\ \
\bar{v}_k:=\varepsilon_{k} v_k(\varepsilon_{k} x+y_k)\to  \frac{\sqrt{1-\gamma}}{\|w\|_2}w(x)
\end{equation}
uniformly in $ \R^2 $ as $k\to\infty$,
where $0<\gamma= \gamma(a_1,a_2)<1$ is given by (\ref{def:beta.k}) and the unique maximum point $(x_k,y_k)$ of $(u_k,v_k)$ satisfies
\begin{equation}\label{lem4.1:3}
  \lim\limits_{k\to\infty}\frac{x_k-x_1}{\varepsilon_{k}}=
  \lim\limits_{k\to\infty}\frac{y_k-x_1}{\varepsilon_{k}}=y_0.
\end{equation}
Here   $\varepsilon_{k}>0$ satisfies
\begin{equation}\label{lem4.1:2}
  \varepsilon_{k}:=\frac{1}{\lambda}(\beta^*-\beta_k)^\frac{1}{2+p_0}>0, \,\ p_0=\min\big\{p_{11}, p_{21}\big\}>0,
\end{equation}
and
\begin{equation}\label{def:beta.gamma.unique}
 \lambda:=\Big[\frac{p}{4\gamma(1-\gamma)}H_1(y_0)\Big]^\frac{1}{2+p}>0.
\end{equation}

\item The minimizer $(u_k,v_k)$ decays exponentially in the sense that
\begin{equation}
\bar{u}_k(x)\le Ce^{-\frac{|x|}{2}} \,\ \text{and}\,\,\ \bar{v}_k(x)\le Ce^{-\frac{|x|}{2}} \,\ \text{as}\,\ |x|\to\infty,
\label{2:conexp}
\end{equation}
and
\begin{equation}
|\nabla \bar{u}_k|\le Ce^{-\frac{|x|}{4}} \,\ \text{and}\,\,\ |\nabla \bar{v}_k|\le Ce^{-\frac{|x|}{4}}\,\ \text{as}\,\ |x|\to\infty,
\label{2:conexp2}
\end{equation}
where the constant $C>0$ is independent of $k$.
\end{enumerate}
\end{lem}

\noindent\textbf{Proof.}  In view of Theorem \ref{Thm1.4},
one can deduce that \eqref{lem4.1:1} holds in $H^1(\R^2)$, where $\varepsilon_k>0$ satisfies (\ref{lem4.1:2}) and $(x_k,y_k)$ satisfies (\ref{lem4.1:3}).
We now address the exponential decay  \eqref{2:conexp} and \eqref{2:conexp2} of nonnegative minimizers $(u_k, v_k)$ as $k\to\infty$.
Without loss of generality, here we just prove the results of $u_k$.
It is easy to verify from the previous section that $\bar{u}_{k}$ satisfies the following equation
\begin{equation}\label{4:star}
-\Delta \bar{u}_{k} +\varepsilon_k^2V_1\big(\varepsilon_kx+ x_k \big)\bar{u}_{k} =\varepsilon_k^2\mu_k \bar{u}_{k} +a_1\bar{u}_{k}^3 +\beta_k \bar{v}_{k}^2 \bar{u}_{k}   \,\ \mbox{in}\,\  \R^2,
\end{equation}
where $\varepsilon_k^2\mu_k\to -1$ as $k\to\infty$.
Following these, one can derive from \eqref{lim:beta.mu}, \eqref{decay:beta.w} and \eqref{lim:beta.zz.bdd}  that there exists a constant $R>0$ large enough that
$$ -\Delta \bar{u}_{k}+\frac{1}{2}\bar{u}_{k}\leq 0 \quad \text{and}\quad \bar{u}_{k}(x)\leq Ce^{-\frac{1}{2}R}\,\ \text{for}\,\ |x|\geq R,$$
where $C>0$ is independent of $k$.
By the comparison principle, comparing $\bar{u}_{k}$ with $Ce^{-\frac{1}{2}|x|}$ yields that
$$ \bar{u}_{k}(x)\leq Ce^{-\frac{1}{2}|x|}\,\ \text{for}\,\ |x|\geq R,$$
which then implies that \eqref{2:conexp} holds true.
Furthermore, 
one can deduce from (\ref{eq1.39}) that
$$ \big|\varepsilon_k^2V_1\big(\varepsilon_kx+ x_k \big)\bar{u}_{k}\big|\leq C' e^{-\frac{1}{4}|x|} \,\  \text{for}\,\   |x|>R, $$
where $C'>0$ is independent of $k$.
Therefore, by the exponential decay \eqref{2:conexp}, applying the local elliptic estimates (cf. $(3.15)$ in \cite{GT}) to \eqref{4:star} yields that
\begin{equation*}
  |\nabla \bar{u}_{k}(x)|\leq Ce^{-\frac{1}{4}|x|}\,\  \text{for}\,\   |x|>R.
\end{equation*}
The above estimates thus imply that the exponential decay \eqref{2:conexp} and \eqref{2:conexp2} hold true.

Since $w$ and $(\bar u_k, \bar v_k)$  decay exponentially as $|x|\to\infty$, using the standard elliptic regularity theory then yields that \eqref{lem4.1:1} holds uniformly in $\R^2$ (e.g. Lemma 4.9 in \cite{M} for similar arguments).
\qed
\vskip .1truein

We now note from Lemma \ref{lem4.1} that $ (u_0 ,v_0 )=\big(\sqrt{\gamma} w , \sqrt{1-\gamma} w \big)$ is a positive solution of the following system
\begin{equation}\label{uniq:limit-1}
\arraycolsep=1.5pt
\left\{\begin{array}{lll}
\Delta u -  u +\displaystyle\frac{a_1}{a^*} u ^3 + \displaystyle\frac{\beta ^*}{a^*}  v ^2 u &=0   \,\ \mbox{in}\,\  \R^2,\\ [2.0mm]
 \Delta v - v +\displaystyle\frac{a_2}{a^*} v ^3 +  \displaystyle\frac{\beta ^*}{a^*} u ^2 v &=0   \,\ \mbox{in}\,\  \R^2,\,\
\end{array}\right.
\end{equation}
It follows from \cite[Lemma 2.2 and Theorem 3.1]{DW} that the positive solution $(u_0,v_0)$ is {\em non-degenerate}, in the sense that the solution set of the linearized system for (\ref{uniq:limit-1}) about $(u_0,v_0)$ satisfying
\bsub \label{uniq:limit-2}
\arraycolsep=1.5pt
\begin{eqnarray}
\mathcal{L}_1(\phi_1,\phi_2):=\Delta \phi_1-  \phi_1 + \displaystyle\frac{3a_1}{a^*} u_0^2\phi_1 + \displaystyle\frac{\beta ^*}{a^*} v_0^2 \phi_1+\displaystyle\frac{2\beta ^*}{a^*}  u_0 v_0\phi_2&=0   \,\ \mbox{in}\,\  \R^2, \label{uniq:limit-2a}\\ [1.0mm]
 \mathcal{L}_2(\phi_2,\phi_1):=\Delta \phi_2-  \phi_2 +\displaystyle\frac{3a_2}{a^*} v_0^2\phi_2 +  \displaystyle\frac{\beta ^*}{a^*} u_0^2 \phi_2+\displaystyle\frac{2\beta ^*}{a^*}  u_0 v_0\phi_1\, &=0   \,\ \mbox{in}\,\  \R^2,\,\ \label{uniq:limit-2b}
 \end{eqnarray}\esub
is exactly $2$-dimensional, namely,
\begin{equation}\label{uniq:limit-3}
\left(\begin{array}{cc} \phi_1\\[2mm]
 \phi_2\end{array}\right)=\sum _{j=1}^{2}b_j\left(\begin{array}{cc} \frac{\partial u_0}{\partial x_j}\\[2mm] \frac{\partial v_0}{\partial x_j}
\end{array}\right)
\end{equation}
for some constants $b_j$.

Under the assumptions of Lemma \ref{lem4.1}, the nonnegative minimizer $(u_k,v_k)$ of $e(a_1,a_2,\beta_k)$ satisfies
\begin{equation}\label{uniq:a-1H}
\arraycolsep=1.5pt
\left\{\begin{array}{lll}
-\Delta u_k +V_1(x)u_k =\mu_{ k} u_k +a_1u_k^3 +\beta_k  v_k^2 u_k   \,\ \mbox{in}\,\  \R^2,\\[2mm]
-\Delta v_k +V_2(x)v_k  =\mu_{k} v_k +a_2v_k^3 +\beta_k  u_k^2 v_k \,  \ \ \mbox{in}\,\  \R^2,\,\
\end{array}\right.
\end{equation}
where $\mu _{k}\in \R$ is a suitable Lagrange multiplier and satisfies
\begin{equation}\label{uniq:a-6H}
 \mu_k=e(a_1,a_2,\beta_k)- \displaystyle \frac{a_1}{2} \inte   u_k^4 - \displaystyle \frac{a_2}{2 }\inte    v_k^4-\displaystyle \beta _k \inte u_k^2v_k^2.
\end{equation}
One can further check from (\ref{uniq:a-6H}) and Lemma \ref{lem4.1} that $\mu _{k}$ satisfies
\begin{equation}\label{uniq:a-6HH}
\lim_{k\to\infty} \mu_k\varepsilon^2_{k}=- 1,
\end{equation}
where $\varepsilon _{k}>0$ is defined by (\ref{lem4.1:2}).

\subsection{Proof of Theorem \ref{Thm1.5}}

\vskip 0.1truein

This subsection is devoted to the proof of Theorem \ref{Thm1.5} on the uniqueness of nonnegative minimizers.
For any given $a_1\in (0, a^*)$ and $a_2\in (0, a^*)$, towards this purpose we suppose that there exist two different nonnegative minimizers $(u_{1,k},v_{1,k})$ and $(u_{2,k},v_{2,k})$ of $e(a_1,a_2,\beta_k)$, where   $\beta_k\nearrow\beta^*$ as $k\to\infty$. Let $(x_{1,k},y_{1,k})$ and $(x_{2,k},y_{2,k})$ be the unique maximum point of $(u_{1,k},v_{1,k})$ and $(u_{2,k},v_{2,k})$, respectively. Note from (\ref{uniq:a-1H}) that the nonnegative minimizer $(u_{i,k},v_{i,k})$ solves the system
\begin{equation*}
 \begin{cases}
-\Delta u_{i,k} +V_1(x)u_{i,k} =\mu_{i,k} u_{i,k} +a_1u_{i,k}^3 +\beta_k  v_{i,k}^2 u_{i,k}   \,\ \mbox{in}\,\  \R^2,\\
-\Delta v_{i,k} +V_2(x)v_{i,k}\, =\mu_{i,k} v_{i,k} +a_2v_{i,k}^3 +\beta_k  u_{i,k}^2 v_{i,k}   \ \ \mbox{in}\,\  \R^2,\,\
\end{cases}
\end{equation*}
where $\mu _{i,k}\in \R$ is a suitable Lagrange multiplier and satisfies (\ref{uniq:a-6H}) and (\ref{uniq:a-6HH}) with $\mu _k=\mu _{i,k}$ for $i=1, 2$.
Define
\begin{equation}\label{uniq:a-2}
\big(\bar u_{i,k}(x), \bar v_{i,k}(x)\big):= \sqrt{a^*}\eps _k \big(  u_{i,k}(\varepsilon_{k} x+x_{2,k}),   v_{i,k}(\varepsilon_{k} x+x_{2,k})\big), \ \ \mbox{where}\ \ i=1,2.
\end{equation}
It then follows from Lemma \ref{lem4.1} that
\begin{equation}\label{uniq:a-3}
\big(\bar u_{i,k}(x), \bar v_{i,k}(x)\big)\to \big(u_0,v_0\big)\equiv \big(\sqrt{\gamma} w , \sqrt{1-\gamma} w \big)
\end{equation}
uniformly in $\R^2$ as $k\to\infty$, and  $\big(\bar u_{i,k}(x), \bar v_{i,k}(x)\big)$ satisfies the system
\begin{equation}\label{uniq:a-4}\arraycolsep=1.5pt
\left\{\begin{array}{lll}
- \Delta \bar u_{i,k} +\varepsilon_{k}^2V_1(\varepsilon_{k} x+x_{2,k})\bar u_{i,k} &=\mu_{i,k} \varepsilon_{k}^2\bar u_{i,k} +\displaystyle\frac{a_1}{a^*}\bar u_{i,k}^3 +\displaystyle\frac{\beta_k }{a^*} \bar v_{i,k}^2 \bar u_{i,k}   \,\ \mbox{in}\,\  \R^2,\\[2.5mm]
- \Delta \bar v_{i,k} +\varepsilon_{k}^2V_2(\varepsilon_{k} x+x_{2,k})\bar v_{i,k}  &=\mu_{i,k} \varepsilon_{k}^2\bar v_{i,k} +\displaystyle\frac{a_2}{a^*}\bar v_{i,k}^3 +\displaystyle\frac{\beta_k }{a^*} \bar u_{i,k}^2 \bar v_{i,k} \,  \,\ \mbox{in}\,\  \R^2.\,\
\end{array}\right.
\end{equation}

\begin{lem}\label{lem4.2} For any given $a_1\in (0, a^*)$ and $a_2\in (0, a^*)$, we have
\begin{equation}\label{step-1:1}
C_1\|\bar v_{2,k}-\bar v_{1,k}\| _{L^\infty(\R^2)}\le \|\bar u_{2,k}-\bar u_{1,k}\| _{L^\infty(\R^2)}\le C_2 \|\bar v_{2,k}-\bar v_{1,k}\| _{L^\infty(\R^2)}  \ \ \mbox{as}\ \ k\to\infty,
\end{equation}
where the  positive constants $C_1>0$ and $C_2>0$ are independent of $k$.
\end{lem}

\noindent\textbf{Proof.}
We first prove that the right inequality of (\ref{step-1:1}) holds. On the contrary, suppose that
\begin{equation}\label{step-1:2}
 \liminf _{k\to\infty}\frac{\|\bar u_{2,k}-\bar u_{1,k}\| _{L^\infty(\R^2)}}{\|\bar v_{2,k}-\bar v_{1,k}\| _{L^\infty(\R^2)} }=+\infty \ \ \mbox{as}\ \ k\to\infty.
\end{equation}
Following the first equation of (\ref{uniq:a-4}), we have
\begin{equation}\label{uniq:a-5}\arraycolsep=1.5pt
 \begin{array}{lll}
&- \Delta (\bar u_{2,k}-\bar u_{1,k}) +\varepsilon_{k}^2V_1(\varepsilon_{k} x+x_{2,k})(\bar u_{2,k}-\bar u_{1,k})\\[2mm]
 =&\mu_{2,k} \varepsilon_{k}^2(\bar u_{2,k}-\bar u_{1,k})+ \bar u_{1,k}\varepsilon_{k}^2(\mu_{2,k}-\mu_{1,k})+\displaystyle\frac{a_1}{a^*}(\bar u_{2,k}^3-\bar u_{1,k}^3) \\[2mm]
 & +\displaystyle\frac{\beta_k }{a^*}\Big[\bar  v_{1,k}^2( \bar u_{2,k} -\bar u_{1,k})+\bar u_{2,k}( \bar v_{2,k}^2- \bar v_{1,k}^2 )\Big]  \,\ \mbox{in}\,\  \R^2,
\end{array}
\end{equation}
where the term $\varepsilon_{k}^2(\mu_{2,k}-\mu_{1,k})$ satisfies
\begin{equation}\label{uniq:a-6}\arraycolsep=1.5pt
 \begin{array}{lll}
\varepsilon_{k}^2(\mu_{2,k}-\mu_{1,k})&=- \displaystyle \frac{a_1}{2}\varepsilon_{k}^2\inte (  u_{2,k}^4- u_{1,k}^4)- \displaystyle \frac{a_2}{2 }\varepsilon_{k}^2\inte (  v_{2,k}^4- v_{1,k}^4)\\[3mm]
&\quad-\displaystyle \beta _k\varepsilon_{k}^2\inte \Big[v_{2,k}^2(  u_{2,k}^2- u_{1,k}^2)+u_{1,k}^2(  v_{2,k}^2- v_{1,k}^2)\Big]\\[3mm]
&=- \displaystyle \frac{a_1}{2(a^*)^2}\inte ( \bar u_{2,k}^4- \bar u_{1,k}^4)- \displaystyle \frac{a_2}{2(a^*)^2}\inte (  \bar v_{2,k}^4- \bar v_{1,k}^4)\\[3mm]
\quad &\quad-\displaystyle\frac{\beta_k}{(a^*)^2}\inte \Big[\bar v_{2,k}^2( \bar  u_{2,k}^2- \bar u_{1,k}^2)+\bar u_{1,k}^2( \bar  v_{2,k}^2- \bar v_{1,k}^2)\Big],
\end{array}
\end{equation}
by applying  (\ref{uniq:a-6H}) and (\ref{uniq:a-2}).
Set
\begin{equation}\label{step-1:3}
\xi_{k}(x)=\displaystyle\frac{\bar  u_{2,k}(x)- \bar  u_{1,k}(x)}{\|\bar u_{2,k}-\bar u_{1,k}\|_{L^\infty(\R^2)}}.
\end{equation}
It then yields from (\ref{uniq:a-5}) and (\ref{step-1:3}) that $\xi_{k}$ satisfies
\begin{equation}\label{uniq:a-6c}\arraycolsep=1.5pt
 \begin{array}{lll}
&- \Delta \xi_{k} +\varepsilon_{k}^2V_1(\varepsilon_{k} x+x_{2,k})\xi_{k}=\mu_{2,k} \varepsilon_{k}^2\xi_{k}+\displaystyle\frac{a_1}{a^*}\big(\bar u_{2,k}^2+\bar u_{2,k}\bar u_{1,k}+\bar u_{1,k}^2\big)\xi_{k} \\[2mm]
 &+\displaystyle \bar u_{1,k}\frac{\varepsilon_{k}^2(\mu_{2,k}-\mu_{1,k})}{\|\bar u_{2,k}-\bar u_{1,k}\|_{L^\infty(\R^2)}} +\displaystyle\frac{\beta_k }{a^*}\Big[\bar  v_{1,k}^2\xi_{k}+\bar u_{2,k}\frac{ \bar v_{2,k}^2- \bar v_{1,k}^2 }{\|\bar u_{2,k}-\bar u_{1,k}\|_{L^\infty(\R^2)}}\Big]  \,\ \mbox{in}\,\  \R^2.
\end{array}
\end{equation}
Since $\|\xi _k\|_{L^\infty(\R^2)}\le 1$, the standard elliptic regularity theory then implies from (\ref{uniq:a-6c}) and (\ref{step-1:2}) that $\|\xi _k\|_{C^{1,\alpha }_{loc}(\R^2)}\le C$ for some $\alp \in (0,1)$, where the constant $C>0$ is independent of $k$. Therefore, there exists a subsequence, still denoted by $\{\beta_k\}$, of $\{\beta_k\}$ and a function $\xi =\xi (x)$ such that $\xi_k \to \xi $ in $C_{loc}^1(\R^2)$ as $k\to\infty$.

Applying (\ref{step-1:2}) again, we now derive from (\ref{uniq:a-6HH}) and (\ref{uniq:a-6c}) that $\xi $ satisfies
\begin{equation}\label{step-1:4}
 \Delta \xi -  \xi+\Big(\displaystyle\frac{3a_1}{a^*}u_0^2+\displaystyle\frac{\beta ^*}{a^*}v_0^2\Big)\xi-\Big[\frac{2a_1}{(a^*)^2}\inte u_0^3\xi+\frac{2\beta ^*}{(a^*)^2}\inte u_0v_0^2\xi\Big]u_0=0 \,\ \mbox{in}\,\  \R^2,
\end{equation}
where the last term follows from (\ref{uniq:a-6}).
Employing (\ref{uniq:a-3}), we can simplify (\ref{step-1:4}) as
\[
 \Delta \xi -  \xi+\Big[\frac{3a_1}{a^*} \gamma + \displaystyle\frac{\beta ^*}{a^*}\big(1-\gamma\big) \Big]w^2\xi=\frac{2}{a^*}\gamma \Big[\frac{a_1}{a^*} \gamma +\displaystyle\frac{\beta ^*}{a^*}\big(1-\gamma\big) \Big]w \inte w^3\xi \,\,\,\ \mbox{in}\,\  \R^2,
\]
which can be further simplified as
\begin{equation}\label{step-1:5}
 \Delta \xi -  \xi+\Big[1+\frac{2a_1}{a^*} \gamma \Big]w^2\xi=\frac{2}{a^*}\gamma w \inte w^3\xi \,\,\,\ \mbox{in}\,\  \R^2.
\end{equation}
Multiplying (\ref{equ:w}) by $\xi$ yields that
$$
\inte \nabla w\nabla \xi +\inte  w \xi =\inte w^3\xi,
$$
and while multiplying (\ref{step-1:5}) by $w$ gives that
\[
\inte \nabla w\nabla \xi +\inte  w \xi = \Big[1+2\gamma \Big(\displaystyle\frac{a_1}{a^*} -1\Big) \Big]\inte w^3\xi.
\]
Since $a_1\in (0, a^*)$, the above two equations imply that $\inte w^3\xi =0$.
The equation (\ref{step-1:5}) is then reduced into
\begin{equation}\label{step-1:7}
 \Delta \xi -  \xi+\Big[1+\frac{2a_1}{a^*} \gamma \Big]w^2\xi=0 \,\ \mbox{in}\,\  \R^2,\quad \text{where}\ \ \Big[1+\frac{2a_1}{a^*} \gamma \Big] \in (1,3).
\end{equation}
We thus conclude from (\ref{step-1:7}) and \cite[Lemma 2.2]{DW} that
$\xi \equiv 0$ in $\R^2$.

On the other hand, let $y_k$ be a point satisfying $|\xi_k(y_k)|=\|\xi_k\|_{L^\infty(\R^2)}=1$. Since both $\bar u_{i,k}$ and $\bar v_{i,k}$ decay exponentially as $|x|\to\infty$, where $i=1, 2$, applying the maximum principle to (\ref{uniq:a-6c}) yields that $|y_k|\le C$ uniformly in $k$. Therefore, we conclude that $\xi_k\to \xi \not\equiv 0$ uniformly on $\R^2$, which however contradicts to the above conclusion that $\xi  \equiv 0$ on $\R^2$. This implies that (\ref{step-1:2}) is false and hence the right inequality of (\ref{step-1:1}) holds true.

By considering the second equation of (\ref{uniq:a-4}), the same argument as above yields that the left inequality of (\ref{step-1:1}) also holds. This completes the proof of Lemma \ref{lem4.2}.
\qed
\vskip .1truein

If  $v_{2,k}\equiv v_{1,k}$ (or $u_{2,k}\equiv u_{1,k}$) in $\R^2$, it then follows from Lemma \ref{lem4.2} that $u_{2,k}\equiv u_{1,k}$ (or $v_{2,k}\equiv v_{1,k}$, respectively) in $\R^2$, and the proof of Theorem \ref{Thm1.5} is thus done. Therefore, the rest is to consider the case where $v_{2,k}\not \equiv v_{1,k}$ and $u_{2,k}\not \equiv u_{1,k}$ in $\R^2$.

We define
\begin{equation*}
\big(\hat u_{i,k}(x), \hat v_{i,k}(x)\big):= \sqrt{a^*}\eps _k \big(  u_{i,k}( x),   v_{i,k}(x)\big), \ \ \mbox{where}\ \ i=1,2.
\end{equation*}
It then follows from Lemma \ref{lem4.1} that
\begin{equation*}
\big(\bar u_{i,k}(x), \bar v_{i,k}(x)\big)=\big(\hat u_{i,k}(\varepsilon_{k} x+x_{2,k}), \hat v_{i,k}(\varepsilon_{k} x+x_{2,k})\big)\to \big(u_0,v_0\big)\equiv \big(\sqrt{\gamma} w , \sqrt{1-\gamma} w \big)
\end{equation*}
uniformly in $\R^2$ as $k\to\infty$.
Note that $\big(\hat u_{i,k}(x), \hat v_{i,k}(x)\big)$ satisfies the system
\begin{equation}\label{5.2:0}\arraycolsep=1.5pt
\left\{\begin{array}{lll}
- \varepsilon_{k}^2\Delta \hat u_{i,k} +\varepsilon_{k}^2V_1(x)\hat u_{i,k} &=\mu_{i,k} \varepsilon_{k}^2\hat u_{i,k} +\displaystyle\frac{a_1}{a^*}\hat u_{i,k}^3 +\displaystyle\frac{\beta_k }{a^*} \hat v_{i,k}^2 \hat u_{i,k}   \,\ \mbox{in}\,\  \R^2,\\[2.5mm]
- \varepsilon_{k}^2\Delta \hat v_{i,k} +\varepsilon_{k}^2V_2(x)\hat v_{i,k}  &=\mu_{i,k} \varepsilon_{k}^2\hat v_{i,k} +\displaystyle\frac{a_2}{a^*}\hat v_{i,k}^3 +\displaystyle\frac{\beta_k }{a^*} \hat u_{i,k}^2 \hat v_{i,k} \,  \,\ \mbox{in}\,\  \R^2.\,\
\end{array}\right.
\end{equation}
Since $v_{2,k}\not \equiv v_{1,k}$ and $u_{2,k}\not \equiv u_{1,k}$ in $\R^2$, we also define
\begin{equation*}
\arraycolsep=1.5pt
 \begin{array}{lll}
  \hat\xi_{1,k}(x)&=\displaystyle\frac{ u_{2,k}(x)-  u_{1,k}(x)}{\| u_{2,k}-  u_{1,k}\|^\frac{1}{2}_{L^\infty(\R^2)}\| v_{2,k}- v_{1,k}\|^\frac{1}{2}_{L^\infty(\R^2)}},\\[4mm]
 \hat\xi_{2,k}(x)&=\displaystyle\frac{ v_{2,k}(x)-  v_{1,k}(x)}{\| u_{2,k}-  u_{1,k}\|^\frac{1}{2}_{L^\infty(\R^2)}\| v_{2,k}-  v_{1,k}\|^\frac{1}{2}_{L^\infty(\R^2)}}.
\end{array}
\end{equation*}
From Lemma \ref{lem4.2}, we deduce  that there exists $C>0$ independent of $k$ such that
\begin{equation}\label{eq4.29}
0\le|\hat \xi_{1,k}(x)|,|\hat \xi_{2,k}(x)|\leq C<\infty \text{ and }|\hat \xi_{1,k}(x) \hat \xi_{2,k}(x)|\le 1 \text{ in } \R^2.
\end{equation}
Moreover, we have the following local estimates of $(\hat\xi_{1,k}, \hat\xi_{2,k})$.

\begin{lem}\label{lem4.3} Assume that $a_1\in (0, a^*)$ and $a_2\in (0, a^*)$ are given. Then for any $x_0\in\R^2$, there exists a small constant $\delta >0$  such that
\begin{equation}
    \int_{\partial B_\delta (x_0)} \Big( \eps ^2_k |\nabla \hat \xi_{i,k}|^2+ \frac{1}{2} |\hat \xi_{i,k}|^2+ \eps ^2_k  V_i(x)|\hat \xi_{i,k}|^2\Big)dS=O( \eps ^2_k)\,\ \text{as}\,\ k\to\infty,\,\   i=1, 2.
\label{5.2:6}
\end{equation}
\end{lem}

The proof of Lemma \ref{lem4.3} is given in Appendix A.3. Following above estimates, we shall carry out the rest proof of Theorem \ref{Thm1.5} by the following three steps:

\vskip 0.1truein

\noindent{\em  Step 1.}
Set
\begin{equation}\label{uniq:B-3}
\xi_{i,k}(x)=\hat\xi_{i,k}(\eps_kx+x_{2,k}),\, \ \text{where} \, \ i=1, 2.
\end{equation}
If $a_1\in (0, a^*)$ and $a_2\in (0, a^*)$, then there exists a subsequence (still denoted by $\{\beta_k\}$) of $\{\beta_k\}$ such that
\[(\xi_{1,k}, \xi_{2,k})\to (\xi_{10}, \xi_{20}) \,\ \text{in}\,\  C_{loc}^1(\R^2)  \,\ \text{in}\,\   k\to\infty,\]
where $(\xi_{10}, \xi_{20})$ satisfies
\begin{equation}\label{uniq:limit-A3}
\left(\begin{array}{cc}   \xi_{10}\\[2mm]
\xi_{20}\end{array}\right)=\displaystyle b_0 \left(\begin{array}{cc}   u_0+x\cdot\nabla u_0 \\[2mm]
 v_0+x\cdot\nabla v_0 \end{array}\right)+\sum _{j=1}^{2}b_j\left(\begin{array}{cc} \frac{\partial u_0}{\partial x_j}\\[2mm] \frac{\partial v_0}{\partial x_j}
\end{array}\right)
\end{equation}
for some constants $b_j$ with $j=0, 1, 2$.

Following (\ref{step-11:9}), one can actually check that $(\xi_{1,k}, \xi_{2,k})$ satisfies
\begin{equation}\label{step-1:9}
\arraycolsep=1.5pt
 \left\{\begin{array}{lll}
&  \Delta \xi_{1,k}-\varepsilon_{k}^2V_1(\varepsilon_{k} x+x_{2,k})\xi_{1,k}+\mu_{2,k} \varepsilon_{k}^2\xi_{1,k}+\displaystyle\frac{a_1}{a^*}\big(\bar u_{2,k}^2+\bar u_{2,k}\bar u_{1,k}+\bar u_{1,k}^2\big)\xi_{1,k} \\[2mm]
 &\qquad +\displaystyle\frac{\beta_k }{a^*}\big[\bar  v_{1,k}^2\xi_{1,k}+\bar u_{2,k}( \bar v_{2,k}+\bar v_{1,k})\xi_{2,k}\big] +\displaystyle c_k\bar u_{1,k}=0 \,\,\,\ \mbox{in}\,\  \R^2,\\[4mm]
 &  \Delta \xi_{2,k}-\varepsilon_{k}^2V_2(\varepsilon_{k} x+x_{2,k})\xi_{2,k}+\mu_{2,k} \varepsilon_{k}^2\xi_{2,k}+\displaystyle\frac{a_2}{a^*}\big(\bar v_{2,k}^2+\bar v_{2,k}\bar v_{1,k}+\bar v_{1,k}^2\big)\xi_{2,k} \\[2mm]
 &\qquad +\displaystyle\frac{\beta_k }{a^*}\big[\bar  u_{1,k}^2\xi_{2,k}+\bar v_{2,k}( \bar u_{2,k}+\bar u_{1,k})\xi_{1,k}\big] +\displaystyle c_k\bar v_{1,k}=0 \,\,\,\ \mbox{in}\,\  \R^2,
\end{array}\right.
\end{equation}
where the coefficient $c_{k}$ satisfies
\begin{equation}\label{step-1:10}
\arraycolsep=1.5pt
 \begin{array}{lll}
c_k:&=\displaystyle\frac{\varepsilon_{k}^2(\mu_{2,k}-\mu_{1,k})}{\|\bar u_{2,k}-\bar u_{1,k}\|^\frac{1}{2}_{L^\infty(\R^2)}\|\bar v_{2,k}-\bar v_{1,k}\|^\frac{1}{2}_{L^\infty(\R^2)}}\\[3mm]
&=\displaystyle\frac{\varepsilon_{k}^2(\mu_{2,k}-\mu_{1,k})}{\|\hat u_{2,k}-\hat u_{1,k}\|^\frac{1}{2}_{L^\infty(\R^2)}\|\hat  v_{2,k}-\hat v_{1,k}\|^\frac{1}{2}_{L^\infty(\R^2)}}\\[3mm]
&=- \displaystyle \frac{a_1}{2(a^*)^2}\inte ( \bar u_{2,k}^2+ \bar u_{1,k}^2)( \bar u_{2,k}+\bar u_{1,k})\xi_{1,k}\\[3mm]
&\quad - \displaystyle \frac{a_2}{2(a^*)^2}\inte ( \bar v_{2,k}^2+ \bar v_{1,k}^2)( \bar v_{2,k}+\bar v_{1,k})\xi_{2,k}\\[3mm]
\quad &\quad-\displaystyle\frac{\beta_k}{(a^*)^2}\inte \Big[\bar v_{2,k}^2( \bar  u_{2,k}+ \bar u_{1,k})\xi_{1,k}+\bar u_{1,k}^2( \bar  v_{2,k}+ \bar v_{1,k})\xi_{2,k}\Big].
\end{array}\end{equation}
From (\ref{eq4.29}) and (\ref{uniq:B-3}) we see that
\begin{equation}\label{eq4.35}
\text{$\xi_{1,k}(x)$ and $\xi_{2,k}(x)$ are both bounded uniformly in $\R^2$},
 \end{equation}
  and thus there exists $C>0$ independent of $k$ such that
  \begin{equation}\label{eq4.36}
  |c_k|\leq C <\infty.
  \end{equation}
The standard elliptic regularity theory then implies from (\ref{step-1:10}) that $\|\xi _{1,k}\|_{C^{1,\alpha }_{loc}(\R^2)}\le C$ for some $\alp \in (0,1)$, where the constant $C>0$ is independent of $k$. Therefore, up to a subsequence if necessary, there holds that  $(\xi_{1,k}, \xi_{2,k})\overset{k}\to (\xi_{10}, \xi_{20})$ in $C_{loc}^1(\R^2)$, where   $(\xi_{10}, \xi_{20})$ satisfies
\begin{equation}\label{uniq:limit-A1}
\arraycolsep=1.5pt
\left\{\begin{array}{lll}
\Delta   \xi_{10}-   \xi_{10} + \displaystyle\frac{3a_1}{a^*} u_0^2  \xi_{10} + \frac{\beta^*}{a^*} v_0^2   \xi_{10}+ \frac{2\beta^*}{a^*}  u_0 v_0  \xi_{20}&=2b_0u_0   \ \ \mbox{in}\,\  \R^2,\\ [1.8mm]
 \Delta  \xi_{20}-    \xi_{20} +\displaystyle\frac{3a_2}{a^*} v_0^2 \xi_{20} +  \frac{\beta^*}{a^*} u_0^2   \xi_{20}+ \frac{2\beta^*}{a^*}  u_0 v_0  \xi_{10}&=2b_0v_0   \ \ \mbox{in}\,\  \R^2,
\end{array}\right.
\end{equation}
by applying (\ref{step-1:9}).
Here the constant $b_0\in \R$ is given by the limit $ c_k\overset{k}\to -2b_0$. Following (\ref{uniq:limit-3}), one can derive from (\ref{uniq:limit-A1}) that $(\xi_{10}, \xi_{20})$ satisfies (\ref{uniq:limit-A3}) for some constants $b_j$ with $j=0, 1, 2$, and Step 1 is thus established.

\vskip 0.05truein

\noindent{\em  Step 2.} The constants $b_0=b_1=b_2=0$ in (\ref{uniq:limit-A3}), i.e., $\xi_{10}=\xi_{20}=0$.

We first claim the following Pohozaev-type identities: if $p_{11}=p_{21},$ then
\begin{equation}
\arraycolsep=1.5pt\begin{array}{lll}
&&b_0 \displaystyle\inte \Big[\frac{\partial V_1(x+y_0)}{\partial x_j}\big(x\cdot \nabla u_0^2\big)+\displaystyle \frac{\partial V_2(x+y_0)}{\partial x_j}\big(x\cdot \nabla v_0^2\big)\Big]\\[4mm]
&&-\displaystyle\sum ^2_{i=1}b_i\inte \Big[\frac{\partial ^2 V_1(x+y_0)}{\partial x_j\partial x_i}u_0^2+  \frac{\partial ^2 V_2(x+y_0)}{\partial x_j\partial x_i}v_0^2\Big]=0,\quad j=1,\,2,
\end{array}
\label{5.2:AA}
\end{equation}
and if   $p_{11}<p_{21},$ then
\begin{equation}
\arraycolsep=1.5pt\begin{array}{lll}
 b_0 \displaystyle\inte  \frac{\partial V_1(x+y_0)}{\partial x_j}\big(x\cdot \nabla u_0^2\big) -\displaystyle\sum ^2_{i=1}b_i\inte  \frac{\partial ^2 V_1(x+y_0)}{\partial x_j\partial x_i}u_0^2 =0,\quad j=1,\,2,
\end{array}
\label{5.2:AA-1}
\end{equation}
and if   $p_{11}>p_{21},$ then
\begin{equation}
\arraycolsep=1.5pt\begin{array}{lll}
 b_0 \displaystyle\inte  \frac{\partial V_2(x+y_0)}{\partial x_j}\big(x\cdot \nabla v_0^2\big) -\displaystyle\sum ^2_{i=1}b_i\inte  \frac{\partial ^2 V_2(x+y_0)}{\partial x_j\partial x_i}v_0^2 =0,\quad j=1,\,2.
\end{array}
\label{5.2:AA-2}
\end{equation}

To derive (\ref{5.2:AA})--(\ref{5.2:AA-2}), multiply the first equation of (\ref{5.2:0}) by $\frac{\partial \hat u_{i,k}}{\partial  x_j}$ and integrate over $B_\delta (x_{2,k})$, where $\delta >0$ is small and given by (\ref{5.2:6}). It then gives that
\begin{equation}\arraycolsep=1.5pt\begin{array}{lll}
&&-\eps _k^2\displaystyle\intB\frac{\partial \hat u_{i,k}}{\partial  x_j}\Delta \hat u_{i,k}+\eps _k^2\displaystyle\intB V_1(x)\frac{\partial \hat u_{i,k}}{\partial  x_j} \hat u_{i,k}\\[4mm]
&=&\mu_{i,k}\eps _k^2\displaystyle\intB \frac{\partial \hat u_{i,k}}{\partial  x_j} \hat u_{i,k}+\displaystyle\frac{a_1}{a^*}\intB \frac{\partial \hat u_{i,k}}{\partial  x_j} \hat u_{i,k}^3
+\displaystyle\frac{\beta _k}{a^*}\intB \frac{\partial \hat u_{i,k}}{\partial  x_j} \hat u_{i,k}\hat v_{i,k}^2\\[4mm]
&=&\displaystyle\frac{\mu_{i,k}}{2}\eps _k^2\intPB \hat u_{i,k}^2\nu _jdS+\displaystyle\frac{a_1}{4a^*}\intPB \hat u_{i,k}^4\nu _jdS+\displaystyle\frac{\beta _k}{2a^*}\intB \frac{\partial \hat u_{i,k}^2}{\partial  x_j}  \hat v_{i,k}^2,
\end{array}\label{5.2:7}
\end{equation}
where $\nu =(\nu _1,\nu _2)$ denotes the outward unit normal of $\partial B_\delta (x_{2,k})$.
Note that
\[\arraycolsep=1.5pt\begin{array}{lll}
 &&-\eps _k^2\displaystyle\intB\frac{\partial \hat u_{i,k}}{\partial  x_j}\Delta \hat u_{i,k}\\[4mm]&=&-\eps _k^2\displaystyle\intPB\frac{\partial \hat u_{i,k}}{\partial  x_j}\frac{\partial \hat u_{i,k}}{\partial  \nu}dS+\eps _k^2\displaystyle\intB\nabla \hat u_{i,k}\cdot\nabla\frac{\partial \hat u_{i,k}}{\partial  x_j}\\[4mm]
 &=&-\eps _k^2\displaystyle\intPB\frac{\partial \hat u_{i,k}}{\partial  x_j}\frac{\partial \hat u_{i,k}}{\partial  \nu}dS+\displaystyle\frac{1}{2}\eps _k^2\intPB |\nabla \hat u_{i,k}|^2\nu _jdS,
\end{array}
\]
and
\[
 \eps _k^2\displaystyle\intB V_1(x)\frac{\partial \hat u_{i,k}}{\partial  x_j} \hat u_{i,k}= \frac{\eps _k^2}{2}\intPB V_1(x)\hat u_{i,k}^2\nu _jdS-\frac{\eps _k^2}{2}\intB \frac{\partial V_1(x)}{\partial  x_j}\hat u_{i,k}^2.
\]
We then derive from (\ref{5.2:7}) that
\begin{equation}\arraycolsep=1.5pt\begin{array}{lll}
&&\displaystyle\eps _k^2 \intB \frac{\partial V_1(x)}{\partial  x_j}\hat u_{i,k}^2+\displaystyle\frac{\beta _k}{a^*}\intB \frac{\partial \hat u^2_{i,k}}{\partial  x_j} \hat v_{i,k}^2\\[4mm]
&=&-2\eps _k^2\displaystyle\intPB\frac{\partial \hat u_{i,k}}{\partial  x_j}\frac{\partial \hat u_{i,k}}{\partial  \nu}dS+\displaystyle \eps _k^2\intPB |\nabla \hat u_{i,k}|^2\nu _jdS \\[4mm]
&& +\displaystyle\eps _k^2 \intPB V_1(x)\hat u_{i,k}^2\nu _jdS-\displaystyle \mu_{i,k}\eps _k^2\intPB \hat u_{i,k}^2\nu _jdS\\[4mm]
&&-\displaystyle\frac{ a_1}{2a^*}\intPB \hat u_{i,k}^4\nu _jdS.
\end{array}\label{5.2:8A}
\end{equation}
Similarly, we derive from the second equation of (\ref{5.2:0}) that
\begin{equation}\arraycolsep=1.5pt\begin{array}{lll}
&&\displaystyle\eps _k^2 \intB \frac{\partial V_2(x)}{\partial  x_j}\hat v_{i,k}^2+\displaystyle\frac{\beta _k}{a^*}\intB \frac{\partial \hat v^2_{i,k}}{\partial  x_j} \hat u_{i,k}^2\\[4mm]
&=&-2\eps _k^2\displaystyle\intPB\frac{\partial \hat v_{i,k}}{\partial  x_j}\frac{\partial \hat v_{i,k}}{\partial  \nu}dS+\displaystyle \eps _k^2\intPB |\nabla \hat v_{i,k}|^2\nu _jdS \\[4mm]
&& +\displaystyle\eps _k^2 \intPB V_2(x)\hat v_{i,k}^2\nu _jdS-\displaystyle \mu_{i,k}\eps _k^2\intPB \hat v_{i,k}^2\nu _jdS\\[4mm]
&&-\displaystyle\frac{ a_2}{2a^*}\intPB \hat v_{i,k}^4\nu _jdS.
\end{array}\label{5.2:8B}
\end{equation}
Note that
\[\displaystyle\intB \frac{\partial \hat v^2_{i,k}}{\partial  x_j} \hat u_{i,k}^2dx=-\intB \frac{\partial \hat u^2_{i,k}}{\partial  x_j} \hat v_{i,k}^2dx+\intPB   \hat v^2_{i,k}  \hat u_{i,k}^2dS.\]
Following (\ref{5.2:8A}) and (\ref{5.2:8B}), we thus have
\begin{equation}
\displaystyle\eps _k^2 \intB \Big[\frac{\partial V_1(x)}{\partial  x_j}\big(\hat u_{2,k}+\hat u_{1,k}\big) \hat \xi _{1,k}+\frac{\partial V_2(x)}{\partial  x_j}\big(\hat v_{2,k}+\hat v_{1,k}\big) \hat \xi _{2,k}\Big] dx:=\mathcal{I}^u_k+\mathcal{I}^v_k,
\label{5.2:8}
\end{equation}
where we denote
\[\arraycolsep=1.5pt\begin{array}{lll}
\mathcal{I}^u_k&=&-2\displaystyle \eps _k^2\intPB \Big[\frac{\partial \hat u_{2,k}}{\partial  x_j}\frac{\partial \hat \xi_{1,k}}{\partial  \nu}+\frac{\partial \hat \xi_{1,k}}{\partial  x_j}\frac{\partial \hat u_{1,k}}{\partial  \nu}\Big]dS\\[4mm]
&&+\eps _k^2\displaystyle\intPB\nabla \hat \xi_{1,k} \cdot\nabla \big(\hat u_{2,k}+\hat u_{1,k}\big)\nu _jdS \\[4mm]
&&+\displaystyle \eps _k^2 \intPB V_1(x)\big(\hat u_{2,k}+\hat u_{1,k}\big) \hat \xi _{1,k} \nu _jdS\\[4mm]
&&-c_k\displaystyle\intPB \hat u_{2,k}^2\nu _jdS-\displaystyle \mu_{1,k}\eps _k^2\intPB \big(\hat u_{2,k}+\hat u_{1,k}\big) \hat \xi _{1,k}\nu _jdS \\[4mm]
&& -\displaystyle\frac{ a_1}{2a^*}\intPB \big(\hat u^2_{2,k}+\hat u^2_{1,k}\big) \big(\hat u_{2,k}+\hat u_{1,k}\big)\hat \xi _{1,k}\nu _jdS,
\end{array}\label{5.2:9A}
\]
and
\[\arraycolsep=1.5pt\begin{array}{lll}
\mathcal{I}^v_k&=&-2\displaystyle \eps _k^2\intPB \Big[\frac{\partial \hat v_{2,k}}{\partial  x_j}\frac{\partial \hat \xi_{2,k}}{\partial  \nu}+\frac{\partial \hat \xi_{2,k}}{\partial  x_j}\frac{\partial \hat v_{1,k}}{\partial  \nu}\Big]dS\\[4mm]
&&+\eps _k^2\displaystyle\intPB\nabla \hat \xi_{2,k} \cdot\nabla \big(\hat v_{2,k}+\hat v_{1,k}\big)\nu _jdS\\[4mm]
&&+\displaystyle \eps _k^2 \intPB V_2(x)\big(\hat v_{2,k}+\hat v_{1,k}\big) \hat \xi _{2,k} \nu _jdS \\[4mm]
&&-c_k\displaystyle\intPB \hat v_{2,k}^2\nu _jdS -\displaystyle \mu_{1,k}\eps _k^2\intPB \big(\hat v_{2,k}+\hat v_{1,k}\big) \hat \xi _{2,k}\nu _jdS \\[4mm]
&& -\displaystyle\frac{ a_2}{2a^*}\intPB \big(\hat v^2_{2,k}+\hat v^2_{1,k}\big) \big(\hat v_{2,k}+\hat v_{1,k}\big)\hat \xi _{2,k}\nu _jdS\\[4mm]
&& -\displaystyle\frac{\beta _k}{a^*}\intPB   \big[\hat u^2_{2,k}  (\hat v_{2,k}+\hat v_{1,k})\hat \xi _{2,k}+\hat v^2_{1,k}  (\hat u_{2,k}+\hat u_{1,k})\hat \xi _{1,k}\big]dS.
\end{array}\label{5.2:9B}
\]
Here $c_k$ is given by (\ref{step-1:10}),
which is bounded  uniformly in $k$.

For the right hand side of (\ref{5.2:8}),
applying Lemma \ref{lem4.3}, we deduce that if $\delta >0$ is small,
\[\arraycolsep=1.5pt\begin{array}{lll}
&&\displaystyle \eps _k^2\intPB \Big|\frac{\partial \hat u_{2,k}}{\partial  x_j}\frac{\partial \hat \xi_{1,k}}{\partial  \nu}\Big|dS\\[4mm]
&\le &\displaystyle \eps _k\Big(\intPB \Big|\frac{\partial \hat u_{2,k}}{\partial  x_j}\Big|^2dS\Big)^{\frac{1}{2}}\Big(\eps _k^2\intPB \Big|\frac{\partial \hat \xi_{1,k}}{\partial  \nu}\Big|^2dS\Big)^{\frac{1}{2}}\le C\eps _k^2e^{-\frac{C\delta}{\eps_k}}\,\ \mbox{as} \,\ k\to\infty,
\end{array}\label{5.2:9a}
\]
due to the fact that $\nabla \hat u_{2,k}(\eps _kx+x_{2,k})$ satisfies the exponential decay (\ref{2:conexp2}), where $C>0$ is independent of $k$. Similarly, employing (\ref{2:conexp}), (\ref{2:conexp2}), (\ref{eq4.36}) and Lemma \ref{lem4.3} again, we can  prove that other terms of $\mathcal{I}^u_k$ and $\mathcal{I}^v_k$ can also  be controlled  by the order $o(e^{-\frac{C\delta}{\eps_k}})$ as $k\to\infty$ for some $C>0$.
We therefore conclude from above that
\begin{equation}\label{5.2:9aB}
\mathcal{I}^u_k+\mathcal{I}^v_k=o(e^{-\frac{C\delta}{\eps_k}}) \,\ \mbox{as} \,\ k\to\infty.
\end{equation}

It now follows from (\ref{eq1.40}),  (\ref{uniq:B-3}), (\ref{5.2:8}) and (\ref{5.2:9aB}) that for small $\delta >0$,
\begin{align}\label{5.2:10}
 o(e^{-\frac{C\delta}{\eps_k}})&=\displaystyle\eps _k^2 \intB \Big[\frac{\partial V_1(x)}{\partial  x_j}\big(\hat u_{2,k}+\hat u_{1,k}\big) \hat \xi _{1,k}+\frac{\partial V_2(x)}{\partial  x_j}\big(\hat v_{2,k}+\hat v_{1,k}\big) \hat \xi _{2,k}\Big] dx\\
 &=\displaystyle\eps _k^4 \int_{B_{\frac{\delta}{\eps_k}}(0)} \Big[\frac{\partial V_1\big(\eps_k[x+(x_{2,k}-x_1)/\eps_k]+x_1\big)}{\partial  x_j}\big(\bar u_{2,k}+\bar u_{1,k}\big) \xi _{1,k}\nonumber\\
 &\quad+\frac{\partial V_2\big(\eps_k[x+(x_{2,k}-x_1)/\eps_k]+x_1\big)}{\partial  x_j}\big(\bar v_{2,k}+\bar v_{1,k}\big)  \xi _{2,k}\Big] dx\nonumber\\
 &=\eps _k^{4} \int_{B_{\frac{\delta}{\eps_k}}(0)} \bigg\{\eps_k^{p_{11}-1}\frac{\partial V_{11}\big(x+(x_{2,k}-x_1)/\eps_k\big)}{\partial  x_j}\big(\bar u_{2,k}+\bar u_{1,k}\big) \xi _{1,k}\nonumber\\
 &\quad+R_{1j}\big(\eps_k x+(x_{2,k}-x_1)\big)\big(\bar u_{2,k}+\bar u_{1,k}\big) \xi _{1,k}\nonumber\\
 &\quad+\eps_k^{p_{21}-1}\frac{\partial V_{21}\big(x+(x_{2,k}-x_1)/\eps_k\big)}{\partial  x_j}\big(\bar v_{2,k}+\bar v_{1,k}\big)  \xi _{2,k}\nonumber\\
 &\quad +R_{2j}\big(\eps_k x+(x_{2,k}-x_1)\big)\big(\bar v_{2,k}+\bar v_{1,k}\big) \xi _{2,k}\bigg\} dx.\label{eq4.56}
  \end{align}
  From (\ref{eq1.40}), (\ref{lem4.1:3}), (\ref{2:conexp}) and (\ref{eq4.35}) we note that
 \begin{equation*}
  \begin{split}
& \Big|\int_{B_{\frac{\delta}{\eps_k}}(0)}R_{1j}\big(\eps_k x+(x_{2,k}-x_1)\big)\big(\bar u_{2,k}+\bar u_{1,k}\big) \xi _{1,k} dx\Big|\\
&\leq C\eps _k^{q_1} \int_{B_{\frac{\delta}{\eps_k}}(0)} \big|x+\frac{x_{2,k}-x_1}{\eps_k}\big|^{q_1}\big(\bar u_{2,k}+\bar u_{1,k}\big) |\xi _{1,k}|dx\leq C\eps _k^{q_1}, \\
 \end{split}
  \end{equation*}
  and
 \begin{equation*}
  \begin{split}
& \Big|\int_{B_{\frac{\delta}{\eps_k}}(0)}R_{2j}\big(\eps_k x+(x_{2,k}-x_1)\big)\big(\bar v_{2,k}+\bar v_{1,k}\big) \xi _{2,k} dx\Big|\\
&\leq C\eps _k^{q_2} \int_{B_{\frac{\delta}{\eps_k}}(0)} \big|x+\frac{x_{2,k}-x_1}{\eps_k}\big|^{q_2}\big(\bar v_{2,k}+\bar v_{1,k}\big) |\xi _{1,k}|dx\leq C\eps _k^{q_2}, \\
 \end{split}
  \end{equation*}
 where $q_i>p_{i1}-1$ for $i=1,2$.
When $p_1=p_2$, it then follows from (\ref{eq4.56})  that
  \begin{equation*}
  \begin{split}
  o(1)&= \int_{B_{\frac{\delta}{\eps_k}}(0)} \Big[\frac{\partial V_{11}(x+(x_{2,k}-x_1)/\eps_k)}{\partial  x_j}\big(\bar u_{2,k}+\bar u_{1,k}\big) \xi _{1,k}\\
 &\quad\quad+\frac{\partial V_{21}(x+(x_{2,k}-x_1)/\eps_k)}{\partial  x_j}\big(\bar v_{2,k}+\bar v_{1,k}\big)  \xi _{2,k}\Big] dx.
 \end{split}
  \end{equation*}
On the other hand, we deduce from  \eqref{2:conexp} and (\ref{eq4.35}) that
 \begin{equation*}
  \begin{split}
  o(e^{-\frac{C\delta}{\eps_k}})&= \int_{B^c_{\frac{\delta}{\eps_k}}(0)} \Big[\frac{\partial V_{11}\big(x+(x_{2,k}-x_1)/\eps_k\big)}{\partial  x_j}\big(\bar u_{2,k}+\bar u_{1,k}\big) \xi _{1,k}\\
 &\quad\quad+\frac{\partial V_{21}\big(x+(x_{2,k}-x_1)/\eps_k\big)}{\partial  x_j}\big(\bar v_{2,k}+\bar v_{1,k}\big)  \xi _{2,k}\Big] dx.
 \end{split}
  \end{equation*}
Therefore, we have
   \begin{equation}\label{eq4.577}
  \begin{split}
  o(1)&= \inte \Big[\frac{\partial V_{11}\big(x+(x_{2,k}-x_1)/\eps_k\big)}{\partial  x_j}\big(\bar u_{2,k}+\bar u_{1,k}\big) \xi _{1,k}\\
 &\quad\quad+\frac{\partial V_{21}\big(x+(x_{2,k}-x_1)/\eps_k\big)}{\partial  x_j}\big(\bar v_{2,k}+\bar v_{1,k}\big)  \xi _{2,k}\Big] dx.
 \end{split}
  \end{equation}
Setting $k\to\infty$, we then derive from \eqref{lem4.1:3}, \eqref{uniq:a-3} and (\ref{eq4.577}) that
\[ \arraycolsep=1.5pt\begin{array}{lll}
0&=&2\displaystyle\inte \Big(\frac{\partial V_{11}(x+y_0)}{\partial x_j}u_0\,\xi_{10}+\frac{\partial V_{21}(x+y_0)}{\partial x_j}v_0\,\xi_{20}\Big)\\[3mm]
&=&\displaystyle2\inte \frac{\partial V_{11}(x+y_0)}{\partial x_j}u_0\Big[b_0\big(u_0+x\cdot \nabla u_0\big)+\sum ^2_{i=1}b_i\frac{\partial u_0}{\partial x_i}\Big]\\[3mm]
&&+\displaystyle2\inte \frac{\partial V_{21}(x+y_0)}{\partial x_j}v_0\Big[b_0\big(v_0+x\cdot \nabla v_0\big)+\sum ^2_{i=1}b_i\frac{\partial v_0}{\partial x_i}\Big]\\[4mm]
&=&b_0 \displaystyle\inte \Big[\frac{\partial V_{11}(x+y_0)}{\partial x_j}\big(x\cdot \nabla u_0^2\big)+\displaystyle \frac{\partial V_{21}(x+y_0)}{\partial x_j}\big(x\cdot \nabla v_0^2\big)\Big]\\[4mm]
&&-\displaystyle\sum ^2_{i=1}b_i\inte \Big[\frac{\partial ^2 V_{11}(x+y_0)}{\partial x_j\partial x_i}u_0^2+  \frac{\partial ^2 V_{21}(x+y_0)}{\partial x_j\partial x_i}v_0^2\Big],\quad j=1,\,2,
\end{array}\]
which then implies that (\ref{5.2:AA}) follows. Similarly, if $p_{11}\not =p_{21}$, then it also follows from (\ref{eq4.56}) that either (\ref{5.2:AA-1}) or (\ref{5.2:AA-2}) holds.

\vskip 0.1truein

We next prove that $b_0=0$ in (\ref{5.2:AA})--(\ref{5.2:AA-2}).
 Using the integration by parts,  we note that
\begin{equation}\arraycolsep=1.5pt\begin{array}{lll}
&&-\displaystyle\eps _k^2 \intB \big[(x-x_{2,k})\cdot \nabla \hat u_{i,k}\big] \Delta \hat u_{i,k} \\[4mm]
&=&- \eps _k^2\displaystyle \intPB \frac{\partial\hat u_{i,k} }{\partial \nu }(x-x_{2,k})\cdot \nabla \hat u_{i,k}
+\displaystyle\eps _k^2 \intB \nabla \hat u_{i,k}\nabla \big[(x-x_{2,k})\cdot \nabla \hat u_{i,k}\big] \\[4mm]
&=& - \eps _k^2\displaystyle \intPB \frac{\partial\hat u_{i,k} }{\partial \nu }(x-x_{2,k})\cdot \nabla \hat u_{i,k}
+\displaystyle \frac{\eps _k^2}{2}\intPB \big[(x-x_{2,k})\cdot \nu \big]|\nabla \hat u_{i,k}|^2.
\end{array}\label{5.3:3}
\end{equation}
Multiplying the first equation of (\ref{5.2:0}) by $ (x-x_{2,k})\cdot \nabla \hat u_{i,k} $, where $i=1,2$, and integrating over $B_\delta (x_{2,k})$, where $\delta >0$ is small as before, we deduce that for $i=1,2,$
\begin{equation}\arraycolsep=1.5pt\begin{array}{lll}
&&-\displaystyle\eps _k^2 \intB \big[(x-x_{2,k})\cdot \nabla \hat u_{i,k}\big] \Delta \hat u_{i,k} \\[4mm]
&=& \displaystyle\eps _k^2 \intB \big[\mu_{i,k}-V_1(x)\big] \hat u_{i,k}\big[(x-x_{2,k})\cdot \nabla \hat u_{i,k}\big]\\[4mm]
&& +
\displaystyle \frac{ a_1}{a^*} \intB \hat u_{i,k}^3\big[(x-x_{2,k})\cdot \nabla \hat u_{i,k}\big]
+\displaystyle \frac{\beta_k}{2a^*} \intB \hat v_{i,k}^2\big[(x-x_{2,k})\cdot \nabla \hat u^2_{i,k}\big]\\[4mm]
&=& -\mu_{i,k}\displaystyle\eps _k^2 \inte \hat u_{i,k}^2+\displaystyle\eps _k^2 \intB \Big(V_1(x)+\frac{x\cdot \nabla V_1(x)}{2}\Big)\hat u_{i,k}^2-\displaystyle \frac{  a_1}{2a^*} \inte \hat u_{i,k}^4\\[4mm]
&&-\displaystyle \frac{\beta_k}{a^*} \inte \hat u_{i,k}^2\hat v_{i,k}^2-\displaystyle \frac{\beta_k}{2a^*} \intB \hat u_{i,k}^2\big[(x-x_{2,k})\cdot \nabla \hat v^2_{i,k}\big]+I_i,
\end{array}\label{5.3:1}
\end{equation}
where the lower order term $I_i$ satisfies
\begin{equation}\arraycolsep=1.5pt\begin{array}{lll}
I_i&=& \mu_{i,k}\displaystyle\eps _k^2 \int _{\R^2\backslash B_\delta (x_{2,k})} \hat u_{i,k}^2+\displaystyle \frac{a_1}{2a^*} \int _{\R^2\backslash B_\delta (x_{2,k})} \hat u_{i,k}^4
\\[4mm]
&&-\displaystyle\frac{1}{2} \eps _k^2\intB \hat u_{i,k}^2 \big[x_{2,k}\cdot \nabla V_1(x)\big]dx\\[4mm]
&&+\displaystyle\frac{\eps _k^2}{2} \intPB\hat u_{i,k}^2\big[\mu_{i,k}-V_1(x)\big](x-x_{2,k})\nu dS\\[4mm]
&&+\displaystyle \frac{  a_1}{4a^*} \intPB \hat u_{i,k}^4(x-x_{2,k})\nu dS+\displaystyle \frac{\beta_k}{a^*} \int _{\R^2\backslash B_\delta (x_{2,k})} \hat u_{i,k}^2\hat v_{i,k}^2\\[4mm]
&&+\displaystyle \frac{\beta_k}{2a^*} \intPB \hat u_{i,k}^2\hat v_{i,k}^2 (x-x_{2,k})\nu dS,\,\ i=1,2.
\end{array}\label{5.3:2}
\end{equation}
Similarly, we have
\begin{equation}\arraycolsep=1.5pt\begin{array}{lll}
&&-\displaystyle\eps _k^2 \intB \big[(x-x_{2,k})\cdot \nabla \hat v_{i,k}\big] \Delta \hat v_{i,k} \\[4mm]
&=& - \eps _k^2\displaystyle \intPB \frac{\partial\hat v_{i,k} }{\partial \nu }(x-x_{2,k})\cdot \nabla \hat v_{i,k}
+\displaystyle \frac{\eps _k^2}{2}\intPB \big[(x-x_{2,k})\cdot \nu \big]|\nabla \hat v_{i,k}|^2,
\end{array}\label{5.3v:3}
\end{equation}
and  the second equation of (\ref{5.2:0}) yields that
\begin{equation}\arraycolsep=1.5pt\begin{array}{lll}
&&-\displaystyle\eps _k^2 \intB \big[(x-x_{2,k})\cdot \nabla \hat v_{i,k}\big] \Delta \hat v_{i,k} \\[4mm]
&=& -\mu_{i,k}\displaystyle\eps _k^2 \inte \hat v_{i,k}^2+\displaystyle\eps _k^2 \intB \Big(V_2(x)+\frac{x\cdot \nabla V_2(x)}{2}\Big)\hat v_{i,k}^2-\displaystyle \frac{  a_2}{2a^*} \inte \hat v_{i,k}^4\\[4mm]
&&+\displaystyle \frac{\beta_k}{2a^*} \intB \hat u_{i,k}^2\big[(x-x_{2,k})\cdot \nabla \hat v^2_{i,k}\big]+II_i,
\end{array}\label{5.3v:1}
\end{equation}
where the lower order term $II_i$ satisfies
\begin{equation}\arraycolsep=1.5pt\begin{array}{lll}
II_i&=& \mu_{i,k}\displaystyle\eps _k^2 \int _{\R^2\backslash B_\delta (x_{2,k})} \hat v_{i,k}^2+\displaystyle \frac{a_2}{2a^*} \int _{\R^2\backslash B_\delta (x_{2,k})} \hat v_{i,k}^4\\[4mm]
&&-\displaystyle\frac{1}{2} \eps _k^2\intB \hat v_{i,k}^2 \big[x_{2,k} \cdot \nabla V_2(x)\big]\\[4mm]
&&+\displaystyle\frac{\eps _k^2}{2} \intPB\hat v_{i,k}^2\big[\mu_{i,k}-V_2(x)\big](x-x_{2,k})\nu dS\\[4mm]
&&+\displaystyle \frac{ a_2}{4a^*} \intPB \hat v_{i,k}^4(x-x_{2,k})\nu dS,\,\ i=1,2.
\end{array}\label{5.3v:2}
\end{equation}

Since it follows from (\ref{uniq:a-6H}) that
\[\arraycolsep=1.5pt\begin{array}{lll}
&&\mu_{i,k}\eps _k^2\Big[\displaystyle \inte \hat u_{i,k}^2+\displaystyle \inte \hat v_{i,k}^2\Big]+\displaystyle \frac{  a_1}{2a^*} \inte \hat u_{i,k}^4+\displaystyle \frac{  a_2}{2a^*} \inte \hat v_{i,k}^4+\displaystyle \frac{\beta _k}{a^*} \inte \hat u_{i,k}^2\hat v_{i,k}^2\\[4mm]
&&= a^*\eps_k^4 e(a_1,a_2,\beta_k),
\end{array}\]
we reduce from (\ref{5.3:3})--(\ref{5.3v:2}) that
\[\arraycolsep=1.5pt\begin{array}{lll}
&&-\displaystyle\eps _k^2 \intB \Big[\Big(V_1(x)+\frac{x\cdot \nabla V_1(x)}{2}\Big)\hat u_{i,k}^2+ \Big(V_2(x)+\frac{x\cdot \nabla V_2(x)}{2}\Big)\hat v_{i,k}^2\Big]dx\\[4mm]
&&+a^*\eps_k^4 e(a_1,a_2,\beta_k)\\[1mm]
&=&I_i+II_i+\eps _k^2\displaystyle \intPB \frac{\partial\hat u_{i,k} }{\partial \nu }(x-x_{2,k})\cdot \nabla \hat u_{i,k}\\[3mm]
&&-\displaystyle \frac{\eps _k^2}{2}\intPB \big[(x-x_{2,k})\cdot \nu \big]|\nabla \hat u_{i,k}|^2+\eps _k^2\displaystyle \intPB \frac{\partial\hat v_{i,k} }{\partial \nu }(x-x_{2,k})\cdot \nabla \hat v_{i,k}\\[3mm]
&&-\displaystyle \frac{\eps _k^2}{2}\intPB \big[(x-x_{2,k})\cdot \nu \big]|\nabla \hat v_{i,k}|^2,\,\ i=1,2,
\end{array}\]
which implies that
\begin{equation}
\begin{split}
\displaystyle -\eps _k^{2}&\intB \Big[\Big(V_1(x)+\frac{x\cdot \nabla V_1(x)}{2}\Big)(\hat u_{1,k}+\hat u_{2,k})\hat\xi_{1,k}\\
+& \Big(V_2(x)+\frac{x\cdot \nabla V_2(x)}{2}\Big)(\hat v_{1,k}+\hat v_{2,k})\hat\xi_{2,k}\Big]dx
=T_k.
\end{split}
 \label{5.3:4}
\end{equation}
Here the term $T_k$ satisfies that for small $\delta >0$,
\begin{equation}\arraycolsep=1.5pt\begin{array}{lll}
T_k&=&\displaystyle\frac{(I_2-I_1)+(II_2-II_1)}{\|\hat  u_{2,k}-\hat  u_{1,k}\|^\frac{1}{2}_{L^\infty(\R^2)}\|\hat  v_{2,k}-\hat  v_{1,k}\|^\frac{1}{2}_{L^\infty(\R^2)}}\\[4mm]
&&-\displaystyle \frac{\eps _k^2}{2}\intPB \big[(x-x_{2,k})\cdot \nu \big]\big(\nabla \hat u_{2,k}+\nabla \hat u_{1,k}\big)\nabla\hat \xi_{1,k}\\[4mm]
&&+\displaystyle  \eps _k^2 \intPB \Big\{\big[(x-x_{2,k})\cdot \nabla \hat u_{2,k}\big]\big(\nu \cdot \nabla\hat \xi_{1,k}\big)+\big(\nu \cdot \nabla \hat u_{1,k}\big)\big[(x-x_{2,k})\cdot \nabla \hat \xi_{1,k}\big]\Big\}\\[4mm]
&&-\displaystyle \frac{\eps _k^2}{2}\intPB \big[(x-x_{2,k})\cdot \nu \big]\big(\nabla \hat v_{2,k}+\nabla \hat v_{1,k}\big)\nabla\hat \xi_{2,k}\\[4mm]
&&+\displaystyle  \eps _k^2 \intPB \Big\{\big[(x-x_{2,k})\cdot \nabla \hat v_{2,k}\big]\big(\nu \cdot \nabla\hat \xi_{2,k}\big)+\big(\nu \cdot \nabla \hat v_{1,k}\big)\big[(x-x_{2,k})\cdot \nabla \hat \xi_{2,k}\big]\Big\}
\\[4mm]
&=&\displaystyle\frac{(I_2-I_1)+(II_2-II_1)}{\|\hat  u_{2,k}-\hat  u_{1,k}\|^\frac{1}{2}_{L^\infty(\R^2)}\|\hat  v_{2,k}-\hat  v_{1,k}\|^\frac{1}{2}_{L^\infty(\R^2)}}+o(e^{-\frac{C\delta}{\eps_k}})\,\ \mbox{as} \,\ k\to\infty,
\end{array} \label{5.3:5}\end{equation}
due to (\ref{5.2:6}), where the second equality follows by applying the argument of estimating  (\ref{5.2:9aB}).

Using the arguments of estimating (\ref{5.2:9aB}) again, along with the exponential decay (\ref{2:conexp}) and (\ref{2:conexp2}), we also derive from (\ref{5.3:2}) that for small $\delta >0$,
\begin{equation}\arraycolsep=1.5pt\begin{array}{lll}
&&\displaystyle\frac{I_2-I_1}{\|\hat  u_{2,k}-\hat  u_{1,k}\|^\frac{1}{2}_{L^\infty(\R^2)}\|\hat  v_{2,k}-\hat  v_{1,k}\|^\frac{1}{2}_{L^\infty(\R^2)}}\\[4mm]
&=&\mu_{2,k}\displaystyle\eps _k^2 \int _{\R^2\backslash B_\delta (x_{2,k})} \big(\hat u_{2,k}+\hat u_{1,k}\big)\hat \xi_{1,k}
+\displaystyle \frac{  a_1}{2a^*} \int _{\R^2\backslash B_\delta (x_{2,k})} \big(\hat u_{2,k}^2+\hat u_{1,k}^2\big)\big(\hat u_{2,k}+\hat u_{1,k}\big)\hat \xi_{1,k}\\[4mm]
&&+\displaystyle \frac{  \beta_k}{a^*} \int _{\R^2\backslash B_\delta (x_{2,k})} \big[\hat v_{2,k}^2(\hat u_{2,k}+\hat u_{1,k})\hat\xi_{1,k}+\hat u_{1,k}^2 (\hat v_{2,k}+\hat v_{1,k}\big)\hat \xi_{2,k}\big]
\\[4mm]
&&+c_k \displaystyle\int _{\R^2\backslash B_\delta (x_{2,k})}  \hat u_{1,k}^2-\displaystyle\frac{1}{2} \eps _k^2\intB  \big[x_{2,k}\cdot \nabla V_1(x)\big]\big(\hat u_{2,k}+\hat u_{1,k}\big)\hat \xi_{1,k}\\[4mm]
&&+\displaystyle \frac{ a_1}{4a^*} \intPB   \big(\hat u_{2,k}^2+\hat u_{1,k}^2\big)\big(\hat u_{2,k}+\hat u_{1,k}\big)\hat \xi_{1,k}(x-x_{2,k})\nu dS\\[4mm]
&&-\displaystyle\frac{\eps _k^2}{2} \intPB \big(\hat u_{2,k}+\hat u_{1,k}\big)\hat \xi_{1,k}V_1(x)(x-x_{2,k})\nu dS\\[4mm]
&&+\displaystyle\frac{\mu_{2,k}\eps _k^2}{2} \intPB   \big(\hat u_{2,k}+\hat u_{1,k}\big)\hat \xi_{1,k}(x-x_{2,k})\nu dS
+c_k \displaystyle\intPB   \hat u_{1,k}^2 (x-x_{2,k})\nu dS\\[4mm]
&&+\displaystyle \frac{  \beta_k}{2a^*} \int _{\R^2\backslash B_\delta (x_{2,k})} \big[\hat v_{2,k}^2(\hat u_{2,k}+\hat u_{1,k})\hat\xi_{1,k}+\hat u_{1,k}^2 (\hat v_{2,k}+\hat v_{1,k}\big)\hat \xi_{2,k}\big](x-x_{2,k})\nu dS\\[4mm]
&=&
-\displaystyle\frac{1}{2} \eps _k^2\intB  \big[x_{2,k}\cdot \nabla V_1(x)\big]\big(\hat u_{2,k}+\hat u_{1,k}\big)\hat \xi_{1,k}+o(e^{-\frac{C\delta}{\eps_k}})\,\ \mbox{as} \,\ k\to\infty,
\end{array} \label{5.3:10}
\end{equation}
where $c_k$ is defined by (\ref{step-1:10}).
Similarly, we  can derive from \eqref{5.3v:2} that
\begin{equation}\label{eq4.65}
\begin{split}
\displaystyle&\frac{II_2-II_1}{\|\hat  u_{2,k}-\hat  u_{1,k}\|^\frac{1}{2}_{L^\infty(\R^2)}\|\hat  v_{2,k}-\hat  v_{1,k}\|^\frac{1}{2}_{L^\infty(\R^2)}} \\
&=-\displaystyle\frac{1}{2} \eps _k^2\intB  \big[x_{2,k}\cdot \nabla V_2(x)\big]\big(\hat v_{2,k}+\hat v_{1,k}\big)\hat \xi_{2,k}+o(e^{-\frac{C\delta}{\eps_k}})\,\ \mbox{as} \,\ k\to\infty.
\end{split}
\end{equation}
It then follows from \eqref{5.2:10}, \eqref{5.3:10} and \eqref{eq4.65}  that
$$\frac{(I_2-I_1)+(II_2-II_1)}{\|\hat  u_{2,k}-\hat  u_{1,k}\|^\frac{1}{2}_{L^\infty(\R^2)}\|\hat  v_{2,k}-\hat  v_{1,k}\|^\frac{1}{2}_{L^\infty(\R^2)}}=o(e^{-\frac{C\delta}{\eps_k}})\,\ \mbox{as} \,\ k\to\infty.$$
We thus obtain from (\ref{5.3:4}) and (\ref{5.3:5}) that for $p_0=\min\{p_{11}, p_{21}\}$,
\begin{equation}
\begin{split}
&\quad\displaystyle o(e^{-\frac{C\delta}{\eps_k}})\\
&=\intB \Big[\Big(V_1(x)+\frac{x\cdot \nabla V_1(x)}{2}\Big)(\hat u_{1,k}+\hat u_{2,k})\hat\xi_{1,k}\\
& \quad+\Big(V_2(x)+\frac{x\cdot \nabla V_2(x)}{2}\Big)(\hat v_{1,k}+\hat v_{2,k})\hat\xi_{2,k}\Big]dx\\
&=\eps_k^2\int_{B_{\frac{\delta}{\eps_k}}(0)} \bigg\{\Big[V_1(\eps_kx+x_{2,k})+\frac{(\eps_kx+x_{2,k})\cdot \nabla V_1(\eps_kx+x_{2,k})}{2}\Big](\bar u_{1,k}+\bar u_{2,k})\xi_{1,k}\\
& \quad+\Big[V_2(\eps_kx+x_{2,k})+\frac{(\eps_kx+x_{2,k})\cdot \nabla V_2(\eps_kx+x_{2,k})}{2}\Big](\bar v_{1,k}+\bar v_{2,k})\xi_{2,k}\bigg\}dx.
\end{split}\label{5.3:4W}
\end{equation}

Applying the dominated convergence theorem, we can deduce from  (\ref{eq1.26}),  (\ref{lem4.1:3}), (\ref{2:conexp}) and (\ref{uniq:a-3})  that
\begin{equation}\label{eq4.611}
\begin{split}
&\lim_{k\to\infty}\eps_k^{-p_{11}}\int_{B_{\frac{\delta}{\eps_k}}(0)} V_1(\eps_kx+x_{2,k})(\bar u_{1,k}+\bar u_{2,k})\xi_{1,k}dx\\
&=\lim_{k\to\infty}\int_{B_{\frac{\delta}{\eps_k}}(0)} \frac{V_1\big(\eps_k[x+(x_{2,k}-x_1)/\eps_k]+x_1\big)}{V_{11}\big(\eps_k[x+(x_{2,k}-x_1)/\eps_k]\big)}V_{11}\Big(x+\frac{x_{2,k}-x_1}{\eps_k}\Big)\\
&\quad\quad\quad\quad\quad\cdot(\bar u_{1,k}+\bar u_{2,k})\xi_{1,k}dx=2\inte V_{11}(x+y_0)u_0 \xi_{10}dx.
\end{split}
\end{equation}
Since $x\cdot\nabla V_{11}(x)=p_{11}V_{11}(x)$,  we can  deduce from (\ref{eq1.40}), (\ref{lem4.1:3}), (\ref{2:conexp}) and (\ref{uniq:a-3}) that
\begin{equation}\label{eq4.622}
\begin{split}
&\lim_{k\to\infty}\frac{\eps_k^{-p_{11}}}{2}\int_{B_{\frac{\delta}{\eps_k}}(0)} (\eps_kx+x_{2,k})\cdot \nabla V_1(\eps_kx+x_{2,k})(\bar u_{1,k}+\bar u_{2,k})\xi_{1,k}\\
&=p_{11}\inte V_{11}(x+y_0)u_0 \xi_{10}dx.
\end{split}
\end{equation}
Similar to   (\ref{eq4.611}) and (\ref{eq4.622}), we also have
\begin{equation*}
\begin{split}
\lim_{k\to\infty}\eps_k^{-p_{21}}\int_{B_{\frac{\delta}{\eps_k}}(0)} V_2(\eps_kx+x_{2,k})(\bar v_{1,k}+\bar v_{2,k})\xi_{2,k}dx=2\inte V_{21}(x+y_0)v_0 \xi_{20}dx,
\end{split}
\end{equation*}
and
\begin{equation*}
\begin{split}
&\lim_{k\to\infty}\frac{\eps_k^{-p_{21}}}{2}\int_{B_{\frac{\delta}{\eps_k}}(0)} (\eps_kx+x_{2,k})\cdot \nabla V_2(\eps_kx+x_{2,k})(\bar v_{1,k}+\bar v_{2,k})\xi_{2,k}\\
&=p_{21}\inte V_{21}(x+y_0)v_0 \xi_{10}dx.
\end{split}
\end{equation*}
%
%
%
%
%
%
%
%
%
%
%
Together with \eqref{5.3:4W}, we then conclude that
\begin{equation}\arraycolsep=1.5pt\begin{array}{lll}
 \quad \displaystyle \inte V_{11}(x+y_0)u_0 \xi_{10}+  \displaystyle \inte V_{21}(x+y_0)v_0 \xi_{20}= 0, \ \  \text{if} \ \ p_{11}=p_{21};\\[3mm]
 \    \displaystyle\inte V_{11}(x)u_0 \xi_{10}=0, \ \   \text{if}\ \ p_{11}<p_{21};\quad \displaystyle\inte V_{21}(x)v_0 \xi_{20}=0,\ \   \text{if} \ \ p_{11}>p_{21}.
\end{array}\label{5.3:4WW}
\end{equation}
We thus derive from (\ref{1:H}) and (\ref{5.3:4WW}) that if $p_{11}=p_{21}$, then
\begin{equation}\label{4.2:*}
\begin{split}
0= & 2\int_{\R^2} V_{11}(x+y_0)u_0\Big[b_0(u_0+x\cdot\nabla u_0)+\sum_{i=1}^{2}b_i\frac{\partial u_0}{\partial x_{i}}\Big]\\
& +2\int_{\R^2} V_{21}(x+y_0)v_0\Big[b_0(v_0+x\cdot\nabla v_0)+\sum_{i=1}^{2}b_i\frac{\partial v_0}{\partial x_{i}}\Big]\\
=& -(b_1,b_2)\cdot\nabla H_1(y_0)\\
&+2b_0\Big\{H_1(y_0)+\frac{1}{2}\int_{\R^2}\Big[V_{11}(x+y_0)(x\cdot\nabla u_0^2)+V_{21}(x+y_0)(x\cdot\nabla v_0^2)\Big]\Big\}\\
 =&-b_0\int_{\R^2}\Big\{u_0^2\big[x\cdot\nabla V_{11}(x+y_0)\big]+
v_0^2\big[x\cdot\nabla V_{21}(x+y_0)\big]\Big\}\\
=&-b_0\Big\{p_{11}H_1(y_0)-\int_{\R^2}\Big[u_0^2\,y_0\cdot\nabla V_{11}(x+y_0)+
v_0^2\,y_0\cdot\nabla V_{21}(x+y_0)\Big]\Big\}\\
=&-b_0p_{11}H_1(y_0)+b_0y_0\cdot\nabla H_1(y_0)=-b_0p_{11}H_1(y_0),
\end{split}
\end{equation}
due to the fact that $(x+y_0)\cdot\nabla V_{i1}(x+y_0)=p_{i1}V_{i1}(x+y_0)$ for $i=1,2$.
Since $H_1(y_0)>0$, we conclude from \eqref{4.2:*} that $b_0=0$. Similarly, if $p_{11}\not =p_{21}$, we can also derive from (\ref{5.3:4WW}) that $-b_0p_0H_1(y_0)=0$, where $p_0=\min\{p_{11},p_{21}\}$, which further implies that $b_0=0$. Therefore, we conclude that $b_0=0$.

Due to the non-degeneracy assumption (\ref{1:H}), setting  $b_0=0$ into (\ref{5.2:AA})--(\ref{5.2:AA-2}) then yields that $b_1=b_2=0$, which thus implies that  $\xi_{10}=\xi_{20}=0$.

 \vskip 0.1truein

\noindent{\em  Step 3.} $\xi_{10}=\xi_{20}=0$ cannot occur.

Finally, let $(x_k, y_k)$   satisfy  $|\xi_{1,k}(x_k)\xi_{2,k}(y_k)|=\|\xi_{1,k}\xi_{2,k}\|_{L^\infty(\R^2)}=1$. By the exponential decay (\ref{2:conexp}), applying the maximum principle to (\ref{step-1:9}) yields that $|x_k|\le C$ and $|y_k|\le C$ uniformly in $k$. We thus conclude that $\xi_{i,k}\to \xi_i\not\equiv 0$ uniformly on $\R^2$ as $k\to\infty$, where $i=1, 2$, which however contradicts to the fact that $\xi_{10}=\xi_{20}=0$ on $\R^2$. This completes the proof of Theorem \ref{Thm1.5}.
\qed
\vskip .1truein
We finally remark that if $V_i(x)\in C^2(\R^2)$ is homogeneous of degree $p_i\ge 2$ for $i=1,2$, Theorem \ref{Thm1.5} is then reduced immediately into the following simplified version.

\begin{cor}\label{prop4.4}
Suppose $\lim_{|x|\to\infty}V_i(x)=\infty$ and $V_i(x)\in C^2(\R^2)$ is homogeneous of degree $p_i\ge 2$, where $i=1,2$. Assume that
\begin{equation*}
y_0 \,\ \text{is the unique and non-degenerate critical point of}\,\ H(y),
\end{equation*}
where $H(y)$ is defined by
\begin{equation*}
H(y)=\begin{cases}
\displaystyle\gamma\displaystyle\inte V_1(x+y)w^2(x)dx\,\  &\text{if }\,\  p_1<p_2,\\[3mm]
\displaystyle\inte \big[\gamma V_1(x+y)+(1-\gamma)V_2(x+y)\big]w^2(x)dx\,\  &\text{if }\,\  p_1=p_2,\\[3mm]
(1-\gamma)\displaystyle\inte V_2(x+y)w^2dx\,\  &\text{if }\,\  p_1>p_2,
\end{cases}
\end{equation*}
and $0<\gamma <1$ is given by \eqref{def:beta.k}.
Then for any given $a_1\in (0, a^*)$ and $a_2\in (0, a^*)$, there exists a unique nonnegative minimizer of $e(a_1,a_2,\beta)$ as $\beta \nearrow \beta^*$.
\end{cor}

\appendix
\section{Appendix}
\subsection{Equivalence between ground states and constraint minimizers}  \label{ap1}

In this appendix, we shall establish the following proposition on the equivalence between ground states of \eqref{equ:CGPS} and constraint minimizers of (\ref{def:e}):

\begin{prop} \label{prop:A1} Suppose $(a_1,a_2,\beta)\in \R^+\times \R^+\times\R^+$ is given, then
any minimizer of (\ref{def:e}) is a ground state of (\ref{equ:CGPS}) for some    $\mu\in \R$; conversely, any  ground state of (\ref{equ:CGPS}) for some $\mu\in \R$ is a minimizer of (\ref{def:e}).
\end{prop}

Given any $(a_1,a_2,\beta)\in \R^+\times \R^+\times\R^+$, the energy  functional of \eqref{equ:CGPS} is defined  by
\begin{equation}\label{eqa2}
\begin{split}
I_\mu(u_1,u_2)=&\sum_{i=1}^2\int_{\R^2}\Big\{\frac{1}{2}\Big[|\nabla u_i|^2+(V_i(x)-\mu)|u_i|^2\Big]-\frac{a_i}{4}| u_i|^4\Big\}dx\\
&-\frac{\beta}{2}\int_{\R^2}|u_1|^2|u_2|^2dx\\
=&\frac{1}{2}E_{a_1,a_2,\beta}(u_1,u_2)-\frac{\mu}{2}\int_{\R^2}(|u_1|^2+|u_2|^2)dx,\,\ (u_1,u_2)\in \mathcal{X},
\end{split}
\end{equation}
where $\mu\in \R$ is a parameter and the energy functional $E_{a_1,a_2,\beta}(u_1,u_2)$ is given by (\ref{def:E}).
The set of all nontrivial weak solutions for (\ref{equ:CGPS}) is then given by
\begin{equation*}
S_\mu:=\Big\{
  (u_1,u_2)\in \mathcal{X}\setminus \{(0,0)\}:
  \ \langle I_\mu'(u_1,u_2),(\varphi_1, \varphi_2)\rangle=0,\, \forall \ (\varphi_1, \varphi_2)\in \mathcal{X}
  \Big\},
\end{equation*}
and the set of all ground states for (\ref{equ:CGPS}) is thus defined as
\begin{equation}
G_\mu:=\Big\{
  (u_1,u_2)\in S_\mu:\,
  I_\mu(u_1,u_2)\leq I_\mu(\bar u_1,\bar u_2)\, \text{ for all }\, (\bar u_1,\bar u_2) \in S_\mu
  \Big\}.
\end{equation}

\noindent {\bf Proof of Proposition \ref{prop:A1}.} For any given $(a_1,a_2,\beta)\in \R^+\times \R^+\times\R^+$, assume   $(u_{1\beta},u_{2\beta})$ is a minimizer of (\ref{def:e}), and suppose
 $(u_1,u_2) $ is a ground state of (\ref{equ:CGPS}) for some $\mu\in \R$. Set
\begin{equation}\label{eqa3}\tilde u_i=\frac{u_i}{\sqrt{\rho}},\,\ \text{ where }\, \rho:=\int_{\R^2}\big(|u_1|^2+|u_2|^2\big)dx>0, \,\ i=1,\, 2,\end{equation}
so that $\int_{\R^2}(|\tilde u_1|^2+|\tilde u_2|^2)dx =1 $. Thus,
\begin{equation*}
E_{a_1,a_2,\beta}(\tilde u_1, \tilde u_2)\geq E_{a_1,a_2,\beta}(u_{1\beta},u_{2\beta})\, \text{ and } \, I_{\mu}(u_{1\beta},u_{2\beta}) \geq I_{\mu}(u_1,u_2).
\end{equation*}
It then follows from (\ref{eqa2}) that
\begin{equation}\label{eqa4}
I_{\mu}(\tilde u_1,\tilde u_2) \geq I_{\mu}(u_{1\beta},u_{2\beta})\geq I_{\mu}(u_1,u_2).
\end{equation}
On the other hand, using (\ref{equ:CGPS}) we derive from (\ref{eqa2}) that
\begin{equation}\label{eqa4:A}
I_{\mu}(u_1,u_2)=\sum_{i=1}^2\frac{a_i}{4}\int_{\R^2}| u_i|^4dx+\frac{\beta}{2}\int_{\R^2}|u_1|^2|u_2|^2dx.
\end{equation}
Using (\ref{equ:CGPS}) and (\ref{eqa3}), we also have
\begin{equation}\label{eqa4:B}
  I_{\mu}(\tilde u_1,\tilde u_2)
=\sum_{i=1}^2\frac{a_i}{4\rho}\Big(2-\frac{1}{\rho}\Big)\int_{\R^2}| u_i|^4dx+\frac{\beta}{2\rho}\Big(2-\frac{1}{\rho}\Big)\int_{\R^2}|u_1|^2|u_2|^2dx.
\end{equation}
We thus conclude from (\ref{eqa4})--(\ref{eqa4:B}) that
\begin{equation*}
\frac{1}{\rho}\Big(2-\frac{1}{\rho}\Big)\geq1,
\end{equation*}
which holds if and only if $\rho=1$. This further implies that (\ref{eqa4}) is indeed an equality, $i.e.,$
$$I_{\mu}(u_1,u_2)=I_{\mu}(u_{1\beta},u_{2\beta}),\
E_{a_1,a_2,\beta}(u_1,u_2)=E_{a_1,a_2,\beta}(u_{1\beta},u_{2\beta}),$$
and hence
$(u_{1\beta},u_{2\beta})$ is a ground state of (\ref{equ:CGPS}) for some    $\mu\in \R$, and  $(u_1,u_2)$ is a minimizer of (\ref{def:e}). The proof is therefore complete.
\qed

\subsection{Gagliardo-Nirenberg type inequality \eqref{Ineq:GN} }

In this appendix, we improve the results obtained in \cite[Section 3]{FM} to derive the following lemma on the Gagliardo-Nirenberg type inequality \eqref{Ineq:GN}:

\begin{lem}\label{Lem:w1w2w}
The Gagliardo-Nirenberg type inequality \eqref{Ineq:GN} holds with the best constant $\frac{2}{ \|w\|_2^2} $ and is attained at $(w\sin \theta , w\cos \theta )$ for any $\theta \in [0,2\pi )$.
\end{lem}

\noindent {\bf Proof.}
Consider the minimization problem
\begin{equation}\label{def:j}
  j:=\inf\limits_{(0,0)\not=(u_1,u_2)\in H^1(\R^2)\times H^1(\R^2)}J(u_1,u_2),
\end{equation}
where $J(u_1,u_2)$ satisfies
 \begin{equation}\label{def:J}
 J(u_1,u_2):=
  \frac{\int_{\R ^2} \big(|\nabla u_1|^2+|\nabla u_2|^2\big) \dx
        \int_{\R ^2}\big(|u_1|^2+|u_2|^2\big) \dx}
       {\int_{\R ^2} \big(|u_1|^2+|u_2|^2  \big)^2 \dx}.
 \end{equation}
Taking the test function $(u_1,u_2)= (\frac{1}{\sqrt{2}}w,\frac{1}{\sqrt{2}}w)$, we derive from \eqref{ide:w} that
\begin{equation}\label{val:J.w}
  j\leq J\Big(\frac{1}{\sqrt{2}}w,\frac{1}{\sqrt{2}}w\Big)=\frac{\|w\|_2^2 }{2}.
\end{equation}
On the other hand,  we have for any $(u_1,u_2)\in H^1(\R^2)\times H^1(\R^2)$,
\begin{equation}\label{sub:J}
\begin{split}
  J(u_1,u_2)
  & =    \frac{\int_{\R ^2} \big(|\nabla u_1|^2+|\nabla u_2|^2\big) \dx
               \int_{\R ^2}\big(|u_1|^2+|u_2|^2\big) \dx}
              {\int_{\R ^2} \big(\sqrt{|u_1|^2+|u_2|^2}\big)^4 \dx} \\
  & \geq \frac{\int_{\R ^2} \big(\nabla\sqrt{|u_1|^2+|u_2|^2}\big)^2 \dx
               \int_{\R ^2}\big(\sqrt{|u_1|^2+|u_2|^2}\big)^2 \dx}
              {\int_{\R ^2} \big(\sqrt{|u_1|^2+|u_2|^2}\big)^4 \dx}\geq \frac {\|w\|_2^{2}}{2},
\end{split}
\end{equation}
where the first inequality follows from  \cite[Theorem 7.8]{LL}, and the second one is obtained by applying \eqref{ineq:GNQ}.
Combining \eqref{val:J.w} and \eqref{sub:J} then yields that the best constant $j$ satisfies $j=\frac {\|w\|_2^{2}}{2}$.

Note from \cite[Section 3]{FM} that $j$ is attained at any $(w_1,w_2)$ satisfying
 \begin{equation}\label{Ineq:GN.w1w2}
  j=\frac{ \|w_1\|_2^2+\|w_2\|_2^2}{2}= J(w_1,w_2),
 \end{equation}
where $(w_1,w_2)$ is a ground state of the following system
\begin{equation}\label{equ:CNSS}
  \begin{cases}
-\Delta u_{1}+u_{1}=u_{1}^3+ u_{2}^2u_{1}\,\,\ \mbox{in}\,\ \R^2,\\
-\Delta u_{2}+u_{2}=u_{2}^3+ u_{1}^2u_{2}\,\,\ \mbox{in}\,\ \R^2.
\end{cases}
\end{equation}
Recall that the existence of ground states for \eqref{equ:CNSS} is given in \cite{MMP}.
Following \eqref{Ineq:GN.w1w2}, we now get that
\begin{equation}\label{ide:w1w2.w}
  \|w_1\|_2^2+\|w_2\|_2^2 = \|w\|_2^2.
\end{equation}
Since all equalities in \eqref{sub:J} hold, we then follow from \cite[Theorem 7.8]{LL} that  there exists a constant $c$, independent of $x$, such that $w_1(x) = cw_2(x)$ holds $a.e.$ in $\R^2$.
Applying \eqref{ide:w1w2.w}, we therefore conclude that, up to scalings, there exists $\theta\in [0,2\pi )$ such that
\begin{equation}\label{val:w1w2.w}
  w_1(x)=w(x)\sin\theta \,\ \text{and}\,\  w_2(x)=w(x)\cos\theta  \,\ a.e.\,\ \text{in} \,\ \R^2.
\end{equation}
Since $(w_1,w_2)$ is arbitrary, this completes the proof of the lemma.
\qed

\subsection{Proof of Lemma \ref{lem4.3}}

This subsection is focussed on the proof of  Lemma \ref{lem4.3}.

\noindent\textbf{Proof of Lemma \ref{lem4.3}.}
We first note from (\ref{5.2:0}) that $(\hat\xi_{1,k}, \hat\xi_{2,k})$ satisfies
\begin{equation}\label{step-11:9}\arraycolsep=1.5pt
 \left\{\begin{array}{lll}
&  \varepsilon_{k}^2\Delta \hat\xi_{1,k}-\varepsilon_{k}^2V_1(x)\hat\xi_{1,k}+\mu_{2,k} \varepsilon_{k}^2\hat\xi_{1,k}+\displaystyle\frac{a_1}{a^*}\big(\hat  u_{2,k}^2+\hat  u_{2,k}\hat  u_{1,k}+\hat  u_{1,k}^2\big)\hat\xi_{1,k} \\[2mm]
 &\qquad +\displaystyle\frac{\beta_k }{a^*}\big[\hat   v_{1,k}^2\hat \xi_{1,k}+\hat  u_{2,k}( \hat  v_{2,k}+\hat  v_{1,k})\hat \xi_{2,k}\big] +\displaystyle  c_k\hat  u_{1,k}=0 \,\ \mbox{in}\,\  \R^2,\\[4mm]
 & \varepsilon_{k}^2 \Delta \hat \xi_{2,k}-\varepsilon_{k}^2V_2(x)\hat \xi_{2,k}+\mu_{2,k} \varepsilon_{k}^2\hat \xi_{2,k}+\displaystyle\frac{a_2}{a^*}\big(\hat  v_{2,k}^2+\hat  v_{2,k}\hat  v_{1,k}+\hat  v_{1,k}^2\big)\hat \xi_{2,k} \\[2mm]
 &\qquad +\displaystyle\frac{\beta_k }{a^*}\big[\hat   u_{1,k}^2\hat \xi_{2,k}+\hat  v_{2,k}( \hat  u_{2,k}+\hat  u_{1,k})\hat \xi_{1,k}\big] +\displaystyle  c_k\hat  v_{1,k}=0 \,\ \mbox{in}\,\  \R^2,
\end{array}\right.
\end{equation}
where the coefficient $c_{k}$ is given by \eqref{step-1:10} and bounded uniformly in $k$ in view of  (\ref{eq4.36}).
Multiplying the first equation of (\ref{step-11:9}) by $\hat \xi_{1,k}$ and integrating over $\R^2$, we then obtain that
\[\arraycolsep=1.5pt\begin{array}{lll}
&&\displaystyle \eps ^2_k\inte |\nabla \hat \xi_{1,k}|^2 -\mu_{2,k}\eps ^2 _k\inte  |\hat \xi_{1,k}|^2+\eps ^2_k\inte V_1(x)|\hat \xi_{1,k}|^2 \\[4mm]
&=&\displaystyle \frac{ a_1}{a^*}\inte \big(\hat  u_{2,k}^2+\hat  u_{2,k}\hat  u_{1,k}+\hat  u_{1,k}^2\big)|\hat \xi_{1,k}|^2 \\[4mm]
&&+\displaystyle \frac{ \beta_k}{a^*}\inte \big[\hat   v_{1,k}^2|\hat \xi_{1,k}|^2+\hat  u_{2,k}( \hat  v_{2,k}+\hat  v_{1,k})\hat \xi_{1,k}\hat \xi_{2,k}\big]+ c_k\inte \hat  u_{1,k}\hat \xi_{1,k}
\\[4mm]
&\le &\displaystyle \frac{ a_1}{a^*}\inte \big(\hat  u_{2,k}^2+\hat  u_{2,k}\hat  u_{1,k}+\hat  u_{1,k}^2\big)+\displaystyle \frac{ \beta ^*}{a^*}\inte \big[\hat   v_{1,k}^2 +\hat  u_{2,k}( \hat  v_{2,k}+\hat  v_{1,k}) \big]+\hat C\inte \hat  u_{1,k}
\\[4mm]
&= &\displaystyle \frac{ a_1\eps_k^2}{a^*}\inte \big(\bar  u_{2,k}^2+\bar u_{2,k}\bar  u_{1,k}+\bar  u_{1,k}^2\big)+\displaystyle \frac{ \beta ^*\eps_k^2}{a^*}\inte \big[\bar  v_{1,k}^2 +\bar  u_{2,k}( \bar  v_{2,k}+\bar  v_{1,k}) \big]+\hat C\eps_k^2\inte \bar u_{1,k}
\\[4mm]
&\le& C\eps ^2_k\,\ \mbox{as} \,\ k\to\infty,
\end{array}\]
since $|\hat \xi_{i,k}|$ is bounded uniformly in $k$, and $\bar u_{i,k}$ and $\bar v_{i,k}$ also decay exponentially as $|x|\to\infty$, $i=1,\, 2$. This implies from \eqref{uniq:a-6HH} that there exists a constant $C_1>0$ such that
\begin{equation}
     I:=\displaystyle \eps ^2_k\inte |\nabla \hat \xi_{1,k}|^2 +\frac{1}{2}\inte  |\hat \xi_{1,k}|^2+\eps ^2_k\inte V_1(x)|\hat \xi_{1,k}|^2<C_1\eps ^2_k \quad \text{as}\ \, k\to\infty.
\label{5.2:5}
\end{equation}
Applying Lemma 4.5 in \cite{Cao}, we then conclude that for any $x_0\in\R^2$, there exist a small constant $\delta >0$ and $C_2>0$  such that
\[
    \int_{\partial B_\delta (x_0)} \Big( \eps ^2_k |\nabla \hat \xi_{1,k}|^2+ \frac{1}{2}  |\hat \xi_{1,k}|^2+ \eps ^2_k  V_1(x)|\hat \xi_{1,k}|^2\Big)dS\le C_2I\le C_1C_2\eps ^2_k\,\ \mbox{as} \,\ k\to\infty,
\]
which therefore implies that (\ref{5.2:6}) holds for $i=1$. Similarly, one can also obtain that  (\ref{5.2:6}) holds for $i=2$, and the proof is thus complete.
\qed

\vskip 0.16truein
\noindent {\bf Acknowledgements:}  Part of this work was finished when Y. J. Guo was visiting   Pacific Institute for Mathematical Sciences (PIMS) and Department of Mathematics at the University of British Columbia from March to April in
2017. He would like to thank them for their warm hospitality.



\begin{thebibliography}{40}
 \bibitem{BC} W. Z. Bao and Y. Y. Cai, {\em Ground states of two-component Bose-Einstein condensates with an internal atomic Josephson junction}, East Asia J. Appl. Math. {\bf 1} (2011), 49--81.

 \bibitem{BW} T. Bartsch and Z. Q. Wang, {\em Existence and multiplicity results for some superlinear elliptic problems on $\R^N$}, Comm. Partial Differential Equations {\bf 20} (1995), 1725--1741.


\bibitem{Cao} D. M. Cao,  S. L. Li and P. Luo, {\em Uniqueness of positive bound states with multi-bump for nonlinear Schr\"odinger equations}, Calc. Var. Partial Differential Equations  {\bf 54} (2015),  no. 4, 4037--4063.




\bibitem{DW} E. N. Dancer and J. C. Wei, {\em Spike solutions in coupled nonlinear Schr\"{o}dinger equations with attractive interaction}, Trans. Amer. Math. Soc.  {\bf 361}  (2009),  no. 3, 1189--1208.

\bibitem{Deng}  Y. B. Deng, C. S. Lin and S. Yan, {\em On the prescribed scalar curvature problem in $\R^N$, local uniqueness and periodicity}, J. Math. Pures Appl.   {\bf 104}  (2015),  no. 6, 1013--1044.


\bibitem{EGBB} B. D. Esry, C. H. Greene,  J. P. Burke  and J. L. Bohn,  {\em Hartree-Fock theory for double condensates}, Phys. Rev. Lett. {\bf78} (1997), 3594--3597.

\bibitem{FM} L. Fanelli and E. Montefusco, {\em On the blow-up threshold for weakly coupled nonlinear Schr\"{o}dinger equations}, J. Phys. A: Math. Theory {\bf 40} (2007), 14139--14150.

\bibitem{GNN}  B. Gidas, W. M. Ni and L. Nirenberg,  {\em Symmetry of positive solutions of nonlinear elliptic equations in $\mathbb{R}^n$}, Mathematical analysis and applications  Part A, Adv. in Math. Suppl. Stud. Vol. {\bf 7}  (1981), 369--402.

\bibitem{GT} D. Gilbarg and  N. S. Trudinger, {\em Elliptic Partial Differential Equations of Second Order}, Springer, (1997).


\bibitem{Grossi} M. Grossi, {\em On the number of single-peak solutions of the nonlinear Schr\"odinger equations}, Ann. Inst H. Poincar¨¦ Anal. Non Lin¨¦aire {\bf 19}  (2002), 261--280.

\bibitem{GLW} Y. J. Guo, C. S. Lin and J. C. Wei, {\em Local uniqueness and refined spike profiles of ground states for two-dimensional attractive Bose-Einstein condensates}, SIAM J. Math. Anal. {\bf 49} (2017), 3671--3715.

\bibitem{GLWZ} Y. J. Guo, S. Li, J. C. Wei and X. Y. Zeng, {\em Ground states of two-component attractive Bose-Einstein condensates II: semi-trivial limit behavior}, arxiv.org/abs/1707.07500, submitted, (2017).

\bibitem{GS}  Y. J. Guo and R. Seiringer, {\em On the mass concentration for Bose-Einstein condensates with attractive interactions}, Lett. Math. Phys. {\bf 104} (2014), 141--156.

\bibitem{GWZZ} Y. J. Guo, Z. Q. Wang, X. Y. Zeng and H. S. Zhou,  {\em  Properties of ground states of attractive Gross-Pitaevskii equations with multi-well potentials}, Nonlinearity, accepted, (2017).

\bibitem{GZZ} Y. J. Guo, X. Y. Zeng and H. S. Zhou, {\em Energy estimates and symmetry breaking  in attractive Bose-Einstein condensates with ring-shaped potentials}, Ann. Inst. H. Poincar\'e Anal. Non Lin\'eaire {\bf 33} (2016), 809--828.

\bibitem{GZZ2} Y. J. Guo, X. Y. Zeng and H. S. Zhou, {\em Blow-up solutions for two coupled Gross-Pitaevskii equations with attractive interactions},  Discrete Contin. Dyn. Syst. A, {\bf 37 (7)} (2017), 3749--3786.




\bibitem {HMEWC} D. S. Hall, M. R. Matthews,  J. R. Ensher,  C. E. Wieman and E. A. Cornell, {\em Dynamics of component separation in a binary mixture of Bose-Einstein condensates}, Phys. Rev. Lett. {\bf 81} (1998), 1539--1542.

\bibitem {HL} Q. Han and F. H. Lin, {\em Elliptic Partial Differential Equations, Second Edition}, Courant Lect. Notes Math Vol. 1, Courant Institute of Mathematical Science/AMS, New York, (2011).

\bibitem {Kwong} M. K. Kwong, {\em Uniqueness of positive solutions of $\Delta u-u+u^p=0$ in $\R^N$}, Arch. Ration. Mech. Anal. {\bf 105} (1989), 243--266.

\bibitem {LN} Y. Li and W. M. Ni, {\em Radial symmetry of positive solutions of nonlinear elliptic equations in $\R^n$}, Comm. Partial Differ. Eqns. {\bf 18} (1993), 1043--1054.

\bibitem {LL} E. H. Lieb and M. Loss, {\em Analysis}, Graduate Studies in Math. 14, Amer. Math. Soc., Providence, RI, second edition (2001).


\bibitem {LWCMP} T. C. Lin and J. C. Wei, {\em  Ground state of N coupled nonlinear Schr\"{o}dinger equations in $\R^N$, $N\le 3$}, Comm. Math. Phys. {\bf 255} (2005), 629--653.  Erratum: Comm. Math. Phys. {\bf 277} (2008), 573--576.


\bibitem {LW} T. C. Lin and J. C. Wei, {\em Spikes in two coupled nonlinear Schr\"{o}dinger equations}, Ann. Inst. H. Poincar¡äe Anal. Non Lin¡äeaire  {\bf 22} (2005), 403--439.

\bibitem {LW2} T. C. Lin and J. C. Wei, {\em Spikes in two-component systems of nonlinear Schr¡§odinger equations with trapping potentials}, J. Diff. Eqns. {\bf  229} (2006), 538--569.

\bibitem {Lions1} P. L. Lions, {\em The concentration-compactness principle in the caclulus of variations. The locally compact case. I}, Ann. Inst H. Poincar\'{e}. Anal. Non Lin\'{e}aire {\bf 1} (1984), 109--145.


\bibitem {M} M. Maeda, {\em On the symmetry of the ground states of nonlinear Schr\"odinger equation with potential}, Adv. Nonlinear Stud. {\bf 10} (2010), 895--925.

 \bibitem{MMP} L. A. Maia,  E. Montefusco and B. Pellacci, {\em Positive solutions for a weakly coupled nonlinear Schr\"{o}dinger system}, J. Diff. Eqns. {\bf 229} (2006), 743--767.




 \bibitem{PW} A. S. Parkins and D. F. Walls, {\em The physics of trapped dilute-gas Bose-Einstein condensates}, Phys. Rep. {\bf303} (1998), 1--80.

 \bibitem{RS} M.  Reed and B. Simon, {\em Methods of Modern Mathematical Physics. IV. Analysis of Operators}, Academic Press, New York-London, (1978).

 \bibitem{Royo} J. Royo-Letelier, {\em Segregation and symmetry breaking of strongly coupled two component Bose-Einstein condensates in a harmonic trap}, Calc. Var. Partial Differential Equations {\bf 49} (2014), 103--124.

 \bibitem{S} M. Struwe, {\em Variational Methods: Applications to Nonlinear Partial Differential Equations and Hamiltonian Systems}, Ergebnisse Math. {34}, Springer, (2008).


 \bibitem{W} M. I. Weinstein, {\em Nonlinear Schr\"odinger equations and sharp interpolations estimates}, Comm. Math. Phys. {\bf 87} (1983), 567--576.


 \end{thebibliography}
\end{document}